%% file: main_v5_revisions_clean.tex
\documentclass[a4paper,11pt]{article}
\usepackage{tmlrmod}

\usepackage{tocloft}  

\usepackage[nottoc]{tocbibind}  

\input{macros_v2.tex}

\newcommand{\revision}[1]{#1}
\newcommand{\minorrev}[1]{#1}

\newcommand{\renaming}[1]{#1}

\newcommand{\Lgrad}{\minorrev{L_{1}}}
\newcommand{\Lhess}{\minorrev{L_{2}}}

\title{Introduction to Model-Based Derivative-Free Optimization}

\author{\name Lindon Roberts \email lindon.roberts@unimelb.edu.au \\
      \addr School of Mathematics and Statistics \& ARC Training Centre in Optimisation Technologies, Integrated Methodologies, and Applications (OPTIMA) \\
      University of Melbourne
      }
      
\begin{document}
\thispagestyle{plain}
\maketitle

\begin{abstract}
    The field of derivative-free optimization (DFO) studies algorithms for nonlinear optimization that do not rely on the availability of gradient or Hessian information.
    It is primarily designed for settings when functions are black-box, expensive to evaluate and/or noisy.
    A widely used and studied class of DFO methods for local optimization is \renaming{model-based DFO}, where the general principles from derivative-based nonlinear optimization algorithms are followed, but local Taylor-type approximations are replaced with alternative local models constructed by interpolation.
    This document provides an overview of the basic algorithms and analysis for \renaming{model-based DFO}, covering worst-case complexity, approximation theory for polynomial interpolation models, and extensions to constrained and noisy problems.

\end{abstract}

\setcounter{tocdepth}{1}
{\small \tableofcontents }

\section{Introduction}

This article is an introduction to the algorithmic ideas and analysis for \renaming{model-based derivative-free optimization (MBDFO)}.
As the name suggests, \renaming{derivative-free optimization (DFO)---model-based or otherwise---}refers to optimization in the absence of derivative information, whether that be for the objective function and/or constraint functions.\footnote{DFO is sometimes called \emph{zeroth order optimization} (ZOO), particularly in machine learning.}
Our focus here is nonlinear optimization, where we aim to minimize a (usually \minorrev{assumed} smooth) function of several real variables, possibly with constraints, which may be specified either explicitly (e.g.~$x_1 \geq 0$) or via other nonlinear function(s) (e.g.~$c(\bx) \geq 0$). 
Most \minorrev{general-purpose} nonlinear optimization algorithms, such as gradient descent or (quasi-)Newton methods, require some derivative information, such as the gradient or Hessian of the objective (and any nonlinear constraints).
DFO methods become relevant in situations where no derivative information is available, and \renaming{MBDFO} methods form an important sub-class of DFO methods.

\subsection{Overview of Derivative-Free Optimization} \label{sec_dfo_overview}

In order to understand the uses/benefits of DFO, we first have to consider how we may obtain derivative information.
If we have a function $f:\R^n\to\R$, such as our objective function, there are three main ways to evaluate or approximate its derivatives \cite[Chapter 8]{Nocedal2006}:
\begin{description}
    \item[Explicit calculation] If the mathematical form of $f$ is known, we can directly compute the relevant derivatives by hand (or using a symbolic computation package);
    \item[Automatic (aka algorithmic) differentiation] If we have computer code to compute $\bx\mapsto f(\bx)$, automatic differentiation software can be used to create code that computes derivatives of $f$ by analyzing the code for $f$ and repeatedly applying the chain rule.\footnote{There are two `modes' of automatic differentiation, forward and reverse. If $f$ has many inputs and few outputs, which is typically the case for optimization, the reverse mode is usually more efficient. In machine learning, the reverse mode is often called \emph{backpropagation}.} This produces exact derivatives, \revision{typically for the cost (in both time and memory) of a small---that is, $\bigO(1)$ and independent of $n$---number of evaluations of $f$ \cite[Chapter 4]{Griewank2008}.} 
    \item[Finite differencing] We can approximate derivatives of $f$ by comparing values of $f$ at nearby points. For example, forward finite differencing estimates first derivatives via
    \begin{align}
        \frac{\partial f}{\partial x_i}(\bx) \approx \frac{f(\bx+h\be_i) - f(\bx)}{h}, \label{eq_fin_diff_example}
    \end{align}
    for some small $h>0$, where $\be_i\in\R^n$ is the $i$-th coordinate vector.\footnote{\revision{If $f$ can be extended to complex-valued inputs, the \emph{complex step} approximation $\frac{\partial f}{\partial x_i}(\bx) \approx \operatorname{Im}(f(\bx+h\be_i))/h$ can be used \cite{Squire1998}. This is less susceptible than \eqref{eq_fin_diff_example} to rounding errors when $h$ is very small.}}
    To evaluate a full gradient $\grad f(\bx)$, we would need to evaluate $f$ at $\bx$ and $\bx+h\be_i$ for all $i=1,\ldots,n$ (i.e.~$n+1$ evaluations of $f$).
    This only requires the ability to evaluate $f$, but is only an approximation. 
\end{description}
In most circumstances, at least one of these three approaches can be successfully used without significant effort.
DFO methods are useful in the situations where none of these approaches are possible or practical.
This typically is in situations where at least two of the following apply:
\begin{description}
    \item[\minorrev{Black-box}] That is, the underlying mathematical structure of $f$ is not available. This could be because the details are unknown (e.g.~legacy or proprietary software is used) or a situation where a clear mathematical description does not exist (e.g.~results of a real-world experiment). This immediately rules out both explicit calculation and automatic differentiation.
    \item[Expensive to evaluate] If computing $\bx\mapsto f(\bx)$ is costly to evaluate---whether that cost is time, effort or money---then any algorithms relying on many evaluations of $f$ are likely to be impractical. Finite differencing  requires at least $n+1$ evaluations of $f$ at each iteration, and so may not be suitable in this case.
    \item[Noisy] In essence, this means that evaluating $f$ at nearby inputs leads to non-negligible differences in outputs. This noise may be deterministic or stochastic, depending on whether evaluating $f$ repeatedly at the same point gives the same result or not. This would usually not refer to rounding errors, but more significant differences such as the results of a Monte Carlo simulation. Here, finite differencing cannot be used without significant care, as it relies on comparing function values at nearby points.
\end{description}
A large survey of practical examples where these situations arise and DFO methods are useful is \cite{Alarie2021}.
However, some specific examples include the following.

\begin{example}[Calibrating Climate Models \cite{Tett2022}]
In many areas of science, mathematical models of real-world phenomena are constructed to help predict outcomes and make decisions.
Often, these models have parameters (e.g.~coefficients of different terms in a differential equation) that cannot be directly measured, and instead must be calibrated using indirect measurements.
This most commonly leads to least-squares minimization problems
\begin{align}
    \min_{\bx\in\R^n} f(\bx) := \sum_{i=1}^{N} (\operatorname{model}(\bx,\by_i) - z_i)^2,
\end{align}
where $\bx$ are the model parameters, $z_1,\ldots,z_N\in\R$ are observations, and $\operatorname{model}(\bx,\by_i)$ represents evaluating the model with parameters $\bx$ and other inputs $\by_i$ to produce a predicted value for $z_i$.
In atmospheric physics, $\operatorname{model}(\bx,\by)$ can involve a global climate simulation, incorporating coupled ocean, atmospheric and sea ice dynamics, and so its mathematical description is too complex to allow for a clear description (i.e.~$f$ is black-box).
Because of this complexity, evaluating $\operatorname{model}(\bx,\by)$ can be computationally expensive, \revision{e.g.~8 hours on a high-performance computing system to simulate 18 months of climate forecasts}. 
Additionally, climate dynamics can exhibit chaotic behavior, which in practice means that very similar parameters $\bx$ can produce very different model results, which effectively mean that $\operatorname{model}(\bx,\by)$ is noisy.
\end{example}

\begin{example}[Quantum Optimization \cite{Abbas2023}]
One promising use of quantum computers is to solve combinatorial optimization problems such as the graph max-cut problem.
In the Quantum Approximate Optimization Algorithm, a particular combinatorial problem is converted into a minimum eigenvalue problem
\begin{align}
    \min_{\bx\in\R^n} f(\bx) := \bm{\psi}(\bx)^* \bH \bm{\psi}(\bx),
\end{align}
for some (Hermitian) matrix $\bH$ and complex vector $\bm{\psi}(\bx)$.
However, $\bm{\psi}$ is described implicitly, and $f(\bx)$ must be evaluated by performing specific calculations on a quantum computer.
However, quantum computers are inherently stochastic, and so we can only ever evaluate random approximations to $f(\bx)$.
Given the limited availability of quantum computing hardware, the cost of evaluating $f(\bx)$ may also be high.
\end{example}

\paragraph{Types of DFO Methods}
There are many different classes of DFO methods, both for local and global optimization.
DFO methods are very common in global optimization (since converging to a stationary point is not sufficient), and the resources \cite{Audet2017,Locatelli2013,Rios2013} provide an overview of different global optimization methods, such as Bayesian optimization, genetic algorithms, branch-and-bound, and Lipschitz optimization \revision{(such as DIRECT \cite{Jones1993})}.
Our focus here is on local optimization, where we do aim to find (approximate) stationary points, just like many popular nonlinear optimization methods such as (quasi-)Newton methods.
\revision{Indeed, global optimization based only on local problem information (such as objective values and/or derivatives) cannot succeed unless an algorithm densely samples the feasible region or some global information is used, such as convexity or global derivative bounds \cite{Stephens1998}.}
The most common DFO methods for local optimization are:
\begin{description}
    \item[\renaming{Model-based DFO (MBDFO)}] The focus of this introduction. Here, the goal is to mimic derivative-based optimization methods, substituting local gradient-based approximations such as Taylor series with local models constructed by interpolation. 
    \item[Direct search methods] Such methods attempt to iteratively improve on a candidate solution by sampling nearby points, but without using any gradient approximations. This definition is broad---see a discussion of this in \cite[Section 1.4]{Kolda2003}---but typically nearby points are selected from a small number of perturbations of the current iterate. A famous example is the Nelder--Mead method, where a simplex of points is iteratively modified to find improved solutions, but more widely studied are Generating Set Search \cite{Kolda2003} \revision{and Mesh Adaptive Direct Search \cite{Audet2006a}} methods, where the perturbations are chosen in an predictable, structured way.
    \item[Finite differencing/implicit filtering] Although finite differencing-based derivative approximations are typically used within standard (derivative-based) optimization algorithms without modification, some works explicitly consider the management of the perturbation size $h$ in \eqref{eq_fin_diff_example} the algorithm. If $h$ is large, this is sometimes called implicit filtering. Building gradient approximations by finite differencing-type approximations along randomly generated directions \cite{Nesterov2017} is particularly popular in machine learning \cite{Ghadimi2013}. 
\end{description}
For readers interested in other classes of local DFO, we refer to the books \cite{Conn2009,Kelley2011,Audet2017} and survey papers \minorrev{\cite{Powell1998,Kolda2003,Custodio2012,Audet2014,Larson2019,Dzahini2025}}.
There are different reasons to prefer these classes---for example, direct search methods are much more developed than \renaming{MBDFO} at handling complex problem structures such as nonsmoothness and discrete variables.
However, the similarity of \renaming{MBDFO} to widely recognized and successful (derivative-based) optimization algorithms is appealing, and they often perform well in practice, \revision{when measured by the total number of objective evaluations required to solve a problem} (see below). 

\revision{In direct search methods, the \emph{search--poll} paradigm \cite{Booker1999} allows for very flexible `search steps', which allow for any procedure that produces potentially good points, to be combined with the rigorous `poll step' from direct search methods.
Search steps using \renaming{MBDFO} ideas is one popular approach in this framework which significantly improves the practical performance of direct search methods \cite{Custodio2010,Conn2013}.
Global surrogate models for the objective such as radial basis functions and Gaussian Processes can also enable a search step \cite{Booker1999,Audet2018,Audet2022a}. 
Here, a surrogate does not need to accurately approximate the objective, only accurately rank the quality of suggested points \cite{Audet2017,Lakhmiri2022}.}

\paragraph{Cheap vs.~Expensive Evaluations}
An extremely important issue in DFO is whether a given problem has functions (objective, constraints, etc.) that are `cheap' or `expensive' to evaluate.
These terms are meant in a (somewhat) relative sense, as different people/contexts have different understandings of cost.\footnote{As above, `cost' could be financial, computational effort or time, for example.} 
\revision{This can vary from seconds (e.g.~valve train design \cite{Choi2000} and financial model calibration \cite{NAG2019}) to hours (e.g.~algorithm tuning \cite{Audet2006}) to days (e.g.~hydrofoil noise reduction \cite{Marsden2007}).}

In the case of expensive evaluations, we typically mean that the cost of evaluating the objective (or other such function) is the dominant cost of the optimization process.
If this is true, the effort required to determine the next candidate point to evaluate is essentially negligible, and so little attention need be given to topics such as efficient linear algebra or subproblem implementations, or memory management.
Performing a significant amount of work within the algorithm in order to \revision{avoid even one extra} evaluation is generally considered worthwhile.
\revision{In this setting, some MBDFO software will store the full history of objective evaluations, which implicitly assumes the number of evaluations is not too large (e.g.~the IBCDFO software tries to construct its models primarily by using already-evaluated points; see \secref{sec_interp_practical_ibcdfo}).}

When evaluations are cheap, we need to pay more attention to these other factors impacting the performance of the optimization.
However, as in many computational tasks, the level of effort that should be devoted to more efficient algorithm implementations depends on the context; whether solving an optimization problem takes 1 second or 10 seconds on a laptop is often not a critical distinction.

As described above, DFO is particularly important/useful in the expensive setting, but may also be used in the cheap setting if the situation requires it.
In line with the expensive regime, numerical comparisons of DFO algorithms/implementations in the research literature are most commonly performed by comparing the number of evaluations required to solve a given problem, rather than other metrics such as number of iterations or runtime \cite{More2009}.
\revision{General advice on comparing optimization algorithms can be found in \cite{Beiranvand2017}.}

\subsection{Model-Based DFO}
In \renaming{MBDFO}, our broad goal is to use the algorithmic structures found in derivative-based algorithms, but replace local, Taylor series-based approximations for the objective (and/or constraints) with local models constructed by interpolation.
Hence, a \renaming{MBDFO} algorithm maintains a small set of points in $\R^n$ (of which one is usually the current iterate), and iteratively moves this set towards a solution.
Although in principle this is a generic framework, in practice the underlying derivative-based algorithms used are \emph{trust-region methods} \cite{Conn2000}.
Here, we control the maximum stepsize we are willing to take in any given iteration (balancing accurate approximations with fast convergence), and minimize Taylor-like models in a  ball around the current iterate.
Trust-region methods are a widely used class of methods for derivative-based optimization, alongside other approaches such as linesearches and regularization methods.
In \renaming{MBDFO}, trust-region methods are much more widely used than these alternatives, because we always know in advance a region where our next iterate will be found, and it is much more natural to construct interpolation models when we know the exact region in which they need to be a good approximation.

\begin{figure}[tb]
  \centering
  \begin{subfigure}[b]{0.48\textwidth}
    \includegraphics[width=\textwidth]{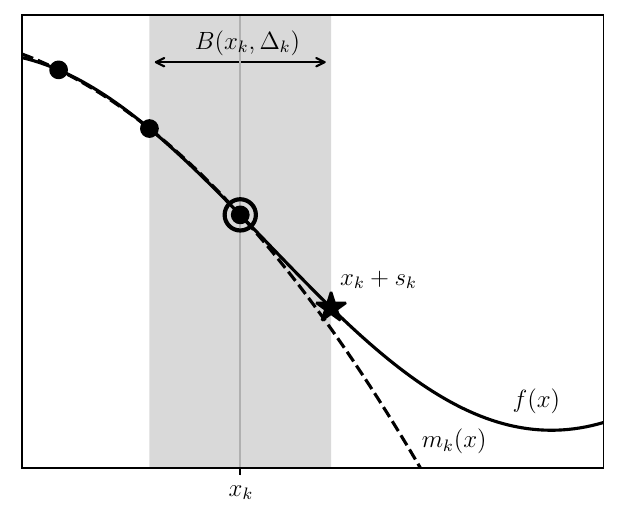}
    \caption{Initial interpolation model and step calculation}
    \label{fig_dfotr_demo_init}
  \end{subfigure}
  ~
  \begin{subfigure}[b]{0.48\textwidth}
    \includegraphics[width=\textwidth]{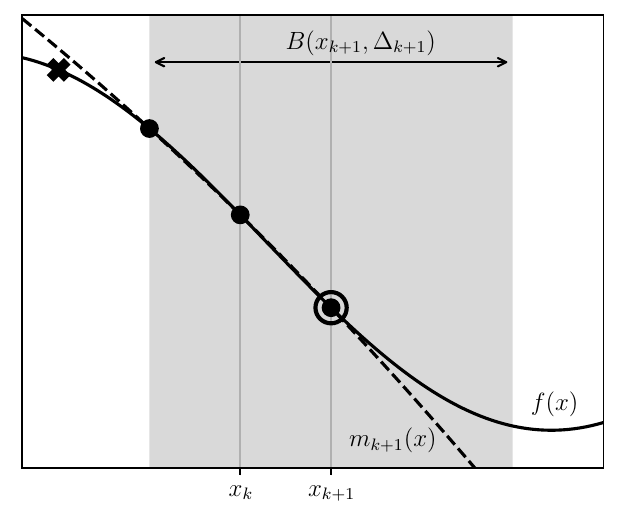}
    \caption{After step acceptance, $x_{k+1}=x_k+s_k$}
    \label{fig_dfotr_demo_after}
  \end{subfigure}
  \caption{Single iteration of \renaming{MBDFO}. \revision{The left-most point is used to build the model in (a), despite being outside the trust region, but is not used to build the model in (b) (indicated with a cross).}}
  \label{fig_dfotr_demo}
\end{figure}

A simple one-dimensional illustration of \renaming{MBDFO} is given in \figref{fig_dfotr_demo}. 
At the start of iteration $k$, we have three interpolation points\revision{---$x_k$ (illustrated with a large circle) and two others (small circles)---}which we use to construct a quadratic approximation $m_k(x)$ (dashed line) to approximate the true objective $f(x)$ (solid line).
\revision{In practice, of course, the full function $f(x)$ is not fully known to the algorithm and can only be sampled.}
We minimize $m_k(x)$ inside the (shaded) trust region, $B(x_k,\Delta_k)$ for some $\Delta_k>0$, to get a tentative new iterate $x_k+s_k$ (illustrated with a star).
We compute $f(x_k+s_k)$, and in this case observe that $f(x_k+s_k) < f(x_k)$, and so we set $x_{k+1}=x_k+s_k$ as the next iterate.
After accepting the step (and setting $\Delta_{k+1}=2\Delta_k$ to help speed up convergence), in \figref{fig_dfotr_demo_after}, we have added $x_{k+1}$ as an interpolation point---and removed the left-most point from the interpolation set (illustrated with a cross)---to build a new model $m_{k+1}$, which will subsequently be minimized inside $B(x_{k+1},\Delta_{k+1})$.

Of course, in this article we will outline exactly how we select/update interpolation points, construct interpolation models, compute steps, and choose $x_{k+1}$ and $\Delta_{k+1}$.
In particular, in higher dimensions we will need to ensure that the interpolation points are well-spaced, while remaining close to the current iterate, to ensure that a suitable model exists and is a sufficiently good approximation to the objective.
Much of our theoretical discussion will be on interpolation set management.

\subsection{Scope of this work}
The goal of this work is to provide an accessible introduction to \renaming{MBDFO}.
The intended audience is graduate students with an interest in nonlinear optimization, and researchers hoping to learn the fundamentals of this topic.
Although the end of each section includes a discussion of important references and related works, it is not intended to be a comprehensive survey (see \cite{Larson2019} for this), but rather to provide an overview of the fundamental algorithms, approximation techniques, and convergence theory used in the field.
It aims to provide a deeper introduction than found in the books \cite{Nocedal2006,Audet2017}.
The most similar existing work to this article is the book \cite{Conn2009}, but we include more recent advances (e.g.~complexity analysis of algorithms, approximation theory in constrained regions, stochastic algorithms), and provide a more targeted analysis of core concepts (for example, we do not cover interpolation theory for higher order polynomial models), including some novel proofs of existing results.
Some familiarity with nonlinear optimization techniques (e.g.~from \cite{Nocedal2006}), especially trust-region methods, would be helpful here but is not essential.

The largest part of our efforts is devoted to:
\begin{itemize}
    \item Introducing the basic (trust region) \renaming{MBDFO} algorithm and proving first- and second-order convergence and worst-case complexity bounds; and
    \item Showing how linear and quadratic interpolation models can be constructed for use in \renaming{MBDFO} methods, and proving results quantifying their accuracy. \revision{We consider structured interpolation set choices yielding optimal algorithm bounds and their practical extension where points are selected from a database of existing evaluations (\secref{sec_model_construction}) and geometric approaches based on incremental/minimal updating of interpolation points (\secref{sec_self_correcting}). These two approaches are both theoretically interesting and form the basis of practical implementations.}
\end{itemize}
For both of these topics, we primarily consider the simplest case of unconstrained optimization problems with access to exact objective values.
These two topics form the core from which more advanced algorithms can be built.
For example, in later sections, we show how these ideas can be extended to handle constrained problems, and problems with inexact/noisy objective evaluations.
At the end of each section, we provide a short overview of key references and related works, including more recent research directions.

Our focus here is on the key theoretical underpinnings of \renaming{MBDFO}.
We largely avoid detailed discussions of practical considerations for software implementations, providing only occasional notes on this important topic (most notably on termination conditions, \secref{sec_termination}) and a list of some notable open-source \renaming{MBDFO} packages in \secref{sec_conclusion}.
In particular, we give minimal attention to the numerical linear algebra required to solve the subproblems that arise in \renaming{MBDFO} (e.g.~step calculation and interpolation model building), in line with the regime of expensive evaluations described in \secref{sec_dfo_overview}. 
All the methods describe here, except for \algref{alg_storm} in \secref{sec_stochastic_noise}, are suitable for the expensive evaluation regime (see \remref{rem_storm_expensive} for a discussion of this point).

For brevity, we also avoid discussion of and comparison with other local DFO methods or global optimization methods, which often avoid using derivative information (see \secref{sec_dfo_overview} and references therein for good resources on these methods).
Although there are a wide range of techniques for global optimization, we note that Bayesian and surrogate optimization \minorrev{\cite{Shahriari2016,Booker1999,Regis2007,Audet2017}} have similar principles to \renaming{MBDFO}, namely iteratively minimizing an easy-to-evaluate approximation to the objective built from evaluations at known points.

\paragraph{Structure}
In \secref{sec_technicalities} we give some \revision{preliminary technical results that we will use throughout, and introduces standard (i.e.~derivative-based) trust-region methods for unconstrained optimization, including the general algorithmic framework, step calculation and convergence guarantees.}
In \secref{sec_dfotr}, we show how this framework extends to the \renaming{MBDFO} case, outlining what changes to the algorithmic framework are required, the requirements on the interpolation models, and the convergence guarantees.
We then introduce the polynomial interpolation theory required to construct the models, and show how they satisfy the requirements from \secref{sec_dfotr}. 
\revision{In particular, \secref{sec_model_construction} considers structured interpolation set construction and \secref{sec_self_correcting} extends these ideas to allow for incremental interpolation set updating.}
We then consider two important extensions of this theory: \secref{sec_constraints} covers problems with simple \minorrev{(e.g.~bounds)} and general nonlinear constraints, and \secref{sec_noisy} covers problems with deterministic and stochastic noise.
In our summary, \secref{sec_conclusion}, we point to a selection of state-of-the-art software implementations of \renaming{MBDFO}.

\paragraph{Notation}
Throughout, vectors and matrices will be written in bold (to distinguish from scalars), and functions will be written in bold if they have multiple outputs, e.g.~$f:\R^n\to\R$ but $\br:\R^n\to\R^m$.
The vector norm $\|\bx\|$ refers to the Euclidean 2-norm, and the matrix norm $\|\bA\|$ refers to the operator 2-norm (i.e.~largest singular value).
Other norms will be indicated explicitly, such as $\|\bx\|_1$ or $\|\bA\|_{\infty}$.
We will use $B(\bx,r)$ to refer to the closed Euclidean ball centered at $\bx\in\R^n$ with radius $r>0$, that is $B(\bx,r) := \{\by\in\R^n : \|\by-\bx\| \leq r\}$.
The standard coordinate vectors in $\R^n$ will be written as $\be_1,\ldots,\be_n$, where $\be_i$ has a 1 in the $i$-th entry, $\be$ is the vector of all ones, $\bI$ will be the identity matrix, and $\bm{0}$ will be a vector or matrix of zeros (depending on the context).
For a matrix $\bA\in\R^{m\times n}$, we denote its Moore-Penrose pseudoinverse by $\bA^{\dagger}\in\R^{n\times m}$.
\revision{We use $\kappa(\bA)$ to denote the (2-norm) condition number of a matrix, $\kappa(\bA) := \|\bA\| \: \|\bA^{\dagger}\|$.}

\section{Technical Preliminaries} \label{sec_technicalities}

\subsection{Stationarity Conditions and Smoothness Assumptions}

The largest focus of this work is to solve the unconstrained nonlinear optimization problem
\begin{align}
    \min_{\bx\in\R^n} \: f(\bx), \label{eq_problem}
\end{align}
where $f:\R^n\to\R$ is a smooth (we will mostly assume continuously differentiable) but possibly nonconvex (i.e.~not necessarily convex\footnote{The function $f$ is convex if $f(t\bx+(1-t)\by) \leq tf(\bx) + (1-t)f(\by)$ for all $\bx,\by\in\R^n$ and $t\in[0,1]$. An important consequence of convexity is that all local minimizers are global minimizers. If $f$ is known to be convex, \eqref{eq_problem} is typically much easier to solve and there are a plethora of specialized methods for this case \cite{Nesterov2004}.}) objective function.
In general, finding a global minimizer of a nonconvex function is very difficult and suffers from the \emph{curse of dimensionality}\footnote{That is, the cost of finding a solution grows exponentially with the problem dimension $n$.} \cite[Chapter 1.2.6]{Cartis2022}, so we will restrict our consideration to finding local minima for \eqref{eq_problem}.

\begin{definition}
    A point $\bx^*\in\R^n$ is a local minimizer of $f:\R^n\to\R$ if there exists $\epsilon>0$ such that $f(\bx^*) \leq f(\bx)$ for all $\bx\in B(\bx^*,\epsilon)$.
    The value $f(\bx^*)$ is called the local minimum of $f$.
\end{definition}

For general problems such as \eqref{eq_problem}, we need to specify what information about $f$ is known to any algorithm we develop, as this clearly constrains our available algorithmic choices.
In nonlinear optimization, we typically use the \emph{oracle model} for optimization algorithms: information about $f$ is available only via a given set of \emph{oracles}, most commonly some/all of:
\begin{enumerate}[label=(\alph*)]
    \item The zeroth-order oracle for $f$ is the map $\bx\mapsto f(\bx)$;
    \item The first-order oracle for $f$ is the map $\bx\mapsto \grad f(\bx)$ \revision{(provided $f$ is differentiable)}; and
    \item The second-order oracle for $f$ is the map $\bx\mapsto \grad^2 f(\bx)$ \revision{(provided $f$ is twice differentiable)}.
\end{enumerate}
When we first explore classical (derivative-based) trust-region methods in \secref{sec_tr}, we will assume access to at least the zeroth- and first-order oracles for $f$, with the second-order oracle often considered optional.
When in \secref{sec_dfotr} we begin to consider \renaming{MBDFO} algorithms, we will only assume access to a zeroth-order oracle.

\begin{remark}
    \revision{For the purposes of the theoretical analysis, we will assume $f$ is differentiable.
    However, all our algorithms are built assuming we cannot evaluate $\grad f$, for one (or more) of the reasons outlined above.}
\end{remark}

The limited problem information available in the oracle model means we need necessary or sufficient conditions for a point to be a local minimizer which can be checked using the oracles available to us.

\begin{proposition}[Theorems 2.2--2.4, \cite{Nocedal2006}] \label{prop_optimality_conditions}
    Suppose $\bx^*\in\R^n$ and $f:\R^n\to\R$ is continuously differentiable.
    Then,
    \begin{enumerate}[label=(\alph*)]
        \item If $\bx^*$ is a local minimizer of $f$, then $\grad f(\bx^*) = \bm{0}$. (first-order necessary condition); \label{first_order_necessary}
        \item If $f$ is twice continuously differentiable and $\bx^*$ is a local minimizer of $f$, then $\grad^2 f(\bx^*)$ is positive semidefinite. (second-order necessary condition); and \label{second_order_necessary}
        \item If $f$ is twice continuously differentiable, $\grad f(\bx^*)=\bm{0}$ and $\grad^2 f(\bx^*)$ is positive definite, then $\bx^*$ is a local minimizer of $f$. (second-order sufficient conditions). \label{second_order_sufficient}
    \end{enumerate}
\end{proposition}

In light of \propref{prop_optimality_conditions}\ref{first_order_necessary} and \ref{second_order_necessary}, we call points $\bx$ such that $\grad f(\bx)=\bm{0}$ \emph{stationary} or \emph{first-order optimal}, and if both $\grad f(\bx)=\bm{0}$ and $\grad^2 f(\bx)$ is positive semidefinite then we say $\bx$ is \emph{second-order optimal}.

\revision{In practice, our algorithms will be able to find points that approximately satisfy the first- or second-order necessary conditions. In the first-order case, this means finding a point with $\|\grad f(\bx)\| \leq \epsilon$ for some (small) tolerance $\epsilon>0$. Since a second-order optimal point is defined by two necessary conditions $\|\grad f(\bx_k)\|=0$ and $\lambda_{\min}(\grad^2 f(\bx_k)) \geq 0$, our measure of approximate second-order optimality will be $\sigma(\bx) \leq \epsilon$, where
\begin{align}
    \sigma(\bx) := \max(\|\grad f(\bx)\|, \tau(\bx)), \qquad \text{with} \qquad \tau(\bx) := \max(-\lambda_{\min}(\grad^2 f(\bx)), 0). \label{eq_optimality_2}
\end{align}
That is, we hope to find a point with $\sigma(\bx) \leq \epsilon$, which implies both $\|\grad f(\bx)\| \leq \epsilon$ and $\tau(\bx) \leq \epsilon$ hold, the latter condition being equivalent to $\lambda_{\min}(\grad^2 f(\bx)) \geq -\epsilon$.
For an algorithm with current iterate $\bx_k$ at iteration $k$, we will denote $\sigma_k := \sigma(\bx_k)$ and $\tau_k := \tau(\bx_k)$.}

We will consider two different assumptions on our objective $f$ in \eqref{eq_problem}, depending on whether we wish to prove convergence to first- or second-order optimal points.
Our main discussion will focus on first-order convergence, in which case we will assume the following.

\begin{assumption} \label{ass_smoothness_1}
    \revision{The function $f:\R^n\to\R$ in \eqref{eq_problem} satisfies:
    \begin{enumerate}[label=(\alph*)]
        \item $f$ is continuously differentiable, \minorrev{and} $\grad f$ is Lipschitz continuous with constant $\Lgrad$; and \label{ass_smoothness_1_smooth}
        \item $f$ \revision{is bounded below, there exists an $\flow\in\R$ such that} $f(\bx) \geq \flow$ for all $\bx\in\R^n$. \label{ass_smoothness_1_bdd_below}
    \end{enumerate}}
\end{assumption}

\begin{remark}
    \revision{The global Lipschitz continuity of $\grad f$ in \assref{ass_smoothness_1}\ref{ass_smoothness_1_smooth} is very strong, excluding simple functions such as $f(x)=x^4$.
    As in \cite{Cartis2022}, we make this assumption for ease of exposition, but in practice it can be weakened to $\grad f$ being $\Lgrad$-Lipschitz continuous on any set $\mathcal{L}$ containing all trust regions $\cup_{k} B(\bx_k,\Delta_k) \subseteq\mathcal{L}$. For example, if $f(\bx_k) \leq f(\bx_0)$ for all iterates (i.e.~our algorithm is \emph{monotone}) and we have an upper bound on all trust-region radii, $\Delta_k \leq \Delta_{\max}$, then we can take $\mathcal{L} = \cup_{\bx : f(\bx) \leq f(\bx_0)} B(\bx,\Delta_{\max})$. For problems with simple constraints (\secref{sec_convex_constraints}), we can restrict $\mathcal{L}$ to the feasible region. Moreover, if $\grad f$ is continuous and $\mathcal{L}$ is bounded, then the existence of a suitable $\Lgrad$ is automatic. This issue is discussed more in, for example, \cite[Chapter 6.2.1]{Conn2000} and \cite[Chapter 10.2]{Conn2009} for deterministic problems, and \cite[Remark 4.2]{Chen2018} for stochastic problems.}
\end{remark}

The most important consequence of \assref{ass_smoothness_1} is a bound on the error in a first-order Taylor series for $f$, specifically:

\begin{lemma}[Theorem A.8.1, \cite{Cartis2022}] \label{lem_lsmooth}
    Suppose \minorrev{\assref{ass_smoothness_1}\ref{ass_smoothness_1_smooth}} holds.
    Then 
    \begin{align}
        \left|f(\by) - f(\bx) - \grad f(\bx)^T (\by-\bx)\right| \leq \frac{\Lgrad}{2} \|\by-\bx\|^2, \label{eq_lsmooth_bound}
    \end{align}
    for all $\bx,\by\in\R^n$.
\end{lemma}

In the case where we are interested in convergence to second-order optimal points, we will use the following alternative smoothness assumption.

\begin{assumption} \label{ass_smoothness_2}
    \revision{The function $f:\R^n\to\R$ in \eqref{eq_problem} satisfies:
    \begin{enumerate}[label=(\alph*)]
        \item $f$ is twice continuously differentiable, \minorrev{and} $\grad^2 f$ is Lipschitz continuous with constant $\Lhess$; and \label{ass_smoothness_2_smooth}
        \item $f$ \revision{is bounded below, there exists an $\flow\in\R$ such that} $f(\bx) \geq \flow$ for all $\bx\in\R^n$. \label{ass_smoothness_2_bdd_below}
    \end{enumerate}}
\end{assumption}

Where \minorrev{\assref{ass_smoothness_1}\ref{ass_smoothness_1_smooth}} allows us to bound the first-order Taylor series error using \lemref{lem_lsmooth}, \minorrev{\assref{ass_smoothness_2}\ref{ass_smoothness_2_smooth}} allows us to bound the second-order Taylor series error via the bound \cite[Corollary A.8.4]{Cartis2022}
\begin{align}
    \left|f(\by) - f(\bx) - \grad f(\bx)^T (\by-\bx) - \frac{1}{2}(\by-\bx)^T \grad^2 f(\bx) (\by-\bx)\right| \leq \frac{\Lhess}{6} \|\by-\bx\|^3, \quad \forall \bx,\by\in\R^n. \label{eq_lh_smooth}
\end{align}
Another consequence of \minorrev{\assref{ass_smoothness_2}\ref{ass_smoothness_2_smooth}} is that $\grad f$ itself has a Lipchitz continuous gradient, and so a result similar to \lemref{lem_lsmooth} applies, namely \cite[Appendix A]{Nocedal2006}
\begin{align}
    \|\grad f(\by) - \grad f(\bx) - \grad^2 f(\bx)(\by-\bx)\| \leq \frac{\Lhess}{2} \|\by-\bx\|^2, \qquad \forall \bx,\by\in\R^n. \label{eq_lh_smooth_gradf}
\end{align}

\revision{Although the main focus of this work is the solution of unconstrained problems \eqref{eq_problem}, we will also consider problems with nonlinear constraints, of the form} 
\begin{subequations} \label{eq_generic_cons}
\begin{align}
    \min_{\bx\in\R^n} &\: f(\bx), \\
    \text{s.t.} &\: c_i(\bx) = 0, \qquad \forall i\in\mathcal{E}, \\
    &\: c_i(\bx) \leq 0, \qquad \forall i\in\mathcal{I},
\end{align}
\end{subequations}
The first-order optimality conditions for \eqref{eq_generic_cons} are the Karush--Kuhn--Tucker (KKT) conditions: defining the Lagrangian as
\begin{align}
    L(\bx,\blambda) := f(\bx) + \sum_{i\in\mathcal{E}\cup\mathcal{I}} \lambda_i c_i(\bx), \label{eq_lagrangian}
\end{align}
where $\bx\in\R^n$ and $\blambda\in\R^{|\mathcal{E}|+|\mathcal{I}|}$, provided a suitable \emph{constraint qualification} holds \revision{(see e.g.~\cite[Chapters 12.2 \& 12.6]{Nocedal2006})}, a necessary condition for $\bx^*$ to be a local minimizer of \eqref{eq_generic_cons} is that there exists $\blambda^*\in\R^{|\mathcal{E}|+|\mathcal{I}|}$ such that
\begin{subequations} \label{eq_kkt}
\begin{align}
    \grad_{\bx} L(\bx^*, \blambda^*) &= \bm{0}, \\
    c_i(\bx^*) &= 0, \qquad \forall i\in\mathcal{E}, \\
    c_i(\bx^*) &\leq 0, \qquad \forall i\in\mathcal{I}, \\
    \lambda^*_i &\geq 0, \qquad \forall i\in\mathcal{I}, \\
    \lambda^*_i c_i(\bx^*) &= 0, \qquad \forall i\in\mathcal{I}.
\end{align}
\end{subequations}
For an overview of introductory ideas from nonlinear optimization and more details on the above, see \cite{Nocedal2006}.

\subsection{Trust-Region Methods} \label{sec_tr}

Trust-region methods are a popular class of algorithms for nonconvex optimization which have proven very successful in practice, featuring in the state-of-the-art codes GALAHAD \cite{Gould2003} and KNitro \cite{Byrd2006}, as well as in many of the methods from MATLAB's Optimization Toolbox and SciPy's optimization library, for example.
Other common classes of algorithm are linesearch and regularization methods, which primarily differ in how \emph{global convergence}\footnote{That is, convergence to \minorrev{a stationary point} is guaranteed regardless of how far the algorithm starts from a solution, in contrast with \emph{local convergence}, which assumes a starting point sufficiently close to a solution.} is guaranteed.
All of these classes are iterative, requiring the user to provide a starting point $\bx_0\in\R^n$ and generating a sequence $\bx_1,\bx_2,\ldots\in\R^n$ which we hope will converge to a minimizer.
They are also all based on iterative minimization of local approximations to the objective $f$, usually Taylor series (or approximations thereof).
In trust-region methods, the crucial ingredient that distinguishes it from the other classes is a positive scalar called the \emph{trust-region radius}, updated at each iteration, that constrains the distance between consecutive iterates.

If we have zeroth- and first-order oracles for the objective $f$, the key steps in one iteration of a trust-region method for solving \eqref{eq_problem} are:
\begin{enumerate}
    \item Build a (usually quadratic) model for $f$ that we expect to be accurate near the current iterate. This typically uses a second-order Taylor series for $f$ based at $\bx_k$, with the Hessian potentially replaced with a \emph{quasi-Newton} approximation such as the symmetric rank-1 update (see \cite[Chapter 6]{Nocedal2006});
    \item Minimize the model in a neighborhood of the current iterate (the \emph{trust region}, i.e.~the region in which we trust the model to be accurate);
    \item Evaluate $f$ at the minimizer found in the previous step and decide whether or not to make this the new iterate, and update the trust-region radius. We accept a new iterate if it sufficiently decreases the objective, and increase the trust-region radius if we have been making good progress (to enable larger steps, hence faster progress to the minimizer), or decrease if we are not making progress (since our model \minorrev{is} more accurate in a smaller neighborhood of the iterate).
\end{enumerate}
This is formalized in \algref{alg_basic_tr}.
The two most important aspects of \algref{alg_basic_tr} are the model construction and step calculation in \eqref{eq_tr_model} and \eqref{eq_trs} respectively.
Typical values for the parameters might be $\gammadec=0.5$, $\gammainc=2$, $\eta_U=0.1$ and $\eta_S=0.7$, but these should not be viewed as prescriptive.

\begin{algorithm}[tb]
\begin{algorithmic}[1]
\Require Starting point $\bx_0\in\R^n$ and trust-region radius $\Delta_0>0$. Algorithm parameters: scaling factors $0 < \gammadec < 1 < \gammainc$ and acceptance thresholds $0 < \eta_U \leq \eta_S < 1$.
\For{$k=0,1,2,\ldots$}
    \State Build a local quadratic Taylor-like model for the objective,
    \begin{align}
        f(\by) \approx m_k(\by) := f(\bx_k) + \bg_k^T (\by-\bx_k) + \frac{1}{2}(\by-\bx_k)^T \bH_k (\by-\bx_k), \label{eq_tr_model}
    \end{align}
    for some $\bg_k\in\R^n$ and (symmetric) $\bH_k \approx \grad^2 f(\bx_k) \in\R^{n\times n}$.
    \State Solve the \emph{trust-region subproblem}: \minorrev{set $\bs_k$ to be an approximate minimizer of}
    \begin{align}
        \min_{\bs\in\R^n} m_k(\bx_k+\bs), \qquad \text{s.t.} \quad \|\bs\| \leq \Delta_k. \label{eq_trs}
    \end{align}
    \State Evaluate $f(\bx_k+\bs_k)$ and calculate the ratio 
    \begin{align}
        \rho_k = \frac{\text{actual reduction}}{\text{predicted reduction}} \defeq \frac{f(\bx_k)-f(\bx_k+\bs_k)}{m_k(\bx_k)-m_k(\bx_k+\bs_k)}. \label{eq_ratio_generic}
    \end{align}
    \If{$\rho_k \geq \eta_S$} \label{ln_update_rule_start}
        \State \textit{(Very successful iteration)} Set $\bx_{k+1}=\bx_k+\bs_k$ and $\Delta_{k+1}=\gammainc\Delta_k$. 
    \ElsIf{$\eta_U \leq \rho_k < \eta_S$}
        \State \textit{(Successful iteration)} Set $\bx_{k+1}=\bx_k+\bs_k$ and $\Delta_{k+1}=\Delta_k$.
    \Else
        \State \textit{(Unsuccessful iteration)} Set $\bx_{k+1}=\bx_k$ and $\Delta_{k+1}=\gammadec\Delta_k$. 
    \EndIf \label{ln_update_rule_end}
\EndFor
\end{algorithmic}
\caption{Derivative-based trust-region method for solving \eqref{eq_problem}.}
\label{alg_basic_tr}
\end{algorithm}

\revision{The core calculation in \algref{alg_basic_tr} is the (approximate) solution of the trust-region subproblem \eqref{eq_trs}.}
Ignoring the specific choice of model and centering the model at the current iterate, this subproblem has the general form
\begin{align}
    \min_{\bs\in\R^n} m(\bs) := c + \bg^T \bs + \frac{1}{2}\bs^T \bH \bs, \qquad \text{s.t.} \quad \|\bs\| \leq \Delta. \label{eq_trs_generic}
\end{align}
At first glance, solving \eqref{eq_trs_generic} appears daunting: instead of our original problem of minimizing $f$, now at each iteration we need to minimize a different function, but with constraints!

However, the special structure of \eqref{eq_trs_generic}---minimizing a (possibly nonconvex) quadratic objective subject to single a Euclidean ball constraint---allows the efficient calculation of exact (global) minimizers \revision{via a one-dimensional search over the Lagrange mulitplier for the (squared) constraint $\|\bs\|^2 \leq \Delta^2$.} 

Implementing an efficient algorithm to do this requires some effort; see \cite[Chapter 7.3]{Conn2000} or the newer \cite{Gould2010} for details.
To find a global minimizer, the dimension $n$ must not be too large, as the subproblem solver requires computing Cholesky factorizations of $\bH+\lambda \bI$ for several different values of $\lambda$, with a cost of $\bigO(n^3)$ \minorrev{operations} each time.

\paragraph{Approximate Subproblem Solutions}
However, we do not necessarily need the global solution to the trust-region subproblem for \algref{alg_basic_tr} to converge.
A very simple approximate solution can be found by performing 1 iteration of gradient descent with exact linesearch.
That is, our approximate solution to \eqref{eq_trs_generic} is of the form $\bs(t) = -t\bg$ for $t\geq 0$.
Restricting the objective to this ray, we get $m(\bs(t)) = c -t\|\bg\|^2 + \frac{\bg^T \bH \bg}{2}t^2$ subject to $0\leq t \leq \Delta/\|\bg\|$.
This is easy to globally minimize in $t$, yielding the so-called \emph{Cauchy point}
\begin{align} \label{eq_cauchy_point}
    \bs_C := -t_C \bg, \qquad \text{where} \qquad t_C := \begin{cases} \min\left(\frac{\|\bg\|}{\bg^T \bH \bg}, \frac{\Delta}{\|\bg\|}\right), & \bg^T \bH \bg > 0, \\ \frac{\Delta}{\|\bg\|}, & \bg^T \bH \bg \leq 0. \end{cases}
\end{align}
By direct computation, we can get a lower bound on how much the Cauchy point decreases the quadratic model.

\begin{lemma}[Theorem 6.3.1, \cite{Conn2000}] 
    The Cauchy point \eqref{eq_cauchy_point} solution to \eqref{eq_trs_generic} satisfies
    \begin{align}
        m(\bm{0}) - m(\bs_C) \geq \frac{1}{2} \|\bg\| \min\left(\Delta, \frac{\|\bg\|}{\|\bH\|+1}\right). \label{eq_cauchy_dec}
    \end{align}
\end{lemma}

It turns out that \eqref{eq_cauchy_dec} is in fact \minorrev{sufficient to} achieve first-order convergence of \algref{alg_basic_tr}.
Hence, after re-introducing $\bx_k$ into the model as per \eqref{eq_tr_model}, we make the following assumption regarding \eqref{eq_trs}.

\begin{assumption} \label{ass_cauchy_decrease}
    The computed step $\bs_k$ in \eqref{eq_trs} satisfies $\|\bs_k\| \leq \Delta_k$ and
    \begin{align}
        m_k(\bx_k) - m_k(\bx_k+\bs_k) \geq \kappa_s \|\bg_k\| \min\left(\Delta_k, \frac{\|\bg_k\|}{\|\bH_k\|+1}\right),
    \end{align}
    for some $\kappa_s\in(0,\frac{1}{2})$.
\end{assumption}

If $n$ is large \minorrev{and} computation of the global minimizer of \eqref{eq_trs_generic} is impractical, there are still alternative methods for improving on the Cauchy point.
The most common is an adaptation of the conjugate gradient (CG) method for solving symmetric positive definite linear systems \cite[Chapter 5.1]{Nocedal2006}.
In this method, sometimes called the \emph{Steihaug-Toint method}, we use CG to solve $\bH\bs=-\bg$ starting from initial iterate $\bs_0=\bm{0}$.
If any iteration would take us outside the feasible region, we truncate the step so we remain on the boundary of the feasible region and terminate.
Similarly, if, at any iteration, we compute a search direction $\bd$ such that $\bd^T \bH \bd \leq 0$, then $\bH$ is not positive definite, and we move from our current CG iterate in the direction $\bd$ until we reach the boundary of the feasible region.
This method has the advantage of only requiring $\bH$ via Hessian-vector products, and so is suitable for large-scale problems, and the first iterate is always the Cauchy point, so \assref{ass_cauchy_decrease} is guaranteed.
More details of this and other approximate subproblem solvers can be found in \cite[Chapter 7.5]{Conn2000}.

\revision{To achieve second-order convergence of \algref{alg_basic_tr}, we first note that the quadratic model $m_k$ \eqref{eq_tr_model} has its own second-order optimality measure, namely
\begin{align}
    \sigma^m_k := \max(\|\bg_k\|, \tau^m_k), \qquad \text{with} \qquad \tau^m_k := \max(-\lambda_{\min}(\bH_k), 0). \label{eq_optimality_model_2}
\end{align}
although here we will assume $\bg_k=\grad f(\bx_k)$ and $\bH_k=\grad^2 f(\bx_k)$, and so $\sigma^m_k = \sigma_k$ and $\tau^m_k=\tau_k$.
Here, we need a stronger assumption on our trust-region subproblem solution, to handle the case where $\bH_k$ is not positive semidefinite.\footnote{For example, if $\bg_k=\bm{0}$ and $\lambda_{\min}(\bH_k) < 0$ then $\bs_k=\bm{0}$ satisfies \assref{ass_cauchy_decrease}.}}
Suppose that $\tau^m_k>0$ (i.e.~$\lambda_{\min}(\bH_k) < 0$), and let $\bu_k$ be a normalized eigenvector corresponding to that eigenvalue, which is also a first-order descent direction; that is,
\begin{align}
    \bH_k \bu_k = \lambda_{\min}(\bH_k) \bu_k, \qquad \|\bu_k\| = 1, \qquad \text{and} \qquad \bu_k^T \bg_k \leq 0.
\end{align}
The last condition may be achieved by judicious sign choice in the eigenvector calculation.
Similar to the definition of the Cauchy point \eqref{eq_cauchy_point}, the \emph{eigenstep} $\bs_k^E$ is \minorrev{the point} in the direction $\bu_k$ that minimizes the model $m_k$ (subject to $\|\bs_E\| \leq \Delta_k$).
Since we have
\begin{align}
    m_k(\bx_k+t\bu_k) = f(\bx_k) + t \bu_k^T \bg_k + \frac{1}{2} \bu_k^T \bH_k \bu_k t^2 = f(\bx_k) + t \bu_k^T \bg_k + \frac{1}{2} \lambda_{\min}(\bH_k) t^2,
\end{align}
we conclude that $m_k(\bx_k+t\bu_k)$ is decreasing as $t\geq 0$ increases, so we have $\bs_k^E = \Delta_k \bu_k$ since $\|\bu_k\|=1$.
We then compute
\begin{align}
    m_k(\bx_k) - m_k(\bx_k+\bs_k^E) = -\Delta_k \bu_k^T \bg_k - \frac{1}{2} \lambda_{\min}(\bH_k) \Delta_k^2 \geq -\frac{1}{2}\lambda_{\min}(\bH_k) \Delta_k^2 \geq 0.
\end{align}
This motivates the following assumption on our step calculation for second-order convergence, where we now require our step to be at least as good as the Cauchy step and, if $\tau^m_k>0$ (i.e.~$\lambda_{\min}(\bH_k) < 0$), the eigenstep.

\begin{assumption} \label{ass_eigenstep_decrease}
    The computed step $\bs_k$ in \eqref{eq_trs} satisfies $\|\bs_k\| \leq \Delta_k$ and
    \begin{align}
        m_k(\bx_k) - m_k(\bx_k+\bs_k) \geq \kappa_s \max\left(\|\bg_k\| \min\left(\Delta_k, \frac{\|\bg_k\|}{\|\bH_k\|+1}\right), \tau^m_k\Delta_k^2\right),
    \end{align}
    for some $\kappa_s\in(0, \frac{1}{2})$, where $\tau^m_k$ is defined in \eqref{eq_optimality_model_2}.
\end{assumption}

\subsection{Convergence of Trust-Region Methods}

\revision{
\revision{We conclude this section by summarizing} the convergence of \algref{alg_basic_tr} to first- and second-order optimal points.
Our results will also show the \emph{worst-case complexity} of \algref{alg_basic_tr}, which essentially is a rate of convergence (i.e.~how many iterations are needed to achieve $\|\grad f(\bx_k)\| < \epsilon$ or $\sigma_k < \epsilon$?).
\revision{We assume the following about the quadratic model \eqref{eq_tr_model}.}

\begin{assumption} \label{ass_model}
    At each iteration $k$ of \algref{alg_basic_tr}, the model $m_k$ \eqref{eq_tr_model} satisfies:
    \begin{enumerate}[label=(\alph*)]
        \item $\bg_k = \grad f(\bx_k)$; \label{ass_model_g}
        \item $\|\bH_k\| \leq \kappa_H - 1$ for some $\kappa_H\geq 1$ (independent of $k$).\footnote{This will be more notationally convenient than the equivalent but more natural $\|\bH_k\| \leq \kappa_H$, for some $\kappa_H\geq 0$.} \label{ass_model_H}
    \end{enumerate}
\end{assumption}

Our main first-order convergence result is the following.

\begin{theorem}[Theorem 2.3.7, \cite{Cartis2022}] \label{thm_wcc}
    Suppose Assumptions~\ref{ass_smoothness_1}, \minorrev{\ref{ass_cauchy_decrease} and \ref{ass_model}} hold and we run \algref{alg_basic_tr}.
    If $k_{\epsilon}$ is the first iteration of \algref{alg_basic_tr} such that $\|\grad f(\bx_{k_\epsilon})\| < \epsilon$, then $k_{\epsilon} = \bigO((\Lgrad+\kappa_H) \epsilon^{-2})$.
    Hence $\liminf_{k\to\infty} \|\grad f(\bx_k)\| = 0$.
\end{theorem}



This tells us that there is a subsequence of iterations whose gradients converge to zero, and that we get $\|\grad f(\bx_k)\|<\epsilon$ for the first time after at most $\bigO(\epsilon^{-2})$ iterations.
In practice, we typically terminate \algref{alg_basic_tr} when $\|\bg_k\|=\|\grad f(\bx_k)\|$ is sufficiently small or we exceed some computational budget (e.g.~maximum number of iterations).

\revision{\algref{alg_basic_tr} to converge to second-order optimal solutions, we have to strengthen \assref{ass_model} to also require $\bH_k=\grad^2 f(\bx_k)$.}

\begin{theorem}[Theorem 3.2.6, \cite{Cartis2022}] \label{thm_wcc_2}
    Suppose Assumptions~\ref{ass_smoothness_2}, \minorrev{\ref{ass_eigenstep_decrease} and \ref{ass_model}} hold, and $\bH_k=\grad^2 f(\bx_k)$ for all $k$, and we run \algref{alg_basic_tr}.
    If $k_{\epsilon}$ is the first iteration of \algref{alg_basic_tr} such that $\sigma_{k_{\epsilon}} \leq \epsilon$, then $k_{\epsilon} = \bigO((\Lhess+ \kappa_H)^2 \epsilon^{-3})$.
    Hence $\liminf_{k\to\infty} \sigma_k = 0$.
\end{theorem}



In summary, if we want to achieve optimality level $\epsilon$, then trust-region methods require $\bigO(\epsilon^{-2})$ iterations to achieve first-order optimality\footnote{This follows from \thmref{thm_wcc}, but the same result holds if we assume exact second-order models, $\bH_k=\grad^2 f(\bx_k)$ \cite[Theorem 3.2.1]{Cartis2022}.} and $\bigO(\epsilon^{-3})$ iterations to achieve second-order optimality, assuming sufficient problem smoothness, model accuracy and subproblem solution quality.
} 

\subsubsection*{Notes and References}
{\small The material in this section is based on the widely used books \cite{Conn2000,Nocedal2006} and the more recent book on complexity theory \cite{Cartis2022}.
A more recent survey of trust-region methods is \cite{Yuan2015}.

Although \thmref{thm_wcc_2} considers the case where we have the same desired accuracy level $\epsilon$ for both first- and second-order optimality, \revision{the proof of \thmref{thm_wcc_2}} gives a complexity of $\bigO(\max(\epsilon_1^{-2} \epsilon_2^{-1}, \epsilon_2^{-3}))$ as $\epsilon_1,\epsilon_2\to 0^{+}$ if we want $\|\grad f(\bx_k)\| \leq \epsilon_1$ and $\tau_k \leq \epsilon_2$.
With minor algorithmic changes we can get the more natural complexity bound $\bigO(\max(\epsilon_1^{-2}, \epsilon_2^{-3}))$ \cite{Gratton2020}.}

\section{Derivative-Free Trust-Region Methods} \label{sec_dfotr}

We now consider how to adapt the generic trust-region method from \secref{sec_tr} to the case where we only have a zeroth order oracle for the objective.
That is, we are solving \eqref{eq_problem}, where $f$ is still differentiable (e.g.~\assref{ass_smoothness_1}), but where we do not have access to $\grad f$ in our algorithm.

As before, at each iteration we will construct a local quadratic model around $\bx_k$, namely
\begin{align}
    f(\by) \approx m_k(\by) := c_k + \bg_k^T (\by-\bx_k) + \frac{1}{2}(\by-\bx_k)^T \bH_k (\by-\bx_k), \label{eq_tr_model_dfo}
\end{align}
noting that, unlike \eqref{eq_tr_model}, we do not require $m_k(\bx_k) = f(\bx_k)$ (i.e.~$c_k \neq f(\bx_k)$ is allowed).
For now, we will not define a specific way to construct the local quadratic model \eqref{eq_tr_model_dfo} using only function values.
Instead, we will focus on simply defining some practical assumptions on the model accuracy (i.e.~replacing~\assref{ass_model}).
In Sections~\ref{sec_model_construction} and \ref{sec_self_correcting} we will outline concrete ways our model accuracy requirements can be specified.

\paragraph{Motivation: finite difference gradients}
In the simplest case, imagine we do not have access to $\grad f$, and simply build $\bg_k$ in \eqref{eq_tr_model_dfo} from \minorrev{(forward)} finite differences:
\begin{align}
    [\bg_k]_i := \frac{f(\bx_k + h\be_i) - f(\bx_k)}{h}, \qquad \forall i=1,\ldots,n, \label{eq_fin_diff}
\end{align}
for some small $h>0$, where $\be_i\in\R^n$ is the $i$-th coordinate vector. 
From standard analysis of finite differencing (e.g.~\cite[Chapter 8.1]{Nocedal2006}), if $f$ is twice continuously differentiable with $\|\grad^2 f(\bx)\| \leq M$ for all $\bx$ in $B(\bx_k,h)$, then $\|\bg_k - \grad f(\bx_k)\|_{\infty} \leq \frac{Mh}{2}$.
That is, the gradient error is of size $\bigO(h)$ as $h\to 0^{+}$.
Although derivative-based trust-region methods can handle inexact gradients \cite[Chapter 8.4]{Conn2000}, they require a bound on the \emph{relative} gradient error, whereas finite differencing gives an \emph{absolute} error in the model gradient.
We can only control the model error by choosing the value of $h$, perhaps differently at each iteration.
It is not obvious how we could pick $h$ if we instead wanted to control the relative gradient error, at least not without already having a good estimate of $\|\bg_k\|$ and $M$ (i.e.~unless we already have the derivative information we are trying to calculate!).

\revision{Now, suppose we build a quadratic model \eqref{eq_tr_model_dfo} using this $\bg_k$ and any bounded Hessian $\|\bH_k\| \leq \kappa_H-1$ (c.f.~\assref{ass_model}\ref{ass_model_H}).}
Given $\|\bg_k-\grad f(\bx_k)\|_{\infty} \leq \frac{Mh}{2}$, and hence $\|\bg_k-\grad f(\bx_k)\| \leq \frac{Mh \sqrt{n}}{2}$, the error in the model \eqref{eq_tr_model_dfo} at a point $\by\in\R^n$ is
\begin{align}
    &\revision{|f(\by) - f(\bx_k) - \bg_k^T (\by-\bx_k) - \frac{1}{2}(\by-\bx_k)^T \bH_k (\by-\bx_k)|} \nonumber \\
    &\qquad \leq |f(\by) - f(\bx_k) - \grad f(\bx_k)^T (\by-\bx_k)| + |(\grad f(\bx_k) - \bg_k)^T (\by-\bx_k)| \revision{+ \frac{\kappa_H}{2}\|\by-\bx_k\|^2}, \\
    &\qquad \leq \frac{\Lgrad \revision{+ \kappa_H}}{2} \|\by-\bx_k\|^2 + \frac{Mh \sqrt{n}}{2} \|\by-\bx_k\|.
\end{align}
So, if we want to approximate the objective within the trust region, $\by\in B(\bx_k,\Delta_k)$, it makes sense to balance both error terms and set $h = \bigO(\Delta_k)$, leading to a model error of size $\bigO(\Delta_k^2)$.

\paragraph{Fully linear models}
In our general \renaming{MBDFO} approach, we will use the trust-region radius $\Delta_k$ to control the size of the model error (e.g.~at iteration $k$, use $h=\Delta_k$ in \eqref{eq_fin_diff}, as suggested above).
This suggests to us the following notion of model accuracy.

\begin{definition} \label{def_fully_linear}
    Suppose we have $\bx\in\R^n$ and $\Delta>0$.
    A local model $m:\R^n\to\R$ approximating $f:\R^n\to\R$ is fully linear in $B(\bx,\Delta)$ if there exist constants $\kappamf,\kappamg>0$, independent of $m$, $\bx$ and $\Delta$, such that
    \begin{subequations} \label{eq_fully_linear}
    \begin{align}
        |m(\by) - f(\by)| &\leq \kappamf \Delta^2, \label{eq_fully_linear_f} \\
        \|\grad m(\by) - \grad f(\by)\| &\leq \kappamg \Delta, \label{eq_fully_linear_g}
    \end{align}
    \end{subequations}
    for all $\by\in B(\bx,\Delta)$.
    Sometimes we will use $\kappam := \max(\kappamf,\kappamg)$ for notational convenience.
\end{definition}


\begin{remark}
    \revision{`Fully linear' here refers not to the model being linear, but the model being as accurate an approximation (up to constants) as a linear Taylor series. Quadratic models as in \eqref{eq_tr_model_dfo} can be fully linear. 
    However, in Sections~\ref{sec_model_construction} and \ref{sec_self_correcting} we will sometimes consider the case where $m_k$ \eqref{eq_tr_model_dfo} is indeed linear (i.e.~$\bH_k=\bm{0}$). For simplicity, we will generically refer to $m_k$ as a `quadratic model' even if $\bH_k=\bm{0}$, unless we specifically wish to draw attention to $m_k$ being linear.}
\end{remark}

If we do indeed have access to a first-order oracle (i.e.~not in the \renaming{MBDFO} case), then we can satisfy \defref{def_fully_linear} using the approaches from \secref{sec_tr}.
For example, if \minorrev{\assref{ass_smoothness_1}\ref{ass_smoothness_1_smooth}} holds (and so we have \lemref{lem_lsmooth}):
\begin{itemize}
    \item The first-order Taylor series $m(\by) = f(\bx) + \grad f(\bx)^T (\by-\bx)$ is fully linear in $B(\bx,\Delta)$ with $\kappamf=\frac{1}{2}\Lgrad$ and $\kappamg = \Lgrad$.
    \item More generally, any quadratic model satisfying \assref{ass_model} is fully linear in $B(\bx,\Delta)$ with $\kappamf = \frac{1}{2}(\Lgrad + \kappa_H)$ and $\kappamg = \Lgrad + \kappa_H$.
\end{itemize}
In Sections~\ref{sec_model_construction} and \ref{sec_self_correcting} we will show how models satisfying \defref{def_fully_linear} may be constructed using only function values.

\begin{assumption} \label{ass_model_dfo}
    At each iteration $k$ of \algref{alg_basic_tr_dfo}, the model $m_k$ \eqref{eq_tr_model_dfo} satisfies:
    \begin{enumerate}[label=(\alph*)]
        \item $m_k$ is fully linear in $B(\bx_k,\Delta_k)$ with constants $\kappamf,\kappamg>0$ independent of $k$; \label{ass_model_dfo_g}
        \item $\|\bH_k\| \leq \kappa_H - 1$ for some \minorrev{fixed} $\kappa_H\geq 1$ (independent of $k$). \label{ass_model_dfo_H}
    \end{enumerate}
\end{assumption}

It is clear from \eqref{eq_fully_linear} that we have to contend with absolute model errors, controlled by the size of $\Delta$ (which will be set to the trust-region radius in our algorithm).
This introduces several difficulties that are not present in the derivative-based case:
\begin{itemize}
    \item $\Delta_k$ now performs two roles: it controls the size of the tentative step ($\|\bs_k\| \leq \Delta_k$ in \eqref{eq_trs}) and the size of the model error \eqref{eq_fully_linear}.\footnote{In some algorithms (e.g.~\cite{Powell2002}), two trust region radii are used to partially decouple these roles.}
    \item If our model suggests we are close to a first-order solution (i.e.~$\|\bg_k\| \approx 0$), this does not mean that we are \emph{actually} near a solution (i.e.~$\|\grad f(\bx_k)\| \approx 0$). 
    By comparison, if we have relative errors in the model gradient, then $\|\bg_k\| \approx 0$ if and only if $\|\grad f(\bx_k)\| \approx 0$ \cite[Lemma 8.4.1]{Conn2000}.
    \item As a consequence of the above point, it is no longer clear when to terminate the algorithm in practice.
    If we had gradients (or relative gradient errors), then terminating when $\|\bg_k\|$ is sufficiently small guarantees we are close to a (first-order) solution. We specifically discuss termination in \secref{sec_termination}.
\end{itemize}
Unfortunately, this additional complexity is a consequence of having a weaker notion of model accuracy, one which can be practically achieved using only function evaluations.

\subsection{First-Order Convergence}
We are now ready to state our \renaming{MBDFO} variant of \algref{alg_basic_tr}, where our Taylor-based models (i.e.~\eqref{eq_tr_model} satisfying \assref{ass_model}) are replaced with fully linear models (i.e.~\eqref{eq_tr_model_dfo} satisfying \assref{ass_model_dfo}).
It is almost identical to the derivative-based version \algref{alg_basic_tr}, including allowing inexact subproblem solutions (satisfying the Cauchy decrease condition, \assref{ass_cauchy_decrease}).

Aside from our new model accuracy condition (\assref{ass_model_dfo} instead of \assref{ass_model}), the only difference is that we need to explicitly check our \emph{criticality measure}\footnote{That is, our measure of optimality, in this case $\|\bg_k\|$ (for first-order convergence).}, and require $\|\bg_k\| \geq \mu_c \Delta_k$ for a step to be declared (very) successful.
This is an important feature of \renaming{MBDFO} methods, and aims to address the decoupling of our measured distance to a solution $\|\bg_k\|$ from our true distance to solution $\|\grad f(\bx_k)\|$.
This mechanism ensures that, for successful iterations, if $\|\grad f(\bx_k)\|$ is large (i.e.~we are far from first-order optimality) then so is $\|\bg_k\|$ (i.e.~the model knows we are far from optimality); see \eqref{eq_crit_comparison} below.

\begin{algorithm}[tb]
\begin{algorithmic}[1]
\Require Starting point $\bx_0\in\R^n$ and trust-region radius $\Delta_0>0$. Algorithm parameters: scaling factors $0 < \gammadec < 1 < \gammainc$, acceptance thresholds $0 < \eta_U \leq \eta_S < 1$, and criticality threshold $\mu_c > 0$.
\For{$k=0,1,2,\ldots$}
    \State Build a local quadratic model $m_k$ \eqref{eq_tr_model_dfo} satisfying \assref{ass_model_dfo}.
    \State Solve the trust-region subproblem \eqref{eq_trs} to get a step $\bs_k$ satisfying \assref{ass_cauchy_decrease}.
    \State Evaluate $f(\bx_k+\bs_k)$ and calculate the ratio $\rho_k$ \eqref{eq_ratio_generic}.
    \If{$\rho_k \geq \eta_S$ and $\|\bg_k\| \geq \mu_c \Delta_k$} 
        \State \textit{(Very successful iteration)} Set $\bx_{k+1}=\bx_k+\bs_k$ and $\Delta_{k+1}=\gammainc\Delta_k$. 
    \ElsIf{$\eta_U \leq \rho_k < \eta_S$ and $\|\bg_k\| \geq \mu_c \Delta_k$}
        \State \textit{(Successful iteration)} Set $\bx_{k+1}=\bx_k+\bs_k$ and $\Delta_{k+1}=\Delta_k$.
    \Else
        \State \textit{(Unsuccessful iteration)} Set $\bx_{k+1}=\bx_k$ and $\Delta_{k+1}=\gammadec\Delta_k$. 
    \EndIf 
\EndFor
\end{algorithmic}
\caption{Simple \renaming{MBDFO} trust-region method for solving \eqref{eq_problem}.}
\label{alg_basic_tr_dfo}
\end{algorithm}

\begin{remark} \label{rem_criticality}
    If $\|\bg_k\| < \mu_c \Delta_k$ at any iteration, then that iteration must be unsuccessful regardless of the value of $\rho_k$.
    So, if this check is performed as soon as the model is constructed, we may save our efforts by not  solving the trust-region subproblem or evaluating $f(\bx_k+\bs_k)$.
\end{remark}

We now present our analysis of \algref{alg_basic_tr_dfo}.

\begin{lemma} \label{lem_very_successful_dfo}
    Suppose Assumptions~\minorrev{\ref{ass_smoothness_1}\ref{ass_smoothness_1_smooth}}, \minorrev{\ref{ass_cauchy_decrease} and \ref{ass_model_dfo}} hold.
    \minorrev{If, on iteration $k$ of \algref{alg_basic_tr_dfo},} $\bg_k\neq\bm{0}$ and
    \begin{align}
        \Delta_k \leq \min\left(\frac{\kappa_s (1-\eta_S)}{2\kappamf}, \frac{1}{\kappa_H}, \frac{1}{\mu_c}\right) \|\bg_k\|, 
    \end{align}
    then $\rho_k \geq \eta_S$ and $\|\bg_k\| \geq \mu_c \Delta_k$ (i.e.~iteration $k$ is very successful).
\end{lemma}
\begin{proof}
    That $\|\bg_k\| \geq \mu_c \Delta_k$ follows immediately by assumption on $\Delta_k$.
    It remains to show $\rho_k \geq \eta_S$. 
    From \assref{ass_model_dfo} we have $\Delta_k \leq \frac{\|\bg_k\|}{\kappa_H} \leq \frac{\|\bg_k\|}{\|\bH_k\|+1}$, and so \assref{ass_cauchy_decrease} gives
    \begin{align}
        m_k(\bx_k) - m_k(\bx_k+\bs_k) \geq \kappa_s \|\bg_k\| \Delta_k.
    \end{align}
    Then using \assref{ass_model_dfo} we can compute
    \begin{align}
        |\rho_k - 1| &\leq \frac{|f(\bx_k+\bs_k)-m_k(\bx_k+\bs_k)| + |f(\bx_k) - m_k(\bx_k)|}{|m_k(\bx_k)-m_k(\bx_k+\bs_k)|}, \\
        &\leq \frac{2\kappamf \Delta_k^2}{\kappa_s \|\bg_k\| \Delta_k}, \\
        &= \frac{2\kappamf \Delta_k}{\kappa_s \|\bg_k\|},
    \end{align}
    and so $|\rho_k-1| \leq 1-\eta_S$ by assumption on $\Delta_k$, which implies $\rho_k \geq \eta_S$ \minorrev{as required}.
\end{proof}

The crucial consequence of \lemref{lem_very_successful_dfo} is that $\Delta_k$ remains large provided we have not yet converged.

\begin{lemma} \label{lem_delta_min_dfo}
    Suppose Assumptions~\minorrev{\ref{ass_smoothness_1}\ref{ass_smoothness_1_smooth}}, \minorrev{\ref{ass_cauchy_decrease} and \ref{ass_model_dfo}} hold and we run \algref{alg_basic_tr_dfo}.
    If $\|\grad f(\bx_k)\| \geq \epsilon$ for all $k=0,\ldots,K-1$, then
    \begin{align}
        \Delta_k \geq \Delta_{\min}(\epsilon) &:= \min\left(\Delta_0, \frac{\gammadec \epsilon}{\max\left(\frac{2\kappamf}{\kappa_s (1-\eta_S)}, \kappa_H, \mu_c\right) + \kappamg}\right), \label{eq_delta_min_dfo}
    \end{align}
    for all $k=0,\ldots,K$.
\end{lemma}
\begin{proof}
    We proceed by induction.
    The result holds trivially for $k=0$, so suppose $\Delta_k \geq \Delta_{\min}(\epsilon)$ for some $k\in\{0,\ldots,K-1\}$.
    To find a contradiction assume that $\Delta_{k+1} < \Delta_{\min}(\epsilon)$.
    Then $\Delta_{k+1} < \Delta_k$, which by the mechanism for updating the trust-region radius means iteration $k$ was unsuccessful, and $\Delta_k = \gammadec^{-1} \Delta_{k+1} < \gammadec^{-1} \Delta_{\min}(\epsilon)$.
    From \lemref{lem_very_successful_dfo}, this means we must have
    \begin{align}
        \Delta_k > \min\left(\frac{\kappa_s (1-\eta_S)}{2\kappamf}, \frac{1}{\kappa_H}, \frac{1}{\mu_c}\right) \|\bg_k\|.
    \end{align}
    But since $m_k$ is fully linear (\assref{ass_model_dfo}), we get
    \begin{align}
        \epsilon \leq \|\grad f(\bx_k)\| \leq \|\bg_k\| + \|\bg_k - \grad f(\bx_k)\| < \max\left(\frac{2\kappamf}{\kappa_s (1-\eta_S)}, \kappa_H, \mu_c\right) \Delta_k + \kappamg \Delta_k,
    \end{align}
    which contradicts $\Delta_k < \gammadec^{-1} \Delta_{\min}(\epsilon)$.
\end{proof}

Our main convergence result is the following.
The arguments here can largely be used to prove the equivalent result for derivative-based methods, \thmref{thm_wcc}.

\begin{theorem} \label{thm_wcc_dfo}
    Suppose Assumptions~\ref{ass_smoothness_1}, \minorrev{\ref{ass_cauchy_decrease} and \ref{ass_model_dfo}} hold and we run \algref{alg_basic_tr_dfo}.
    If $\|\grad f(\bx_k)\| \geq \epsilon$ for all $k=0,\ldots,K-1$, then 
    \begin{align}
        K \leq \frac{\log(\Delta_0 / \Delta_{\min}(\epsilon))}{\log(\gammadec^{-1})} + \left(1 + \frac{\log(\gammainc)}{\log(\gammadec^{-1})}\right) \revision{\frac{f(\bx_0)-\flow}{\eta_U \kappa_s \mu_c \min(1, \mu_c/\kappa_H) \Delta_{\min}(\epsilon)^2}},
    \end{align}
    where $\Delta_{\min}(\epsilon)$ is defined in \lemref{lem_delta_min_dfo}.
\end{theorem}
\begin{proof}
    We first partition the iterations $\{0,\ldots,K-1\}=\mathcal{S}\cup\mathcal{U}$, where $\mathcal{S}$ is the set of successful or very successful iterations and $\mathcal{U}$ is the set of unsuccessful iterations.
    From \lemref{lem_delta_min_dfo} we have $\Delta_k \geq \Delta_{\min}(\epsilon)$ for all $k\in\{0,\ldots,K-1\}$.

    Since $\bx_{k+1}\neq\bx_k$ if and only if $k\in\mathcal{S}$, we have
    \begin{align}
        f(\bx_0) - \flow \geq f(\bx_0) - f(\bx_K) \geq \sum_{k=0}^{K-1} f(\bx_k) - f(\bx_{k+1}) &= \sum_{k\in\mathcal{S}} f(\bx_k) - f(\bx_{k+1}), \\
        &\geq \eta_U \sum_{k\in\mathcal{S}} m_k(\bx_k) - m_k(\bx_k+\bs_k), \label{eq_wcc_success}
    \end{align}
    where the last inequality uses $\rho_k \geq \eta_U$ for all $k\in\mathcal{S}$.
    \revision{If $k\in\mathcal{S}$, then Assumptions~\ref{ass_cauchy_decrease} and \ref{ass_model_dfo}\ref{ass_model_dfo_H}, together with the criticality requirement $\|\bg_k\| \geq \mu_c \Delta_k$, imply
    \begin{align}
        m_k(\bx_k) - m_k(\bx_k+\bs_k) &\geq \kappa_s \|\bg_k\| \min\left(\Delta_{\min}(\epsilon), \frac{\|\bg_k\|}{\|\bH_k\|+1}\right), \\
        &\geq \kappa_s \mu_c \Delta_{\min}(\epsilon)\min\left(\Delta_{\min}(\epsilon), \frac{\mu_c \Delta_{\min}(\epsilon)}{\kappa_H}\right). \label{eq_wcc_dfo_model_dec}
    \end{align}
    We conclude that
    \begin{align}
        f(\bx_0) - \flow \geq \eta_U \kappa_s \mu_c \min(1, \mu_c/\kappa_H) \Delta_{\min}(\epsilon)^2 \cdot |\mathcal{S}|, \label{eq_wcc_dfo_tmp1a}
    \end{align}
    or
    \begin{align}
        |\mathcal{S}| \leq \frac{f(\bx_0)-\flow}{\eta_U \kappa_s \mu_c \min(1, \mu_c/\kappa_H) \Delta_{\min}(\epsilon)^2}. \label{eq_wcc_dfo_tmp1}
    \end{align}
    }
    Separately, the mechanism for updating $\Delta_k$ ensures that $\Delta_{k+1} = \gammadec \Delta_k$ if $k\in\mathcal{U}$ and $\Delta_{k+1} \leq \gammainc \Delta_k$ if $k\in\mathcal{S}$.
    Hence
    \begin{align}
        \Delta_{\min}(\epsilon) \leq \Delta_{K} = \Delta_0 \gammadec^{|\mathcal{U}|} \gammainc^{|\mathcal{S}|}, \label{eq_wcc_tmp1a}
    \end{align}
    which gives
    \begin{align}
        |\mathcal{U}| \leq \frac{\log(\Delta_0 / \Delta_{\min}(\epsilon))}{\log(\gammadec^{-1})} + \frac{\log(\gammainc)}{\log(\gammadec^{-1})}|\mathcal{S}|. \label{eq_wcc_tmp2}
    \end{align}
    The result then follows from $K=|\mathcal{S}| + |\mathcal{U}|$.
\end{proof}

We again summarize our result in terms of the key quantities of interest.
Here, we recall $\kappam := \max(\kappamf, \kappamg)$ (see \defref{def_fully_linear}) and note that $\Delta_{\min}(\epsilon) = \Theta(\revision{(\kappam+\kappa_H)^{-1}} \epsilon)$ as $\epsilon\to 0^{+}$ and $\kappam,\kappa_H\to\infty$, to get the following first-order convergence result.

\begin{corollary} \label{cor_wcc_dfo}
    Suppose the assumptions of \thmref{thm_wcc_dfo} hold.
    If $k_{\epsilon}$ is the first iteration of \algref{alg_basic_tr_dfo} such that $\|\grad f(\bx_{k_\epsilon})\| < \epsilon$, then $k_{\epsilon} = \bigO(\revision{\kappa_H (\kappam+\kappa_H)^2} \epsilon^{-2})$.
    Hence $\liminf_{k\to\infty} \|\grad f(\bx_k)\| = 0$.
\end{corollary}

\revision{In the derivative-based case (\thmref{thm_wcc}), using $\kappam=\bigO(\Lgrad)$ corresponding to a Taylor model, the worst-case complexity was $\bigO((\kappam+\kappa_H) \epsilon^{-2})$ iterations. The bound in \corref{cor_wcc_dfo} is worse by a factor of $\kappa_H (\kappam+\kappa_H)$, which arises from using $\|\bg_k\| \geq \mu_c \Delta_{\min}(\epsilon)$ to get \eqref{eq_wcc_dfo_model_dec} (see \remref{rem_crit_options}).}

\revision{We also have to assess the impact of using interpolation-based models on $\kappam$. In \secref{sec_model_construction}, we show that $\kappam$ is usually of size $\bigO(\Lgrad)$, the same as for Taylor models, but with a constant that depends explicitly on the problem dimension $n$, e.g.~$\kappam=\bigO(\sqrt{n}\: \Lgrad)$ in \thmref{thm_fully_linear}. In some special cases, though, we can recover dimension-independent $\kappam$ (\corref{cor_fully_linear_wcc_improved}). However, we will always get an explicit dimension dependency of at least $\bigO(n)$ if we count the number of (possibly expensive) objective evaluations.\footnote{\revision{In the derivative-based setting, this does not arise because every evaluation of $\grad f$ gives $n$ pieces of information}.}}

\revision{
\begin{remark} \label{rem_crit_options}
There are multiple ways to establish a model decrease of at least $\bigO(\epsilon^{-2})$ on successful iterations, as in \eqref{eq_wcc_dfo_model_dec} above. 
In derivative-based trust-region methods, when $\bg_k=\grad f(\bx_k)$, we have
\begin{align}
    m_k(\bx_k) - m_k(\bx_k+\bs_k) \geq \kappa_s \epsilon \min\left(\Delta_{\min}(\epsilon), \frac{\epsilon}{\kappa_H}\right) = \kappa_s \epsilon \Delta_{\min}(\epsilon),
\end{align}
immediately from the assumption $\|\grad f(\bx_k)\| \geq \epsilon$, and we never need to check $\|\bg_k\| \geq \mu_c \Delta_k$ (i.e.~we may set $\mu_c=0$).
In the DFO setting, if we know that $m_k$ is fully linear then we may instead use
\begin{align}
    \epsilon \leq \|\grad f(\bx_k)\| \leq \|\bg_k\| + \|\grad f(\bx_k) - \bg_k\| \leq \|\bg_k\| + \kappamg \Delta_k \leq (1+\kappamg \mu_c^{-1}) \|\bg_k\|, \label{eq_crit_comparison}
\end{align}
to get $\|\bg_k\| \geq \bigO(\epsilon)$.
The approach in \eqref{eq_wcc_dfo_model_dec} is the most general, as it only requires the $\|\bg_k\| \geq \mu_c \Delta_k$ and does not require anything regarding model accuracy.
The same flexibility appears in the second-order convergence theory (below).
\end{remark}
}

\subsection{Second-Order Convergence}
Similar to \thmref{thm_wcc_2}, we now consider how to extend the above results to converge to second-order critical points.
Just as in the derivative-based case, we need to assume $f$ is twice continuously differentiable (specifically, \assref{ass_smoothness_2}), and our trust-region subproblem solver needs to achieve at least the eigenstep decrease (\assref{ass_eigenstep_decrease}).

However, motivated by \eqref{eq_lh_smooth}, we need a stricter requirement on our model accuracy than fully linear.
Where fully linear models match (up to constants) the error from a first-order Taylor series, we now require models which match the error from a second-order Taylor series.

\begin{definition} \label{def_fully_quadratic}
    Suppose we have $\bx\in\R^n$ and $\Delta>0$.
    A local model $m:\R^n\to\R$ approximating $f:\R^n\to\R$ is fully quadratic in $B(\bx,\Delta)$ if there exist constants $\kappamf,\kappamg,\minorrev{\kappamh}>0$, independent of $m$, $\bx$ and $\Delta$, such that
    \begin{subequations} \label{eq_fully_quadratic}
    \begin{align}
        |m(\by) - f(\by)| &\leq \kappamf \Delta^3, \label{eq_fully_quadratic_f} \\
        \|\grad m(\by) - \grad f(\by)\| &\leq \kappamg \Delta^2, \label{eq_fully_quadratic_g} \\
        \|\grad^2 m(\by) - \grad^2 f(\by)\| &\leq \kappamh \Delta, \label{eq_fully_quadratic_h}
    \end{align}
    \end{subequations}
    for all $\by\in B(\bx,\Delta)$.
    Sometimes we will use $\kappam := \max(\kappamf,\kappamg,\kappamh)$ for notational convenience.
\end{definition}


If \assref{ass_smoothness_2} and hence \eqref{eq_lh_smooth} and \eqref{eq_lh_smooth_gradf} hold, it is not hard to verify that the second-order Taylor series $m(\by) = f(\bx) + \grad f(\bx)^T (\by-\bx) + \frac{1}{2}(\by-\bx)^T \grad^2 f(\bx) (\by-\bx)$ is fully quadratic in $B(\bx,\Delta)$ with $\kappamf=\frac{1}{6}\Lhess$, $\kappamg=\frac{1}{2}\Lhess$ and $\kappamh=\Lhess$.

Our model assumptions mimic \assref{ass_model_dfo} closely.

\begin{assumption} \label{ass_model_dfo_2}
    At each iteration $k$ of \algref{alg_basic_tr_dfo_2}, the model $m_k$ \eqref{eq_tr_model_dfo} satisfies:
    \begin{enumerate}[label=(\alph*)]
        \item $m_k$ is fully quadratic in $B(\bx_k,\Delta_k)$ with constants $\kappamf,\kappamg,\kappamh>0$ independent of $k$; \label{ass_model_dfo_2_g}
        \item $\|\bH_k\| \leq \kappa_H - 1$ for some \minorrev{fixed} $\kappa_H\geq 1$ (independent of $k$). \label{ass_model_dfo_2_H}
    \end{enumerate}
\end{assumption}

Then, our second-order algorithm is given in \algref{alg_basic_tr_dfo_2}.
Compared to the first-order algorithm \algref{alg_basic_tr_dfo}, we require fully quadratic models rather than fully linear models, assume the eigenstep decrease from the trust-region subproblem solver, and compare $\Delta_k$ with $\sigma^m_k$ \eqref{eq_optimality_model_2} to declare a step successful, instead of comparing just with $\|\bg_k\|$.
Additionally, we cap $\Delta_k$ at a maximum level $\Delta_{\max}$, a technicality needed \revision{to bound the error between the true and estimated criticality measures (\lemref{lem_dfo_2_crit_measures}). Ultimately, this comes from the need to control the size of the gradient error \eqref{eq_fully_quadratic_g}: reflecting the `overloaded' role of the trust-region radius, we have to cap the radius (and hence the size of the step $\bs_k$) solely because the radius also acts as the measure of model accuracy.} 

\begin{algorithm}[tb]
\begin{algorithmic}[1]
\Require Starting point $\bx_0\in\R^n$ and trust-region radius $\Delta_0>0$. Algorithm parameters: maximum trust-region radius $\Delta_{\max} \geq \Delta_0$, scaling factors $0 < \gammadec < 1 < \gammainc$, acceptance thresholds $0 < \eta_U \leq \eta_S < 1$, and criticality threshold $\mu_c > 0$.
\For{$k=0,1,2,\ldots$}
    \State Build a local quadratic model $m_k$ \eqref{eq_tr_model_dfo} satisfying \assref{ass_model_dfo_2}.
    \State Solve the trust-region subproblem \eqref{eq_trs} to get a step $\bs_k$ satisfying \assref{ass_eigenstep_decrease}.
    \State Evaluate $f(\bx_k+\bs_k)$ and calculate the ratio $\rho_k$ \eqref{eq_ratio_generic}.
    \If{$\rho_k \geq \eta_S$ and $\sigma^m_k \geq \mu_c \Delta_k$} 
        \State \textit{(Very successful iteration)} Set $\bx_{k+1}=\bx_k+\bs_k$ and $\Delta_{k+1}=\min(\gammainc\Delta_k, \Delta_{\max})$. 
    \ElsIf{$\eta_U \leq \rho_k < \eta_S$ and $\sigma^m_k \geq \mu_c \Delta_k$}
        \State \textit{(Successful iteration)} Set $\bx_{k+1}=\bx_k+\bs_k$ and $\Delta_{k+1}=\Delta_k$.
    \Else
        \State \textit{(Unsuccessful iteration)} Set $\bx_{k+1}=\bx_k$ and $\Delta_{k+1}=\gammadec\Delta_k$. 
    \EndIf 
\EndFor
\end{algorithmic}
\caption{Second-order \renaming{MBDFO} trust-region method for solving \eqref{eq_problem}.}
\label{alg_basic_tr_dfo_2}
\end{algorithm}

First, we show that fully quadratic models yield sufficiently good estimates of all the different criticality measures under consideration \minorrev{(see \eqref{eq_optimality_2} and \eqref{eq_optimality_model_2})}.

\begin{lemma} \label{lem_dfo_2_crit_measures}
    Suppose Assumptions~\minorrev{\ref{ass_smoothness_2}\ref{ass_smoothness_2_smooth}} and \ref{ass_model_dfo_2} hold.
    \minorrev{On iteration $k$ of \algref{alg_basic_tr_dfo_2}, we have} $\|\bg_k-\grad f(\bx_k)\| \leq \kappa_{\sigma} \Delta_k$, $|\tau^m_k - \tau_k| \leq \kappa_{\sigma} \Delta_k$ and
    $|\sigma_k-\sigma^m_k| \leq \kappa_{\sigma}\Delta_k$ all hold,
    where $\kappa_{\sigma} := \max(\kappamg \Delta_{\max}, \kappamh)$.
\end{lemma}
\begin{proof}
    Note that the trust-region updating mechanism in \algref{alg_basic_tr_dfo_2} and $\Delta_0 \leq \Delta_{\max}$ ensures $\Delta_k \leq \Delta_{\max}$ for all iterations $k$.
    
    First, we have $\|\bg_k - \grad f(\bx_k)\| \leq \kappamg \Delta_k^2 \leq \kappamg \Delta_{\max} \Delta_k \leq \kappa_{\sigma} \Delta_k$.
    Secondly, if $\bv$ is a normalized eigenvector corresponding to $\lambda_{\min}(\grad^2 f(\bx_k))$, then
    \begin{align}
        \lambda_{\min}(\bH_k) - \lambda_{\min}(\grad^2 f(\bx_k)) &\leq \bv^T [\bH_k - \grad^2 f(\bx_k)] \bv \leq \|\bH_k - \grad^2 f(\bx_k)\| \leq \kappamh \Delta_k \leq \kappa_{\sigma} \Delta_k,
    \end{align}
    and instead taking $\bv$ to be a normalized eigenvector for $\lambda_{\min}(\bH_k)$ we get the reverse result $\lambda_{\min}(\grad^2 f(\bx_k)) - \lambda_{\min}(\bH_k) \leq \kappa_{\sigma} \Delta_k$, and so all together we get $|\lambda_{\min}(\bH_k) - \lambda_{\min}(\grad^2 f(\bx_k))| \leq \kappa_{\sigma} \Delta_k$.
    The result then follows from the identity\footnote{To show this, suppose $a_{i^*} = \max_i a_i$. Then $\max_i a_i = a_{i^*} \leq b_{i^*} + |a_{i^*}-b_{i^*}| \leq \max_i b_i + \max_i |a_i-b_i|$ and a similar argument for the opposite direction.} $|\max_i a_i - \max_i b_i| \leq \max_i |a_i-b_i|$.
\end{proof}

We now state the second-order equivalent of \lemref{lem_very_successful_dfo}, which is split into two results for the two criticality measures $\|\bg_k\|$ and $\tau^m_k$.

\begin{lemma} \label{lem_very_successful_dfo_2_g}
    Suppose Assumptions~\minorrev{\ref{ass_smoothness_2}\ref{ass_smoothness_2_smooth}},  \minorrev{\ref{ass_eigenstep_decrease} and \ref{ass_model_dfo_2}}.
    \minorrev{If, on iteration $k$ of \algref{alg_basic_tr_dfo_2},} $\bg_k \neq \bm{0}$ and
    \begin{align}
        \Delta_k \leq \min\left(\frac{\kappa_s (1-\eta_S)}{2\kappamf \Delta_{\max}}, \frac{1}{\kappa_H}, \frac{1}{\mu_c}\right) \|\bg_k\|, 
    \end{align}
    then $\rho_k \geq \eta_S$ and $\sigma^m_k \geq \mu_c \Delta_k$ (i.e.~iteration $k$ is very successful).
\end{lemma}
\begin{proof}
    First, we have $\sigma^m_k \geq \|\bg_k\| \geq \mu_c \Delta_k$ by assumption on $\Delta_k$, so it remains to show $\rho_k \geq \eta_S$.
    From Assumptions~\ref{ass_eigenstep_decrease} and \ref{ass_model_dfo_2}, we have
    \begin{align}
        |\rho_k-1| \leq \frac{|f(\bx_k+\bs_k)-m_k(\bx_k+\bs_k)| + |f(\bx_k) - m_k(\bx_k)|}{|m_k(\bx_k)-m_k(\bx_k+\bs_k)|} \leq \frac{2\kappamf \Delta_k^3}{\kappa_s \|\bg_k\| \min\left(\Delta_k, \frac{\|\bg_k\|}{\kappa_H}\right)} = \frac{2\kappamf \Delta_k^3}{\kappa_s \|\bg_k\| \Delta_k},
    \end{align}
    since $\Delta_k \leq \frac{\|\bg_k\|}{\kappa_H}$ by assumption.
    Therefore we have
    \begin{align}
        |\rho_k-1| \leq \frac{2\kappamf \Delta_k^2}{\kappa_s \|\bg_k\|}  \leq \frac{2\kappamf \Delta_{\max}\Delta_k}{\kappa_s \|\bg_k\|},
    \end{align}
    and so we conclude $|\rho_k-1| \geq 1-\eta_S$ and hence $\rho_k \geq \eta_S$.
\end{proof}

\begin{lemma} \label{lem_very_successful_dfo_2_H}
    Suppose Assumptions~\minorrev{\ref{ass_smoothness_2}\ref{ass_smoothness_2_smooth}},  \minorrev{\ref{ass_eigenstep_decrease} and \ref{ass_model_dfo_2}} hold.
    \minorrev{If, on iteration $k$ of \algref{alg_basic_tr_dfo_2},} $\tau^m_k > 0$ and
    \begin{align}
        \Delta_k \leq \min\left(\frac{\kappa_s (1-\eta_S)}{2\kappamf}, \frac{1}{\mu_c} \right) \tau^m_k, 
    \end{align}
    then $\rho_k \geq \eta_S$ and $\sigma^m_k \geq \mu_c \Delta_k$ (i.e.~iteration $k$ is very successful).
\end{lemma}
\begin{proof}
    First, note that $\sigma^m_k \geq \tau^m_k \geq \mu_c \Delta_k$.
    Separately, from Assumptions~\ref{ass_eigenstep_decrease} and \ref{ass_model_dfo_2} we have
    \begin{align}
        |\rho_k-1| \leq \frac{|f(\bx_k+\bs_k)-m_k(\bx_k+\bs_k)| + |f(\bx_k) - m_k(\bx_k)|}{|m_k(\bx_k)-m_k(\bx_k+\bs_k)|} \leq \frac{2\kappamf \Delta_k^3}{\kappa_s \tau^m_k \Delta_k^2} = \frac{2\kappamf \Delta_k}{\kappa_s \tau^m_k} \leq 1-\eta_S,
    \end{align}
    and so $\rho_k \geq \eta_S$.
\end{proof}

We now get our usual result, that $\Delta_k$ is bounded away from zero provided we are not close to optimality.

\begin{lemma} \label{lem_delta_min_dfo_2}
    Suppose Assumptions~\minorrev{\ref{ass_smoothness_2}\ref{ass_smoothness_2_smooth}},  \minorrev{\ref{ass_eigenstep_decrease} and \ref{ass_model_dfo_2}} hold and we run \algref{alg_basic_tr_dfo_2}.
    If $\sigma_k \geq \epsilon$ for all $k=0,\ldots,K-1$, then
    \begin{align}
        \Delta_k \geq \Delta_{\min}(\epsilon) := \min\left(\Delta_0, \frac{\gammadec \epsilon}{\max\left(\frac{2\kappamf \Delta_{\max}}{\kappa_s (1-\eta_S)}, \kappa_H, \mu_c\right) + \kappa_{\sigma}}, \frac{\gammadec \epsilon}{\max\left(\frac{2\kappamf}{\kappa_s (1-\eta_S)}, \mu_c \right) + \kappa_{\sigma}}\right),
    \end{align}
    for all $k=0,\ldots,K$.
\end{lemma}
\begin{proof}
    We again proceed by induction.
    The result holds trivially for $k=0$, so suppose $\Delta_k \geq \Delta_{\min}(\epsilon)$ for some $k\in\{0,\ldots,K-1\}$.
    To find a contradiction assume that $\Delta_{k+1} < \Delta_{\min}(\epsilon)$.
    Then $\Delta_{k+1} < \Delta_k$, which by the mechanism for updating the trust-region radius means iteration $k$ was unsuccessful, and $\Delta_k = \gammadec^{-1} \Delta_{k+1} < \gammadec^{-1} \Delta_{\min}(\epsilon)$.
    
    At iteration $k$, since $\sigma_k \geq \epsilon$ we either have $\|\grad f(\bx_k)\| \geq \epsilon$ or $\tau_k \geq \epsilon$.
    If $\|\grad f(\bx_k)\| \geq \epsilon$, since iteration $k$ is unsuccessful, from \lemref{lem_very_successful_dfo_2_g} we have
    \begin{align}
        \Delta_k > \min\left(\frac{\kappa_s (1-\eta_S)}{2\kappamf \Delta_{\max}}, \frac{1}{\kappa_H}, \frac{1}{\mu_c}\right) \|\bg_k\|,
    \end{align}
    and so from \lemref{lem_dfo_2_crit_measures} we get
    \begin{align}
        \epsilon \leq \|\grad f(\bx_k)\| \leq \|\bg_k\| + \|\grad f(\bx_k) - \bg_k\| < \max\left(\frac{2\kappamf \Delta_{\max}}{\kappa_s (1-\eta_S)}, \kappa_H, \mu_c\right) \Delta_k + \kappa_{\sigma} \Delta_k,
    \end{align}
    which contradicts $\Delta_k < \gammadec^{-1} \Delta_{\min}(\epsilon)$.
    If instead $\tau_k \geq \epsilon$, from \lemref{lem_very_successful_dfo_2_H} we get
    \begin{align}
        \Delta_k > \min\left(\frac{\kappa_s (1-\eta_S)}{2\kappamf}, \frac{1}{\mu_c} \right) \tau^m_k,
    \end{align}
    which again from \lemref{lem_dfo_2_crit_measures} implies
    \begin{align}
        \epsilon \leq \tau_k \leq \tau^m_k + |\tau_k - \tau^m_k| < \max\left(\frac{2\kappamf}{\kappa_s (1-\eta_S)}, \mu_c \right) \Delta_k + \kappa_{\sigma} \Delta_k,
    \end{align}
    and we again get a contradiction with $\Delta_k < \gammadec^{-1} \Delta_{\min}(\epsilon)$.
\end{proof}

We now get our main second-order complexity result.

\begin{theorem} \label{thm_wcc_dfo_2}
    Suppose Assumptions~\ref{ass_smoothness_2},  \minorrev{\ref{ass_eigenstep_decrease} and \ref{ass_model_dfo_2}} hold and we run \algref{alg_basic_tr_dfo_2}.
    If $\sigma_k \geq \epsilon$ for all $k=0,\ldots,K-1$, then
    \begin{align}
        K \leq \frac{\log(\Delta_0 / \Delta_{\min}(\epsilon))}{\log(\gammadec^{-1})} + \left(1 + \frac{\log(\gammainc)}{\log(\gammadec^{-1})}\right) \revision{\frac{f(\bx_0)-\flow}{\eta_U \kappa_s \mu_c \min(1, \mu_c/\kappa_H) \Delta_{\min}(\epsilon)^3},}
    \end{align}
    where $\Delta_{\min}(\epsilon)$ is defined in \lemref{lem_delta_min_dfo_2}.
\end{theorem}
\begin{proof}
    We again partition the iterations $\{0,\ldots,K-1\}=\mathcal{S}\cup\mathcal{U}$, where $\mathcal{S}$ is the set of successful or very successful iterations and $\mathcal{U}$ is the set of unsuccessful iterations.
    From \lemref{lem_delta_min_dfo_2}, we have $\Delta_k \geq \Delta_{\min}(\epsilon)$ for all $k\in\{0,\ldots,K\}$.

    \revision{If $k\in\mathcal{S}$, then the criticality requirement $\sigma^m_k \geq \mu_c \Delta_k$ implies that either $\|\bg_k\| \geq \mu_c \Delta_k$ or $\tau^m_k \geq \mu_c \Delta_k$.
    First, if $\|\bg_k\| \geq \mu_c \Delta_k$, then Assumptions~\ref{ass_eigenstep_decrease} and \ref{ass_model_dfo_2}\ref{ass_model_dfo_2_H} give
    \begin{align}
        m_k(\bx_k) - m_k(\bx_k+\bs_k) \geq \kappa_s \mu_c \Delta_{\min}(\epsilon) \min\left(\Delta_{\min}(\epsilon), \frac{\mu_c \Delta_{\min}(\epsilon)}{\kappa_H}\right). 
    \end{align}
    Instead, if $\tau^m_k \geq \mu_c \Delta_k$, then  \assref{ass_eigenstep_decrease} gives
    \begin{align}
        m_k(\bx_k) - m_k(\bx_k+\bs_k) \geq \kappa_s \mu_c \Delta_{\min}(\epsilon)^3.
    \end{align}
    }
    In either case, we always have
    \begin{align}
        m_k(\bx_k) - m_k(\bx_k+\bs_k) \geq \revision{\kappa_s \mu_c \min(1, \mu_c/\kappa_H) \Delta_{\min}(\epsilon)^3.}
    \end{align}
    By the same reasoning as \eqref{eq_wcc_success} in the proof of first-order convergence (\thmref{thm_wcc_dfo}), we have
    \begin{align}
        f(\bx_0) - \flow \geq \eta_U \sum_{k\in\mathcal{S}} m_k(\bx_k) - m_k(\bx_k+\bs_k) \geq \revision{\eta_U \kappa_s \mu_c \min(1, \mu_c/\kappa_H) \Delta_{\min}(\epsilon)^3 |\mathcal{S}|},
    \end{align}
    or 
    \begin{align}
        |\mathcal{S}| \leq \revision{\frac{f(\bx_0)-\flow}{\eta_U \kappa_s \mu_c \min(1, \mu_c/\kappa_H) \Delta_{\min}(\epsilon)^3}.} \label{eq_wcc_s_dfo2}
    \end{align}
    The trust-region updating mechanism again gives us \eqref{eq_wcc_tmp2}, and the result follows from $K=|\mathcal{S}|+|\mathcal{U}|$.
\end{proof}

Recalling here the notation $\kappam := \max(\kappamf, \kappamg,\kappamh)$ (see \defref{def_fully_quadratic}) and noting that $\kappa_{\sigma}=\Theta(\kappam)$ and $\Delta_{\min}(\epsilon) = \Theta((\kappam+\kappa_H)^{-1} \epsilon)$ as $\epsilon\to 0^{+}$ and $\kappam,\kappa_H\to\infty$, we get the following second-order convergence result.

\begin{corollary} \label{cor_wcc_dfo_2}
    Suppose the assumptions of \thmref{thm_wcc_dfo_2} hold.
    If $k_{\epsilon}$ is the first iteration of \algref{alg_basic_tr_dfo_2} such that $\sigma_k < \epsilon$, then $k_{\epsilon} = \bigO(\revision{\kappa_H (\kappam + \kappa_H)^3} \epsilon^{-3})$.
    Hence $\liminf_{k\to\infty} \sigma_k = 0$.
\end{corollary}

\revision{This bound again matches $\bigO(\epsilon^{-3})$ from the derivative-based case (\thmref{thm_wcc_2}), but is again larger by a factor of $\kappa_H (\kappam + \kappa_H)^2$ arising from the same issue as the first-order case (\remref{rem_crit_options}). Again, although $\kappam$---and now also $\kappa_H$---will usually depend on $n$ explicitly (\thmref{thm_fully_quadratic}), they can be made dimension-independent (\corref{cor_fully_quadratic_wcc_improved}), but we will always need at least $\bigO(n^2)$ objective evaluations (to compensate for the loss of information from the first- and second-order oracles).}

\subsection{Termination} \label{sec_termination}
We end this section by briefly discussing how to terminate \renaming{MBDFO} algorithms in practice.

In all nonlinear optimization methods, in practice it is common to terminate after some fixed budget.
For derivative-based methods, this may be a maximum number of iterations $k$ or (particularly for large-scale problems) a maximum runtime.
In \renaming{MBDFO}, where the objective may be very expensive to evaluate, it is more common to have a maximum number of objective evaluations (since we may have multiple evaluations per iteration, depending on how exactly the model \eqref{eq_tr_model_dfo} is constructed).

Separately, it is important to have a termination condition based on optimality; i.e.~the algorithm terminates because it has reached a sufficiently accurate solution.
In the derivative-based case, we have direct access to optimality measures such as $\|\grad f(\bx_k)\|$, and so we can simply terminate if $\|\grad f(\bx_k)\| \leq \epsilon$ for some user-defined tolerance $\epsilon$.
Of course, this is not possible in the \renaming{MBDFO} case, as we no longer have $\grad f(\bx_k)$.

Although we have access to $\|\bg_k\|$ at each iteration, a small value of $\|\bg_k\|$ does not always imply a small value of $\|\grad f(\bx_k)\|$ (the  condition we actually want to be true).
If the model $m_k$ is fully linear and $\|\bg_k\| \geq \mu_c \Delta_k$, then we do have \revision{$\|\grad f(\bx_k) = \bigO(\|\bg_k\|)$ from \eqref{eq_crit_comparison}},
and so $\|\bg_k\| \leq \epsilon$ (again where $\epsilon$ is a user-defined tolerance) would imply $\|\grad f(\bx_k)\|=\bigO(\epsilon)$.
For this to be an appropriate termination condition, we require a certification that $m_k$ is fully linear. 

In practice, however, \renaming{MBDFO} methods usually terminate with `optimality' if $\Delta_k$ is sufficiently small.
This is theoretically justified, as we know that $\Delta_k$ never gets too small if we remain far from optimality (e.g.~\lemref{lem_delta_min_dfo}).
We also have the following result, to provide further comfort that we can terminate on small $\Delta_k$ (noting that all our main results above are related to convergence over a subsequence of iterations, not the full sequence of iterations).

\begin{lemma} \label{lem_delta_decreases}
    Suppose Assumptions~\ref{ass_smoothness_1}, \minorrev{\ref{ass_cauchy_decrease} and \ref{ass_model_dfo}} hold.
    Then the iterates of \minorrev{\algref{alg_basic_tr_dfo}} satisfy $\lim_{k\to\infty} \Delta_k = 0$.
\end{lemma}
\begin{proof}
    We first consider the case where \minorrev{there} are finitely many (very) successful iterations.
    \minorrev{This} means that iteration $k$ is unsuccessful for all $k$ sufficiently large, say $k\geq K$ for fixed $K$.
    That is, $\Delta_k = \gammadec^{k-K} \Delta_K$ for  for some $\Delta_K \leq \gammainc^K \Delta_0$.
    So, $\Delta_k \to 0$ as $k\to\infty$.
    
    Instead, suppose there are infinitely many (very) successful iterations.
    The below reasoning applies equally to either algorithm.
    For any such iteration $k$, from Assumptions~\ref{ass_model_dfo} and \ref{ass_cauchy_decrease}, \revision{the same reasoning as used to get \eqref{eq_wcc_dfo_tmp1a}---without replacing $\Delta_k$ with the lower bound $\Delta_{\min}(\epsilon)$---gives}
    \begin{align}
        f(\bx_0) - \flow \geq \eta_U \kappa_s \mu_c \min\left(1, \frac{\mu_c }{\kappa_H}\right)\sum_{k\in\mathcal{S}} \Delta_k^2,
    \end{align}
    where here $\mathcal{S}$ is the (infinite) set of all successful or very successful iterations.
    That is, we have $\sum_{k\in\mathcal{S}} \Delta_k^2 < \infty$.
    Since $\mathcal{S}$ is infinite, it must be that $\lim_{k\in \mathcal{S}} \Delta_k = 0$.
    For any $k\notin \mathcal{S}$ sufficiently large (i.e.~so that there was at least one (very) successful iteration before $k$), the trust-region updating mechanism guarantees $\Delta_k \leq \Delta_{s_k} \to 0$, where $s_k\in\mathcal{S}$ is the last (very) successful iteration before $k$, and so the result holds.
\end{proof}

The below result additionally tells us that terminating after the \emph{first} iteration when $\Delta_k$ drops below some termination threshold $\Delta_{\min}$ guarantees that the current iterate $\bx_k$ is a good solution.
Specifically, we get that $\|\grad f(\bx_k)\|=\bigO(\minorrev{(\kappam + \kappa_H)}\Delta_{\min})$.

\begin{lemma} \label{lem_terminate_good}
    Suppose Assumptions~\minorrev{\ref{ass_smoothness_1}\ref{ass_smoothness_1_smooth}}, \minorrev{\ref{ass_cauchy_decrease} and \ref{ass_model_dfo}} hold, and $\Delta_{\min} < \Delta_0$.
    If $\Delta_{k+1} < \Delta_{\min}$ occurs for the first time at the end of iteration $k$ of \minorrev{\algref{alg_basic_tr_dfo}}, then
    \begin{align}
        \|\grad f(\bx_k)\| \leq  \gammadec^{-1} \left[\max\left(\frac{2\kappamf}{\kappa_s (1-\eta_S)}, \kappa_H, \mu_c\right) + \kappamg\right] \Delta_{\min}.
    \end{align}
\end{lemma}
\begin{proof}
    Since $\Delta_{\min} < \Delta_0$, we know that $k\geq 0$.
    Moreover, since $\Delta_k$ is only decreased on unsuccessful iterations, we know that iteration $k$ is unsuccessful and $\Delta_k < \gammadec^{-1} \Delta_{\min}$.
    \minorrev{From \lemref{lem_delta_min_dfo}}, for iteration $k$ to be unsuccessful we must have
    \begin{align}
        \Delta_k > \min\left(\frac{\kappa_s (1-\eta_S)}{2\kappamf}, \frac{1}{\kappa_H}, \frac{1}{\mu_c}\right) \|\bg_k\|.
    \end{align}
    Combined with $\Delta_k < \gammadec^{-1}\Delta_{\min}$ we get
    \begin{align}
        \|\bg_k\| < \gammadec^{-1} \max\left(\frac{2\kappamf}{\kappa_s (1-\eta_S)}, \kappa_H, \mu_c\right) \Delta_{\min}.
    \end{align}
    Therefore, since $m_k$ is fully linear, we get
    \begin{align}
        \|\grad f(\bx_k)\| \leq \|\bg_k\| + \|\bg_k - \grad f(\bx_k)\| &< \gammadec^{-1} \max\left(\frac{2\kappamf}{\kappa_s (1-\eta_S)}, \kappa_H, \mu_c\right) \Delta_{\min} + \kappamg \Delta_k, \\
        &< \gammadec^{-1} \max\left(\frac{2\kappamf}{\kappa_s (1-\eta_S)}, \kappa_H, \mu_c\right) \Delta_{\min} + \kappamg \gammadec^{-1} \Delta_{\min},
    \end{align}
    \minorrev{as required}.
\end{proof}

Of course, \lemref{lem_delta_decreases} ensures that we will always terminate in finite time, for any choice of termination value $\Delta_{\min}$.

\revision{
\begin{remark}
    In \secref{sec_geom_incr_algo}, we will consider an algorithm where $m_k$ is not fully linear at every iteration.
    However, \lemref{lem_terminate_good} still holds in that case, because $\Delta_k$ is only ever decreased when $m_k$ is fully linear.
\end{remark}
}

\subsubsection*{Notes and References}
{\small The first convergence theory for an algorithm similar to \algref{alg_basic_tr_dfo} was given in \cite{Conn1997}, but the framework studied here using fully linear/\minorrev{fully} quadratic models was introduced later in \cite{Conn2009a} and expanded in \cite{Conn2009}.
The extension of these results to include worst-case complexity bounds was first given in the PhD thesis \cite{Judice2015}, with the first-order complexity later published in \cite{Garmanjani2016}. 
To the best of \minorrev{the author's} knowledge, the first use of the simplified algorithmic framework assuming fully linear models at all iterations (and checking criticality within the trust-region updating) was in \cite{Bandeira2014}.}

{\small Currently, one of the most active research directions in \renaming{MBDFO} is adapting this framework to be better suited to large-scale problems (i.e.~where the dimension $n$ is large, roughly $n\gg 100$). 
The most promising direction seems to be methods that construct models in low-dimensional subspaces at each iteration, where the subspaces may be selected randomly \cite{Cartis2023,Dzahini2024,Cartis2024} or using deterministic methods \cite{Xie2024,Zhang2025}.
}

\section{Interpolation Model Construction} \label{sec_model_construction}

In \minorrev{\secref{sec_dfotr},} we introduced a basic \renaming{MBDFO} trust-region method. 
To make this method concrete, we need a way to construct \minorrev{quadratic} models \eqref{eq_tr_model_dfo} at each iteration that are sufficiently accurate to guarantee convergence (i.e.~satisfying Assumptions~\ref{ass_model_dfo} or \ref{ass_model_dfo_2}).
In particular we focus on the fully linear/\minorrev{fully} quadratic requirements Assumptions~\ref{ass_model_dfo}\ref{ass_model_dfo_g} and \ref{ass_model_dfo_2}\ref{ass_model_dfo_2_g}.
Of course, we must be able to form these models using only zeroth-order information about the objective $f$, which we will achieve by constructing models that interpolate $f$ over carefully chosen collections of points.

We first consider generating good collections of points for linear/quadratic interpolation to guarantee fully linear/\minorrev{fully} quadratic models.
\revision{We derive generic error bounds in terms of the linear algebra associated with the interpolation system, and consider special cases where points are sampled in a regular way around the current iterate (e.g.~perturbations along coordinate axes). These special cases are useful to understand the theoretical capabilities of MBDFO methods, but this process of choosing a structured interpolation set is also practically useful.
As outlined in \secref{sec_interp_practical_ibcdfo}, in the IBCDFO software collection (see \secref{sec_conclusion}), an interpolation set is built by appending suitable points one-by-one from the full history of the solver; the error bounds we derive are used to ensure `suitability' and to complete the set of points if not enough historical points are suitable.}


To simplify the derivation of error bounds, we will use the following result, which says that fully linear/\minorrev{fully} quadratic models can be achieved by ensuring sufficiently accurate approximations to a Taylor series.

\begin{lemma} \label{lem_fully_linear_quadratic_from_taylor_approx}
    Suppose we have a quadratic model $m:\R^n\to\R$ approximating $f:\R^n\to\R$ near a point $\bx\in\R^n$,
    \begin{align}
        f(\by) \approx m(\by) := c + \bg^T (\by-\bx) + \frac{1}{2} (\by-\bx)^T \bH (\by-\bx).
    \end{align}
    Then, for any $\Delta>0$, the following results hold.
    \begin{enumerate}[label=(\alph*)]
        \item If \minorrev{\assref{ass_smoothness_1}\ref{ass_smoothness_1_smooth}} holds and there exists $\kappa>0$ such that
        \begin{align}
            |m(\by) - f(\bx) - \grad f(\bx)^T (\by-\bx)| \leq \kappa \Delta^2,
        \end{align}
        for all $\by\in B(\bx,\Delta)$, then $m$ is fully linear in $B(\bx,\Delta)$ with constants (c.f.~\eqref{eq_fully_linear})
        \begin{align}
            \kappamf = \kappa + \frac{\Lgrad}{2}, \qquad \text{and} \qquad \kappamg = 2\kappa + \Lgrad + 2\kappa_H,
        \end{align}
        where $\kappa_H$ is any upper bound on $\|\bH\|$.
        \item If \minorrev{\assref{ass_smoothness_2}\ref{ass_smoothness_2_smooth}} holds and there exists $\kappa>0$ such that
        \begin{align}
            \left|m(\by) - f(\bx) - \grad f(\bx)^T (\by-\bx) - \frac{1}{2}(\by-\bx)^T \grad^2 f(\bx) (\by-\bx)\right| \leq \kappa \Delta^3, \label{eq_model_quad_taylor_diff}
        \end{align}
        for all $\by\in B(\bx,\Delta)$, then $m$ is fully quadratic in $B(\bx,\Delta)$ with constants (c.f.~\eqref{eq_fully_quadratic})
        \begin{align}
            \kappamf = \kappa + \frac{\Lhess}{6}, \qquad \kappamg = 34\kappa + \frac{\Lhess}{2}, \qquad \text{and} \qquad \kappamh = 24 \kappa + \Lhess.
        \end{align}
    \end{enumerate}
\end{lemma}
\begin{proof}
    We will use the technical results in \appref{app_technical_results}, which say that if two linear/quadratic functions are close in $B(\bx,\Delta)$, then their coefficients are also close.
    
    (a) Fix any $\by\in B(\bx,\Delta)$, and use \lemref{lem_lsmooth} to get
    \begin{align}
        |m(\by) - f(\by)| &\leq |m(\by) - f(\bx) - \grad f(\bx)^T (\by-\bx)| + |f(\by) - f(\bx) - \grad f(\bx)^T (\by-\bx)|, \\
        &\leq \kappa \Delta^2 + \frac{\Lgrad}{2} \|\by-\bx\|^2, \label{eq_model_taylor_tmp1}
    \end{align}
    and we get the value for $\kappamf$ after $\|\by-\bx\| \leq \Delta$.
    Next, for any $\by\in B(\bx,\Delta)$ we have 
    \begin{align}
        |c + \bg^T (\by-\bx) - f(\bx) - \grad f(\bx)^T (\by-\bx)| \leq \kappa \Delta^2 + \frac{1}{2}\|\bH\| \: \|\by-\bx\|^2 \leq \left(\kappa + \frac{1}{2}\kappa_H\right)\Delta^2,
    \end{align}
    and so \lemref{lem_difference_of_linear} gives $\|\bg-\grad f(\bx)\| \leq (2\kappa + \kappa_H) \Delta$.
    So, for an arbitrary $\by\in B(\bx,\Delta)$, we use $\grad m(\by) = \bg + \bH(\by-\bx)$ to get
    \begin{align}
        \|\grad m(\by) - \grad f(\by)\| &\leq \|\bH(\by-\bx)\| + \|\bg - \grad f(\bx)\| + \|\grad f(\by) - \grad f(\bx)\|, \\
        &\leq \kappa_H \|\by-\bx\| + (2\kappa + \kappa_H)\Delta + \Lgrad\|\by-\bx\|,
    \end{align}
    and we get the value of $\kappamg$.

    (b) First, to get $\kappamf$, we follow the same reasoning as for \eqref{eq_model_taylor_tmp1} but replacing \lemref{lem_lsmooth} with \eqref{eq_lh_smooth} to get
    \begin{align}
        |m(\by) - f(\by)| \leq \kappa \Delta^3 + \frac{\Lhess}{6} \|\by-\bx\|^3,
    \end{align}
    for any $\by\in B(\bx,\Delta)$.
    Next, we use \lemref{lem_difference_of_quadratic} on \eqref{eq_model_quad_taylor_diff} to get
    \begin{align}
        \|\bg - \grad f(\bx)\| \leq 10 \kappa \Delta^2, \qquad \text{and} \qquad \|\bH - \grad^2 f(\bx)\| \leq 24 \kappa \Delta.
    \end{align}
    So, for any $\by\in B(\bx,\Delta)$ we have
    \begin{align}
        \|\grad m(\by) - \grad f(\by)\| &\leq \|\bg - \grad f(\bx)\| + \|\bH(\by-\bx) - \grad^2 f(\bx) (\by-\bx)\| \nonumber \\
        &\qquad\qquad\qquad\qquad + \|\grad f(\by) - \grad f(\bx) - \grad^2 f(\bx) (\by-\bx)\|, \\
        &\leq 10 \kappa \Delta^2 + 24 \kappa \Delta \|\by-\bx\| + \frac{\Lhess}{2} \|\by-\bx\|^2,
    \end{align}
    where the last inequality uses \eqref{eq_lh_smooth_gradf}, and we get the value for $\kappamg$.
    Lastly, we have
    \begin{align}
        \|\grad^2 m(\by) - \grad^2 f(\by)\| \leq \|\bH - \grad^2 f(\bx)\| + \|\grad^2 f(\by) - \grad^2 f(\bx)\| \leq 24 \kappa \Delta + \Lhess\|\by-\bx\|,
    \end{align}
    and we get $\kappamh$.
\end{proof}

Our goal is now to show that linear/quadratic interpolation models for $f$ can provide sufficiently good approximations to a Taylor series for $f$, from which \lemref{lem_fully_linear_quadratic_from_taylor_approx} tells us that the model is fully linear/\minorrev{fully} quadratic, as needed for our algorithms.

\subsection{Linear Interpolation} \label{sec_linear_interp}

In the first instance, we will try to construct fully linear models (\defref{def_fully_linear}) using linear interpolation.
Here, for our objective function $f:\R^n\to\R$, we will build a local linear model around a base point $\bx\in\R^n$ (which will be $\bx=\bx_k$ at the $k$-th iteration of our algorithm),
\begin{align}
    f(\by) \approx m(\by) := c + \bg^T (\by-\bx), \label{eq_linear_model_generic}
\end{align}
for some $c\in\R$ and $\bg\in\R^n$; i.e.~taking $\bH_k=\bm{0}$ in \eqref{eq_tr_model_dfo}.\footnote{This then satisfies \assref{ass_model_dfo}\ref{ass_model_dfo_H} with $\kappa_H=1$.}
We choose the $p := n+1$ unknowns ($c$ and $\bg$) to interpolate known values of $f$ at $p$ points, $\by_1,\ldots,\by_p$.
That is, we pick $c$ and $\bg$ such that $f(\by_i) = m(\by_i)$ for all $i=1,\ldots,p$.
This reduces to solving the following $p\times p$ linear system:
\begin{align}
    \underbrace{\begin{bmatrix} 1 & (\by_1-\bx)^T \\ \vdots & \vdots \\ 1 & (\by_p-\bx)^T \end{bmatrix}}_{=: \bM} \begin{bmatrix} c \\ \bg \end{bmatrix} = \begin{bmatrix} f(\by_1) \\ \vdots \\ f(\by_p) \end{bmatrix}. \label{eq_linear_interp_system}
\end{align}

We shall see that, for the resulting $c$ and $\bg$ to give a fully linear model (\defref{def_fully_linear}) in a ball around the base point, $B(\bx,\Delta)$---where at iteration $k$ we will usually have $\bx=\bx_k$, our current iterate and $\Delta=\Delta_k$, our current trust-region radius---the interpolation points $\by_1,\ldots,\by_p$ will have to be close to the base point, $\|\by_i-\bx\| \lesssim \Delta$ (i.e.~$\|\by_i-\bx\| \leq c \Delta$ for a constant $c$ not much larger than 1).
However, as we saw in \lemref{lem_delta_decreases}, the trust-region radius $\Delta_k \to 0^{+}$ as $k\to\infty$, and so the matrix $\bM$ in \eqref{eq_linear_interp_system} becomes increasingly ill-conditioned as the algorithm progresses (as each column of $\bM$ except the first approaches zero).

To avoid this ill-conditioning, we can rescale the second through last columns\footnote{This is a form of \emph{equilibration} \cite[Chapter 3.5.2]{Golub1996}.} of the linear system \eqref{eq_linear_interp_system} by $\Delta^{-1}$ to replace $\by_i-\bx$ with $\hat{\bs}_i := (\by_i-\bx)/\Delta$, to get
\begin{align}
    \underbrace{\begin{bmatrix} 1 & \hat{\bs}_1^T \\ \vdots & \vdots \\ 1 & \hat{\bs}_p^T \end{bmatrix}}_{=: \hat{\bM}} \begin{bmatrix} c \\ \hat{\bg} \end{bmatrix} = \begin{bmatrix} f(\by_1) \\ \vdots \\ f(\by_p) \end{bmatrix}, \qquad \text{where $\hat{\bg} := \Delta \: \bg$.} \label{eq_linear_interp_system_scaled}
\end{align}
Now, we expect that all $\|\hat{\bs}_i\| \lesssim 1$, and we avoid ill-conditioning resulting from $\Delta_k\to 0^{+}$.
This rescaling does not \emph{guarantee} that the new linear system $\hat{\bM}$ \eqref{eq_linear_interp_system_scaled} is well-conditioned---we will ensure this by judicious choice of the $\hat{\bs}_i$---but it does ensure that the magnitude of $\Delta$ is not a source of ill-conditioning.

In this setting, since all $\|\hat{\bs}_i\| \lesssim 1$, we know that $\|\hat{\bM}\|$ is not too large, since $\|\hat{\bM}\| \leq \|\hat{\bM}\|_F \lesssim \sqrt{2p}$.
Hence if $\hat{\bM}$ is invertible we have $\kappa(\hat{\bM}) \lesssim \sqrt{2p} \|\hat{\bM}^{-1}\|$, \revision{where $\kappa(\cdot)$ denotes the (2-norm) matrix condition number}, so any ill-conditioning in $\hat{\bM}$ is entirely driven by the size of $\|\hat{\bM}^{-1}\|$.
\revision{This is reflected in the following interpolation error bound, where we assume $\|\hat{\bs}_i\| \leq \beta$ for some $\beta>0$, and both $\beta$ and $\|\hat{\bM}^{-1}\|$ appear explicitly. }

\begin{theorem} \label{thm_fully_linear}
    Suppose $f$ satisfies \minorrev{\assref{ass_smoothness_1}\ref{ass_smoothness_1_smooth}} and we construct a linear model \eqref{eq_linear_model_generic} for $f$ by solving \eqref{eq_linear_interp_system_scaled}, where we assume $\hat{\bM}$ is invertible.
    If $\|\by_i-\bx\| \leq \beta \Delta$ for some $\beta>0$ and all $i=1,\ldots,p$ (where $p=n+1$ is the number of interpolation points), then the model is fully linear in $B(\bx,\Delta)$ with constants
    \begin{align}
        \kappamf = \frac{\Lgrad}{2} (1+\sqrt{n}) \beta^2 \|\hat{\bM}^{-1}\|_{\infty} + \frac{\Lgrad}{2}, \quad \text{and} \quad \kappamg = 2\kappamf. \label{eq_thm_fully_linear_constants} 
    \end{align}
\end{theorem}
\begin{proof}
    For $i\in\{1,\ldots,p\}$, the $i$th row of \eqref{eq_linear_interp_system_scaled} gives $m(\by_i) = f(\by_i)$, and so
    \begin{align}
        |m(\by_i) - f(\bx) - \grad f(\bx)^T (\by_i-\bx)| \leq \frac{\Lgrad}{2} \|\by_i-\bx\|^2 \leq \frac{\Lgrad}{2} \beta^2 \Delta^2, \label{eq_linear_interp_error_tmp1}
    \end{align}
    using \lemref{lem_lsmooth}.
    Now, consider some arbitrary point $\by\in B(\bx,\Delta)$ and define $\hat{\bs}:=(\by-\bx)/\Delta$, giving $\|\hat{\bs}\| \leq 1$.
    Since $\hat{\bM}$ is invertible, so too is $\hat{\bM}^T$, and there exists a unique solution $\bv\in\R^p$ to $ \hat{\bM}^T \bv = \begin{bmatrix} 1 \\ \hat{\bs} \end{bmatrix}$.
    Multiplying the last $n$ rows of this equation by $\Delta$, we get
    \begin{align}
        \begin{bmatrix} 1 \\ \by-\bx \end{bmatrix} = \sum_{i=1}^{p} v_i \begin{bmatrix} 1 \\ \by_i-\bx \end{bmatrix}, \qquad \text{where} \qquad \bv = \hat{\bM}^{-T} \begin{bmatrix} 1 \\ \hat{\bs} \end{bmatrix}. \label{eq_linear_interp_error_tmp3}
    \end{align}

    Considering the value of the model at this point $\by$, we get
    \begin{align}
        |m(\by) - f(\bx) - \grad f(\bx)^T (\by-\bx)| &= \left| \begin{bmatrix} c \\ \bg \end{bmatrix}^T \begin{bmatrix} 1 \\ \by-\bx \end{bmatrix} - \begin{bmatrix} f(\bx) \\ \grad f(\bx) \end{bmatrix}^T \begin{bmatrix} 1 \\ \by-\bx \end{bmatrix}\right|, \\
        &= \left|\sum_{i=1}^{p} v_i \left(\begin{bmatrix} c \\ \bg \end{bmatrix}^T \begin{bmatrix} 1 \\ \by_i-\bx \end{bmatrix} - \begin{bmatrix} f(\bx) \\ \grad f(\bx) \end{bmatrix}^T \begin{bmatrix} 1 \\ \by_i-\bx \end{bmatrix}\right)\right|, \\
        &= \left|\sum_{i=1}^{p} v_i \left(m(\by_i) - f(\bx) - \grad f(\bx)^T (\by_i-\bx)\right)\right|, \label{eq_linear_interp_error_tmp2a} \\
        &\leq \frac{\Lgrad}{2} \beta^2 \Delta^2 \|\bv\|_1, \label{eq_linear_interp_error_tmp2}
    \end{align}
    where the last line comes from the triangle inequality and \eqref{eq_linear_interp_error_tmp1}.
    From $\|\bv\|_1 \leq \|\hat{\bM}^{-T}\|_1 (1 + \|\hat{\bs}\|_1) \leq \|\hat{\bM}^{-1}\|_{\infty} (1 + \sqrt{n})$, a consequence of the identity $\|\bA^T\|_{1} = \|\bA\|_{\infty}$ and $\|\hat{\bs}\|_1 \leq \sqrt{n} \|\hat{\bs}\| \leq \sqrt{n}$, we get
    \begin{align}
        |m(\by) - f(\bx) - \grad f(\bx)^T (\by-\bx)| \leq \frac{\Lgrad}{2} (1+\sqrt{n}) \beta^2 \|\hat{\bM}^{-1}\|_{\infty} \Delta^2. \label{eq_linear_interp_error_tmp4}
    \end{align}
    The result then follows from \lemref{lem_fully_linear_quadratic_from_taylor_approx}(a) with $\kappa_H=0$.
\end{proof}

\revision{\begin{remark}
    \thmref{thm_fully_linear} gives the error bounds in terms of $\|\hat{\bM}^{-1}\|_{\infty}$, as opposed to the 2-norm which occurs in the corresponding bounds in earlier works, e.g.~\cite{Conn2009}.
    This does not significantly change our results, but the proof approach where we control the model error relative to a Taylor series via \eqref{eq_linear_interp_error_tmp2} generalizes to interpolation theory in constrained regions (\secref{sec_convex_constraints}) and naturally yields bounds in terms of $\|\hat{\bM}^{-1}\|_{\infty}$.
\end{remark}}

All that remains is to find a set of interpolation points so that $\|\by_i-\bx\| \leq \beta \Delta$ and $\|\hat{\bM}^{-1}\|_{\infty}$ is not too large.
The simplest way to do this is to use $\bx$ and $n$ perturbations around $\bx$ of size $\Delta$.
\revision{The most natural way to do this is via coordinate perturbations around $\bx$, i.e.~taking our interpolation set to be 
\begin{align}
    \{\bx_k, \bx_k + \Delta_k \bm{e}_1, \ldots, \bx_k + \Delta_k \bm{e}_n\}. \label{eq_linear_interp_coords}
\end{align}
}
This gives us the following worst-case complexity bound for \algref{alg_basic_tr_dfo}, building on the general result \corref{cor_wcc_dfo}.

\begin{corollary} \label{cor_fully_linear_wcc}
    Under the assumptions of \thmref{thm_wcc_dfo}, if the \minorrev{model $m_k$ \eqref{eq_tr_model_dfo}} at each iteration \minorrev{of} \algref{alg_basic_tr_dfo} is generated by linear interpolation to the points \revision{in \eqref{eq_linear_interp_coords}}, then the iterates of \algref{alg_basic_tr_dfo} achieve $\|\grad f(\bx_k)\| < \epsilon$ for the first time after at most $\bigO(n \epsilon^{-2})$ iterations and $\bigO(n^2 \epsilon^{-2})$ objective evaluations.
\end{corollary}
\begin{proof}
    \revision{Given the interpolation points \eqref{eq_linear_interp_coords}, we have
    \begin{align}
     \hat{\bM} = \begin{bmatrix} 1 & \bm{0}^T \\ \bm{e} & \bI \end{bmatrix}, \qquad \text{and so} \qquad \hat{\bM}^{-1} = \begin{bmatrix} 1 & \bm{0}^T \\ -\bm{e} & \bI \end{bmatrix}, \label{eq_linear_sampling}
    \end{align}
    so $\|\hat{\bM}^{-1}\|_{\infty} = 2$ and hence the fully linear constants from \thmref{thm_fully_linear} satisfy $\kappamf, \kappamg = \bigO(\sqrt{n} \Lgrad)$.}
    We then apply \corref{cor_wcc_dfo} with $\kappa_H=1$ in \assref{ass_model_dfo}\ref{ass_model_dfo_H} to get the iteration bound.
    The objective evaluation bound comes from observing that our interpolation mechanism requires evaluating $f$ at no more than $n+1$ points per iteration.
\end{proof}

\revision{Using the interpolation set \eqref{eq_linear_interp_coords}} is equivalent to building a linear Taylor model using finite differencing \eqref{eq_fin_diff} with the dynamically adjusted stepsize $\Delta_k$ to estimate $\grad f(\bx_k)$. 
\revision{However, it is this dynamic adjustment that is at the core of MBDFO, reflecting that we need to adjust our interpolation set depending on the required accuracy in the model (controlled by $\Delta_k$).}

\revision{The dimension dependency in \corref{cor_fully_linear_wcc} can be improved by considering coordinate perturbations of size $\Delta_k/\sqrt{n}$.

\begin{corollary} \label{cor_fully_linear_wcc_improved}
    Under the assumptions of \thmref{thm_wcc_dfo}, if the \minorrev{model $m_k$ \eqref{eq_tr_model_dfo}} at each iteration \minorrev{of} \algref{alg_basic_tr_dfo} is generated by linear interpolation to the points $\{\bx_k, \bx_k + \Delta_k \bm{e}_1 / \sqrt{n}, \ldots, \bx_k + \Delta_k \bm{e}_n / \sqrt{n}\}$, then the iterates of \algref{alg_basic_tr_dfo} achieve $\|\grad f(\bx_k)\| < \epsilon$ for the first time after at most $\bigO(\epsilon^{-2})$ iterations and $\bigO(n \epsilon^{-2})$ objective evaluations.
\end{corollary}
\begin{proof}
    For this interpolation set, we have $\|\hat{\bM}^{-1}\|_{\infty} = \bigO(\sqrt{n})$ and $\beta=1/\sqrt{n}$, so the model is fully linear with $\kappamf,\kappamg=\bigO(\Lgrad)$ (i.e.~no explicit dependence on $n$) by \thmref{thm_fully_linear}.
    The remainder of the proof is identical to \corref{cor_fully_linear_wcc}.
\end{proof}

So, at least in terms of the explicit dependence on $n$, \corref{cor_fully_linear_wcc_improved} gives an iteration bound that matches the derivative-based case, with a $\bigO(n)$ evaluation bound that comes from the fact that we need $n$ objective evaluations to match the same amount of problem information as one gradient evaluation (i.e.~$n$ real numbers).

However, taking an even smaller perturbation size cannot reduce $\kappamf$ to a smaller value (e.g.~$\kappamf=\bigO(1/n)$) to further improve the worst-case complexity, since $\kappamf$ always has a component of size $\bigO(\Lgrad)$ from \assref{ass_smoothness_1} that is independent of the model construction.
}


\subsection{Quadratic Interpolation} \label{sec_quadratic_interp}
Now that we have a method for constructing fully linear models, which can be used to achieve convergence to first-order stationary points, we now consider the construction of fully quadratic models (\defref{def_fully_quadratic}) using quadratic interpolation, which will allow us to find second-order stationary points using \algref{alg_basic_tr_dfo_2}.

In this case, we have an objective $f:\R^n\to\R$ and we now build a local quadratic model around a base point $\bx\in\R^n$ (which again will usually be $\bx=\bx_k$ in the $k$-th iteration of our algorithm),
\begin{align}
    f(\by) \approx m(\by) := c + \bg^T (\by-\bx) + \frac{1}{2} (\by-\bx)^T \bH (\by-\bx), \label{eq_quadratic_model_generic}
\end{align}
for some $c\in\R$, $\bg\in\R^n$ and $\bH\in\R^{n\times n}$ symmetric.
This model has $p := 1 + n + n(n+1)/2 = (n+1)(n+2)/2$ degrees of freedom\footnote{Note we always take $p$ to be the number of interpolation points, so the value of $p$ here is different to the linear case in \secref{sec_linear_interp}.}, and so we again construct $m$ to interpolate $f$ at $p$ points $\by_1,\ldots,\by_p$.

To write this interpolation problem efficiently, we introduce the natural basis for quadratic functions $\R^n\to\R$, which we denote $\bphi:\R^n\to\R^p$, defined as
\begin{align}
    \bphi(\bx) := \left[1, x_1, \ldots, x_n, \minorrev{\tfrac{1}{2}\,x_1^2}, x_1 x_2, \ldots, x_1 x_n, \minorrev{\tfrac{1}{2}\,x_2^2}, x_2 x_3, \ldots, x_{n-1} x_n, \minorrev{\tfrac{1}{2}\,x_n^2}\right]^T. \label{eq_natural_quadratic_basis}
\end{align}
With this basis, and noting $\bH=\bH^T$, the model \eqref{eq_quadratic_model_generic} can be written as
\begin{align}
    m(\by) = \bphi(\by-\bx)^T \begin{bmatrix} c \\ \bg \\ \upper(\bH) \end{bmatrix}, \label{eq_quad_model_natural_basis}
\end{align}
where $\upper(\bH)\in\R^{n(n+1)/2}$ contains the upper triangular elements of $\bH$, arranged as
\begin{align}
    \upper(\bH) = \left[H_{1,1}, H_{1,2}, \ldots, H_{1,n}, H_{2,2}, H_{2,3}, \ldots, H_{n-1,n}, H_{n,n}\right]^T.
\end{align}
With these definitions, the interpolation problem is given by the $p\times p$ linear system
\begin{align}
    \underbrace{\begin{bmatrix} \bphi(\by_1-\bx)^T \\ \vdots \\ \bphi(\by_p-\bx)^T \end{bmatrix}}_{=: \bQ} \begin{bmatrix} c \\ \bg \\ \upper(\bH) \end{bmatrix} = \begin{bmatrix} f(\by_1) \\ \vdots \\ f(\by_p) \end{bmatrix}. \label{eq_quadratic_interp_system}
\end{align}
Just like the linear case, the matrix $\bQ$ can become very ill-conditioned if $\Delta$, a typical value for $\|\by_i-\bx\|$, gets very small.
To address this, we define $\hat{\bs}_i := (\by_i-\bx)/\Delta$ and instead solve
\begin{align}
    \underbrace{\begin{bmatrix} \bphi(\hat{\bs}_1)^T \\ \vdots \\ \bphi(\hat{\bs}_p)^T \end{bmatrix}}_{=: \hat{\bQ}} \begin{bmatrix} c \\ \hat{\bg} \\ \upper(\hat{\bH}) \end{bmatrix} = \begin{bmatrix} f(\by_1) \\ \vdots \\ f(\by_p) \end{bmatrix}, \label{eq_quadratic_interp_system_scaled}
\end{align}
where $\hat{\bg} := \Delta \: \bg$ and $\hat{\bH} := \Delta^2 \: \bH$.

As before, since $\|\hat{\bs}_i\|\lesssim 1$, we know that $\|\hat{\bQ}\|$ is not too large, and so---provided $\hat{\bQ}$ is invertible---the only potential source of ill-conditioning is if $\hat{\bQ}^{-1}$ is large.
As such, the size of $\hat{\bQ}^{-1}$ determines the fully quadratic error bounds.

\begin{theorem} \label{thm_fully_quadratic}
    Suppose $f$ satisfies \minorrev{\assref{ass_smoothness_2}\ref{ass_smoothness_2_smooth}} and we construct a quadratic model $m$ \eqref{eq_quadratic_model_generic} for $f$ by solving \eqref{eq_quadratic_interp_system_scaled}, where we assume $\hat{\bQ}$ is invertible.
    If $\|\by_i-\bx\| \leq \beta \Delta$ for some $\beta>0$ and all $i=1,\ldots,p$ (where $p=(n+1)(n+2)/2$ is the number of interpolation points), then the model is fully quadratic in $B(\bx,\Delta)$ with constants 
    \begin{subequations}
    \begin{align}
        \kappamf &= \frac{1}{6} \Lhess \left(n + \sqrt{n} + \frac{3}{2}\right) \beta^3 \|\hat{\bQ}^{-1}\|_{\infty} + \frac{\Lhess}{6}, \\
        \kappamg &= \frac{17}{3} \Lhess \left(n + \sqrt{n} + \frac{3}{2}\right) \beta^3 \|\hat{\bQ}^{-1}\|_{\infty} + \frac{\Lhess}{2}, \qquad \text{and} \\
        \kappamh &= 4 \Lhess \left(n + \sqrt{n} + \frac{3}{2}\right) \beta^3 \|\hat{\bQ}^{-1}\|_{\infty} + \Lhess.
    \end{align}
    \end{subequations}
\end{theorem}
\begin{proof}
    For ease of notation, define $T_{f,\bx,2}(\by) := f(\bx) + \grad f(\bx)^T (\by-\bx) + \frac{1}{2}(\by-\bx)^T \grad^2 f(\bx)(\by-\bx)$ as the second-order Taylor series for $f$ at $\bx$.
    For $i\in\{1,\ldots,p\}$, the $i$th row of \eqref{eq_quadratic_interp_system_scaled} gives $f(\by_i) = m(\by_i)$ and so
    \begin{align}
        \left|m(\by_i) - T_{f,\bx,2}(\by_i)\right| \leq \frac{\Lhess}{6} \|\by_i-\bx\|^3 \leq \frac{\Lhess}{6} \beta^3 \Delta^3, \label{eq_fully_quadratic_tmp1}
    \end{align}
    using \eqref{eq_lh_smooth}.
    Now we choose $\by\in B(\bx,\Delta)$ and define $\hat{\bs} := (\by-\bx)/\Delta$, giving $\|\hat{\bs}\| \leq 1$.
    Since $\hat{\bQ}$ is invertible, so is $\hat{\bQ}^T$, and so we let $\bv\in\R^p$ be the unique solution to $\hat{\bQ}^T \bv = \bphi(\hat{\bs})$.
    Scaling rows 2 through \minorrev{$n+1$} of this equation by $\Delta$ and rows \minorrev{$n+2$} through $p$ by $\Delta^2$, we get
    \begin{align}
        \bphi(\by-\bx) = \sum_{i=1}^{p} v_i \bphi(\by_i-\bx), \qquad \text{where} \qquad \bv = \hat{\bQ}^{-T} \bphi(\hat{\bs}). \label{eq_fully_quadratic_tmp2}
    \end{align}
    We then get
    \begin{align}
        \left|m(\by) - T_{f,\bx,2}(\by)\right| &= \left|\bphi(\by-\bx)^T \begin{bmatrix} c \\ \bg \\ \upper(\bH) \end{bmatrix} - \bphi(\by-\bx)^T \begin{bmatrix} f(\bx) \\ \grad f(\bx) \\ \upper(\grad^2 f(\bx))\end{bmatrix}\right|, \\
        &= \left|\sum_{i=1}^{p} v_i \left(\bphi(\by_i-\bx)^T \begin{bmatrix} c \\ \bg \\ \upper(\bH) \end{bmatrix} - \bphi(\by_i-\bx)^T \begin{bmatrix} f(\bx) \\ \grad f(\bx) \\ \upper(\grad^2 f(\bx))\end{bmatrix}\right)\right|, \\
        &= \left|\sum_{i=1}^{p} v_i (m(\by_i) - T_{f,\bx,2}(\by_i))\right|, \\
        &\leq \frac{\Lhess}{6} \beta^3 \Delta^3 \|\bv\|_1, \label{eq_fully_quadratic_tmp3}
    \end{align}
    where the inequality comes from \eqref{eq_fully_quadratic_tmp1}.
    We then use $\|\bv\|_1 \leq \|\hat{\bQ}^{-T}\|_1 \|\bphi(\hat{\bs})\|_1 = \|\hat{\bQ}^{-1}\|_{\infty} \|\bphi(\hat{\bs})\|_1$, and it remains to bound $\|\bphi(\hat{\bs})\|_1$.
    To do this, we use the definition of $\bphi$ \eqref{eq_natural_quadratic_basis} to write
    \begin{align}
        \|\bphi(\hat{\bs})\|_1 = 1 + \|\hat{\bs}\|_1 + \frac{1}{2}\hat{s}_1^2 + \cdots + \frac{1}{2}\hat{s}_n^2 + \sum_{i=1}^{n} \sum_{j=i+1}^{n} |\hat{s}_i \hat{s}_j|.
    \end{align}
    Using Young's inequality, $|ab| \leq \frac{|a|^2 + |b|^2}{2} = \frac{a^2 + b^2}{2}$ for $a,b\in\R$, we get
    \begin{align}
        \|\bphi(\hat{\bs})\|_1 &\leq 1 + \|\hat{\bs}\|_1 + \frac{1}{2} \|\hat{\bs}\|^2 + \frac{1}{2}\sum_{i=1}^{n} \sum_{j=i+1}^{n} (\hat{s}_i^2 + \hat{s}_j^2), \\
        &\leq 1 + \|\hat{\bs}\|_1 + \frac{1}{2} \|\hat{\bs}\|^2 + \frac{1}{2}\sum_{i=1}^{n} \sum_{j=1}^{n} (\hat{s}_i^2 + \hat{s}_j^2), \\
        &= 1 + \|\hat{\bs}\|_1 + \frac{1}{2} \|\hat{\bs}\|^2 + \frac{1}{2}\sum_{i=1}^{n} (n \hat{s}_i^2 + \|\hat{\bs}\|^2), \\
        &= 1 + \|\hat{\bs}\|_1 + \frac{1}{2} \|\hat{\bs}\|^2 + n\|\hat{\bs}\|^2,
    \end{align}
    and so combining with $\|\hat{\bs}\| \leq 1$ and $\|\hat{\bs}\|_1 \leq \sqrt{n}$, we get $\|\bphi(\hat{\bs})\|_1 \leq n + \sqrt{n} + \frac{3}{2}$, which gives
    \begin{align}
        |m(\by) - T_{f,\bx,2}(\by)| \leq \frac{\Lhess}{6} \left(n + \sqrt{n} + \frac{3}{2}\right) \beta^3 \|\hat{\bQ}^{-1}\|_{\infty}  \Delta^3.
    \end{align}
    The result then follows from \lemref{lem_fully_linear_quadratic_from_taylor_approx}(b).
\end{proof}

One way to choose the $p$ sample points to ensure that $\|\hat{\bQ}^{-1}\|$ is not too large is to use $\bx$ along with the  positive and negative perturbations along the coordinate axes $\bx \pm \Delta \be_i$ for all $i\in\{1,\ldots,n\}$, and $\bx + \Delta(\be_i + \be_j)$ for all combinations $i\neq j \in \{1,\ldots,n\}$.
That is,
\begin{align}
    \{\by_1, \ldots, \by_p\} &= \{\bx, \bx + \Delta\be_1, \ldots, \bx + \Delta\be_n, \bx - \Delta\be_1, \ldots, \bx-\Delta\be_n, \nonumber \\ 
    & \bx + \Delta(\be_1+\be_2), \ldots, \bx + \Delta (\be_1 + \be_n), \bx + \Delta(\be_2 + \be_3), \ldots, \bx + \Delta(\be_{n-1} + \be_n)\}. \label{eq_fully_quadratic_example_points}
\end{align}

\begin{lemma} \label{lem_quadratic_sample_estimate}
    For the interpolation set given by \eqref{eq_fully_quadratic_example_points}, we have $\|\hat{\bQ}^{-1}\|_{\infty} \leq 8$.
\end{lemma}
\begin{proof}
    Suppose we solve $\hat{\bQ}\bu = \bv$ for some right-hand side $\bv$ with $v_i = f(\by_i)$ and $\|\bv\|_{\infty}=1$, and where $\bu := \begin{bmatrix} c \\ \hat{\bg} \\ \upper(\hat{\bH}) \end{bmatrix}$.
    Then $\max_{\bv} \|\bu\|_{\infty} = \|\hat{\bQ}^{-1}\|_{\infty}$ and so it suffices to find an upper bound on $\|\bu\|_{\infty}$.
    
    The choice of interpolation points \eqref{eq_fully_quadratic_example_points} gives
    \begin{align}
        \{\hat{\bs}_1, \ldots, \hat{\bs}_p\} &= \{\bm{0}, \be_1, \ldots, \be_n, -\be_1, \ldots, -\be_n,  \be_1+\be_2, \ldots, \be_1 + \be_n, \be_2 + \be_3, \ldots, \be_{n-1} + \be_n\}.
    \end{align}
    in the linear system \eqref{eq_quadratic_interp_system_scaled}.
    The first row of \eqref{eq_quadratic_interp_system_scaled} immediately gives $c = v_1$ and so $|c| \leq 1$ since $\|\bv\|_{\infty} \leq 1$.

    Next, for $i\in\{1,\ldots,n\}$, the values $v_{i+1}$ and $v_{n+i+1}$ (corresponding to $\hat{\bs} = \pm \be_i$) completely determine $\hat{g}_i$ and $\hat{H}_{i,i}$ via rows $i+1$ and $n+i+1$ of \eqref{eq_quadratic_interp_system_scaled}, namely
    \begin{align}
        c + \hat{g}_i + \frac{1}{2}\hat{H}_{i,i} = v_{i+1}, \qquad \text{and} \qquad c - \hat{g}_i + \frac{1}{2}\hat{H}_{i,i} = v_{n+i+1},
    \end{align}
    which gives $\hat{g}_i = \frac{1}{2}(v_{i+1} - v_{n+i+1})$ and $\hat{H}_{i,i} = v_{i+1} + v_{n+i+1} - 2c$, and so $|\hat{g}_i| \leq 1$ and $|\hat{H}_{i,i}| \leq 4$ from $|c|,|v_i|,|v_{n+i+1}| \leq 1$.
    
    Lastly, in $\hat{\bQ}$, the column corresponding to $\hat{H}_{i,j}$ for $i<j$ only has one nonzero entry, a 1 in the $t$-th row corresponding to $\hat{\bs}_t = \be_i + \be_j$.
    Hence we get $\hat{H}_{i,j} = v_t - c - \hat{g}_i - \hat{g}_j - \frac{1}{2}\hat{H}_{i,i} - \frac{1}{2}\hat{H}_{j,j}$ and so $|\hat{H}_{i,j}| \leq 8$.
\end{proof}

All together, if we run our \renaming{MBDFO} algorithm \algref{alg_basic_tr_dfo_2}  with quadratic models given by solving \eqref{eq_quadratic_interp_system} with points \eqref{eq_fully_quadratic_example_points} at each iteration, we get the below second-order complexity result.

\begin{corollary} \label{cor_fully_quadratic_wcc}
    Under the assumptions of \thmref{thm_wcc_dfo_2} and \minorrev{the sequence $\{\|\grad^2 f(\bx_k)\|\}_{k=0}^{\infty}$ is} uniformly bounded\footnote{This is required for the derivative-based theory, via \assref{ass_model}\ref{ass_model_H} and $\bH_k=\grad^2 f(\bx_k)$ in \thmref{thm_wcc_2}.}, if the local quadratic model at each iteration of \algref{alg_basic_tr_dfo_2} is generated by quadratic interpolation to the points \eqref{eq_fully_quadratic_example_points} (with $\bx=\bx_k$ and $\Delta=\Delta_k$), then the iterates of \algref{alg_basic_tr_dfo_2} achieve second-order optimality $\sigma_k < \epsilon$ for the first time after at most $\bigO(\revision{n^4} \epsilon^{-3})$ iterations and $\bigO(\revision{n^6} \epsilon^{-3})$ objective evaluations.
\end{corollary}
\begin{proof}
    We note that all the interpolation points \eqref{eq_fully_quadratic_example_points} satisfy $\|\by_i-\bx\| \leq \sqrt{2} \Delta$ (i.e.~$\beta=\sqrt{2}$ in the assumptions of \thmref{thm_fully_quadratic}).
    From \thmref{thm_fully_quadratic} and \lemref{lem_quadratic_sample_estimate}, the model at each iteration is fully quadratic with $\kappamf, \kappamg, \kappamh = \bigO(n \Lhess)$.
    We also get $\|\bH_k\| \leq \|\grad^2 f(\bx_k)\| + \kappamh \Delta_{\max}$ from \eqref{eq_fully_quadratic_h} and so $\kappa_H = \bigO(\kappamh)$.
    We then apply \corref{cor_wcc_dfo_2} to get the iteration bound.
    The objective evaluation bound comes from observing that our interpolation mechanism requires evaluating $f$ at no more than $p = (n+1)(n+2)/2$ points per iteration.
\end{proof}

\revision{Just as in the linear interpolation case (\corref{cor_fully_linear_wcc_improved}), we can improve the explicit dependence on dimension in the worst-case complexity bound by making the coordinate perturbations depend on the dimension.
However, this construction does not align with a single interpolation set, rather it uses different objective evaluations to construct $\bg_k$ as $\bH_k$ in iteration $k$.
Specifically, we use the central difference approximation to the gradient with stepsize $\Delta_k/n^{1/4}$, and forward difference approximation to the Hessian with stepsize $\Delta_k/n$:
\begin{subequations} \label{eq_fully_quadratic_example_points_improved}
\begin{align}
    [\bg_k]_i &= \frac{f(\bx_k+\Delta_k \be_i / n^{1/4}) - f(\bx_k-\Delta_k \be_i / n^{1/4})}{2(\Delta_k / n^{1/4})}, \\
    [\bH_k]_{i,j} &= \frac{f(\bx_k+\Delta_k(\be_i+\be_j)/n) - f(\bx_k+\Delta_k \be_i/n) - f(\bx_k+\Delta_k \be_j/n) + f(\bx_k)}{(\Delta_k/n)^2},
\end{align}
\end{subequations}
for all $i,j=1,\ldots,n$, which requires $(n+1)(n+2)/2+2n$ objective evaluations because the calculation of $\bg_k$ cannot re-use any objective evaluations used to construct $\bH_k$.

\begin{corollary} \label{cor_fully_quadratic_wcc_improved}
    Under the assumptions of \thmref{thm_wcc_dfo_2} and the sequence $\{\|\grad^2 f(\bx_k)\|\}_{k=0}^{\infty}$ is uniformly bounded, if the local quadratic model at each iteration of \algref{alg_basic_tr_dfo_2} is generated as in \eqref{eq_fully_quadratic_example_points_improved}, then the iterates of \algref{alg_basic_tr_dfo_2} achieve second-order optimality $\sigma_k < \epsilon$ for the first time after at most $\bigO( \epsilon^{-3})$ iterations and $\bigO(n^2 \epsilon^{-3})$ objective evaluations.
\end{corollary}
\begin{proof}
    From \cite[Theorem 2.2]{Berahas2022} or \cite[Lemma 5]{Doikov2025} we have $\|\bg_k-\grad f(\bx_k)\| \leq \frac{\Lhess}{6} \Delta_k^2$ and from \cite[Proposition 2.7]{Cao2024} or \cite[Lemma 6]{Doikov2025} we have $\|\bH_k-\grad^2 f(\bx_k)\| \leq \frac{(\sqrt{2}+1) \Lhess}{3} \Delta_k$.
    By arguments similar to the proof of \lemref{lem_fully_linear_quadratic_from_taylor_approx}(b), the model is fully quadratic with $\kappamf,\kappamg,\kappamh=\bigO(\Lhess)$ (i.e.~no explicit dependence on $n$).
    The remainder of the proof follows that of \corref{cor_fully_quadratic_wcc}, but with the slightly larger $(n+1)(n+2)/2+2n$ objective evaluations per iteration.
\end{proof}

In terms of the dependency on $n$, this objective evaluation bound again matches what we may expect from a derivative-based method using (genuine small-$h$) finite differencing to approximate gradients and Hessians.
}

\subsection{Underdetermined Quadratic Interpolation} \label{sec_min_frob}

So far, we have the option of constructing either fully linear models with linear interpolation (using $n+1$ points), or constructing fully quadratic models with quadratic interpolation (using $(n+1)(n+2)/2$ points).
In most cases, quadratic models give much better practical performance because they can adapt to the local curvature of $f$, and so should be preferred over linear models, but here they come with the downside of requiring significantly more interpolation points---and hence objective evaluations---at each iteration.
In this section, we outline an approach for constructing quadratic models \eqref{eq_quadratic_model_generic} with $p\in\{n+2, \ldots, (n+1)(n+2)/2-1\}$ interpolation points.
In this case we have more degrees of freedom in the model than available evaluations, and so the interpolation problem is underdetermined (i.e.~there are infinitely many quadratic functions satisfying the interpolation conditions).
We describe one particular approach, minimum Frobenius norm underdetermined quadratic interpolation.

\minorrev{To that end}, suppose we again have $p$ interpolation points $\by_1,\ldots,\by_p$, but now with $p\in\{n+2, \ldots, (n+1)(n+2)/2-1\}$.
In general, there are infinitely many quadratic models $m$ \eqref{eq_quadratic_model_generic} satisfying $m(\by_i) = f(\by_i)$ for all $i=1,\ldots,p$.
Of these, we pick the model for which the Hessian $\bH$ has minimal Frobenius norm:
\begin{align}
    \min_{c,\bg,\bH} \frac{1}{4}\|\bH\|_F^2, \qquad \text{s.t.} \quad \text{$m(\by_i) = f(\by_i)$ for all $i=1,\ldots,p$}, \quad \text{and} \quad \bH = \bH^T. \label{eq_min_frob_problem}
\end{align}
Since $m$ is linear in the coefficients of $c$, $\bg$ and $\bH$, this is a convex quadratic program with equality constraints, and so can be solved via a linear system.

\begin{lemma}[Section 2, \cite{Powell2004}]
    The solution to \eqref{eq_min_frob_problem} can be obtained by solving the $(p+n+1)\times (p+n+1)$ linear system\footnote{\revision{We use $\bF$ to denote the matrix in \eqref{eq_min_frob_interp_system} as we are doing Frobenius-norm interpolation.}}
    \begin{align}
       \underbrace{\left[\begin{array}{c|c} \bP & \bM \\ \hline \bM^T & \bm{0} \end{array}\right]}_{=: \bF} \left[\begin{array}{c} \lambda_1 \\ \vdots \\ \lambda_p \\ \hline c \\ \bg \end{array}\right] = \left[\begin{array}{c} f(\by_1) \\ \vdots \\ f(\by_p) \\ \hline 0 \\ \b{0} \end{array}\right], \label{eq_min_frob_interp_system}
    \end{align}
    where $\bP\in\R^{p\times p}$ has entries $P_{i,j} = \frac{1}{2}[(\by_i-\bx)^T (\by_j-\bx)]^2$ and 
    \revision{\begin{align}
        \bM = \begin{bmatrix} 1 & (\by_1-\bx)^T \\ \vdots & \vdots \\ 1 & (\by_p-\bx)^T \end{bmatrix} \in \R^{p\times (n+1)},
    \end{align}
    c.f.~\eqref{eq_linear_interp_system}.}
    The model Hessian is given by $\bH = \sum_{i=1}^{p} \lambda_i (\by_i-\bx) (\by_i-\bx)^T$.
\end{lemma}


\begin{remark} \label{rem_fully_quad_is_min_frob}
    If $p=(n+1)(n+2)/2$, then we can still get a model by solving \eqref{eq_min_frob_problem} via \eqref{eq_min_frob_interp_system}.
    However, since we have sufficient degrees of freedom that the interpolation conditions uniquely define the model, there is only one valid $\bH$ (corresponding to fully quadratic interpolation), and so we are better off solving the smaller system \eqref{eq_quadratic_interp_system}.
    However, this does mean that all the below results also apply to fully quadratic models.
\end{remark}

Once again, if $\Delta$ is a typical size of $\|\by_i-\bx\|$ and we define $\hat{\bs}_i := (\by_i-\bx)/\Delta$ for all $i\in\{1,\ldots,p\}$, we get the scaled interpolation system
\begin{align}
   \underbrace{\left[\begin{array}{c|c} \hat{\bP} & \hat{\bM} \\ \hline \hat{\bM}^T & \bm{0} \end{array}\right]}_{=: \hat{\bF}} \left[\begin{array}{c} \hat{\lambda}_1 \\ \vdots \\ \hat{\lambda}_p \\ \hline c \\ \hat{\bg} \end{array}\right] = \left[\begin{array}{c} f(\by_1) \\ \vdots \\ f(\by_p) \\ \hline 0 \\ \b{0} \end{array}\right], \label{eq_min_frob_interp_system_scaled}
\end{align}
where $\hat{P}_{i,j} = \frac{1}{2}[\hat{\bs}_i^T \hat{\bs}_j]^2$ and $\hat{\bM}$ \revision{has $i$-th row $[1 \: \hat{\bs}_i^T]$ (c.f.~\eqref{eq_linear_interp_system_scaled})} 
We get the coefficients $\hat{\bg} = \Delta \: \bg$ and $\hat{\lambda}_i = \Delta^4 \lambda_i$, and so $\hat{\bH} = \sum_{i=1}^{p} \hat{\lambda}_i \hat{\bs}_i \hat{\bs}_i^T$ gives $\hat{\bH} = \Delta^2 \bH$.

\begin{remark} \label{rem_P_semidefinite}
We note that $\bF$ and $\hat{\bF}$ are symmetric, and both $\bP$ and $\hat{\bP}$ are positive semidefinite.
This can be seen, e.g.~for $\hat{\bP}$, by observing that for any $\bv\in\R^p$ we have
\begin{align}
    \bv^T \hat{\bP} \bv = \frac{1}{2} \sum_{i,j=1}^{p} v_i v_j \left(\sum_{k=1}^{n} \hat{s}_{i,k} \hat{s}_{j,k} \right)^2 = \frac{1}{2} \sum_{i,j=1}^{p} \sum_{k,\ell=1}^{n} v_i v_j \hat{s}_{i,k} \hat{s}_{j,k} \hat{s}_{i,\ell} \hat{s}_{j,\ell} = \frac{1}{2} \sum_{k,\ell=1}^{n} \left(\sum_{i=1}^{p} v_i \hat{s}_{i,k} \hat{s}_{i,\ell}\right)^2 \geq 0.
\end{align}
Hence $\bF$ and $\hat{\bF}$ yield \emph{saddle point} linear systems, which have been widely studied \cite{Benzi2005}. 
\end{remark}

Since we choose to ensure $\|\bH\|_F$ is small, we can provide an upper bound on the size of the model Hessian.
This will be necessary for us to apply \lemref{lem_fully_linear_quadratic_from_taylor_approx}(a) to show that the model is fully linear.

\begin{lemma} \label{lem_min_frob_bounded_hessian}
    Suppose $f$ satisfies \minorrev{\assref{ass_smoothness_1}\ref{ass_smoothness_1_smooth}} and we construct a quadratic model $m$ \eqref{eq_quadratic_model_generic} for $f$ by solving \eqref{eq_min_frob_interp_system_scaled}, where we assume $\hat{\bF}$ is invertible. 
    If $\|\by_i-\bx\| \leq \beta \Delta$ for some $\beta>0$ and all $i=1,\ldots,p$, then the model Hessian satisfies
    \begin{align}
        \|\bH\| \leq \kappa_H := \frac{\Lgrad}{2} p \beta^4 \|\hat{\bF}^{-1}\|_{\infty}.
    \end{align}
\end{lemma}
\begin{proof}
    By suitably modifying $c$ and $\bg$, the model Hessian does not change if the objective function is changed by adding a linear function to it.\footnote{From changing $c$ and $\bg$ suitably, any Hessian that interpolates $f$ also interpolates $f$ plus a linear function, so the minimizer $\bH$ must be the same.}
    Hence, we consider interpolating to the function $\tilde{f}(\by) := f(\by) - f(\bx) - \grad f(\bx)^T (\by-\bx)$, and so by \lemref{lem_lsmooth} we have $|\tilde{f}(\by_i)| \leq \frac{\Lgrad}{2}\|\by_i-\bx\|^2$, and so all terms in the right-hand side of \eqref{eq_min_frob_interp_system_scaled} are bounded by $\frac{\Lgrad}{2}\beta^2\Delta^2$.
    This means that for all $i=1,\ldots,p$ we have
    \begin{align}
        |\lambda_i| = \Delta^{-4} |\hat{\lambda}_i| \leq \frac{\Lgrad}{2} \beta^2 \Delta^{-2} \|\hat{\bF}^{-1}\|_{\infty}.
    \end{align}
    Finally, we have
    \begin{align}
        \|\bH\| \leq \sum_{i=1}^{p} |\lambda_i| \: \|\by_i-\bx\|^2 \leq p\left(\frac{\Lgrad}{2} \beta^2 \Delta^{-2} \|\hat{\bF}^{-1}\|_{\infty}\right) \beta^2 \Delta^2,
    \end{align}
    and the result follows.
\end{proof}

We can now establish that our minimum Frobenius norm quadratic interpolation model is indeed fully linear.

\begin{theorem} \label{thm_min_frob_fully_linear}
    Suppose $f$ satisfies \minorrev{\assref{ass_smoothness_1}\ref{ass_smoothness_1_smooth}} and we construct a quadratic model $m$ \eqref{eq_quadratic_model_generic} for $f$ by solving \eqref{eq_min_frob_interp_system_scaled}, where we assume $\hat{\bF}$ is invertible.
    If $\|\by_i-\bx\| \leq \beta \Delta$ for some $\beta>0$ and all $i=1,\ldots,p$, then the model is fully linear with constants
    \begin{align}
        \kappamf = \frac{\Lgrad + \kappa_H}{2}(1+\sqrt{n}) \beta^2 \|\hat{\bM}^{\dagger}\|_{\infty} + \frac{\Lgrad + \kappa_H}{2}, \qquad \text{and} \qquad \kappamg = 2\kappamf + 2\kappa_H,
    \end{align}
    using the value of $\kappa_H$ defined in \lemref{lem_min_frob_bounded_hessian}.
\end{theorem}
\begin{proof}
    Fix $\by\in B(\bx,\Delta)$ and define $\hat{\bs} := (\by-\bx)/\Delta$.
    Since $\hat{\bF}$ is invertible, $\hat{\bM}$ has full column rank \cite[Theorem 3.3]{Benzi2005}, \revision{so the rows of $\hat{\bM}$ span $\R^{n+1}$.
    Hence, the system $\hat{\bM}^T \bv=\begin{bmatrix} 1 \\ \hat{\bs} \end{bmatrix}$ is consistent. 
    So, we take $\bv=(\hat{\bM}^T)^{\dagger} \begin{bmatrix} 1 \\ \hat{\bs} \end{bmatrix}$ to be the minimal norm solution, with $\|\bv\|_1 \leq \|\hat{\bM}^{\dagger}\|_{\infty} (1+\sqrt{n})$ by similar reasoning as in the proof of \thmref{thm_fully_linear}.
    Multiplying the last $n$ rows of this (consistent) equation by $\Delta$, we again get \eqref{eq_linear_interp_error_tmp3} for our new $\bv$.}
    By following the same argument as for \revision{\eqref{eq_linear_interp_error_tmp2a}} and using the interpolation condition $m(\by_i)=f(\by_i)$ we get
    \begin{align}
        |c + \bg^T (\by-\bx) - f(\bx) - \grad f(\bx)^T (\by-\bx)| &= \left|\sum_{i=1}^{p} v_i \left(c + \bg^T (\by_i-\bx) - f(\bx) - \grad f(\bx)^T (\by_i-\bx)\right)\right|, \\
        &\leq \left|\sum_{i=1}^{p} v_i (m(\by_i) - f(\bx) - \grad f(\bx)^T (\by_i-\bx))\right| \nonumber \\
        &\qquad\qquad\qquad + \left|\sum_{i=1}^{p} \frac{1}{2}v_i (\by_i-\bx)^T \bH (\by_i-\bx)\right|, \\
        &\leq \frac{\Lgrad}{2} \beta^2 \Delta^2 \|\bv\|_1 + \frac{1}{2} \|\bH\| \beta^2 \Delta^2 \|\bv\|_1, \label{eq_min_frob_fully_linear_tmp1}
    \end{align}
    and so
    \begin{align}
        |c + \bg^T (\by-\bx) - f(\bx) - \grad f(\bx)^T (\by-\bx)| \leq \frac{\Lgrad + \kappa_H}{2}(1+\sqrt{n}) \beta^2 \|\hat{\bM}^{\dagger}\|_{\infty}  \Delta^2, \label{eq_min_frob_fully_linear_tmp2}
    \end{align}
    which finally gives
    \begin{align}
        |m(\by) - f(\bx) - \grad f(\bx)^T (\by-\bx)| &\leq \frac{\Lgrad + \kappa_H}{2}(1+\sqrt{n}) \beta^2 \|\hat{\bM}^{\dagger}\|_{\infty}  \Delta^2 + \frac{1}{2}|(\by-\bx)^T \bH (\by-\bx)|, \\
        &\leq \left(\frac{\Lgrad + \kappa_H}{2}(1+\sqrt{n}) \beta^2 \|\hat{\bM}^{\dagger}\|_{\infty}  + \frac{\kappa_H}{2}\right) \Delta^2. \label{eq_min_frob_fully_linear_tmp3}
    \end{align}
    The result then follows from \lemref{lem_fully_linear_quadratic_from_taylor_approx}(a).
\end{proof}

\revision{A popular structured interpolation set used for minimum Frobenius norm interpolation is $\by_1=\bx$, $\by_{i+1}=\bx+\Delta\be_i$ and $\by_{n+i+1}=\bx-\Delta\be_i$ for $i=1,\ldots,n$, with $p=2n+1$.
This gives the associated interpolation matrix}
\begin{align}
    \hat{\bF} = \left[\begin{array}{ccc|cc} 0 & \bm{0}^T & \bm{0}^T & 1 & \bm{0}^T \\ \bm{0} & \frac{1}{2} \bI & \frac{1}{2} \bI & \bm{e} & \bI \\ \bm{0} & \frac{1}{2} \bI & \frac{1}{2} \bI & \bm{e} & -\bI \\ \hline 1 & \bm{e}^T & \bm{e}^T & 0 & \bm{0}^T \\ \bm{0} & \bI & -\bI & \bm{0} & \bm{0} \end{array}\right], \qquad \text{with} \qquad \hat{\bF}^{-1} = \left[\begin{array}{ccc|cc} 2n & -\bm{e}^T & -\bm{e}^T & 1 & \bm{0}^T \\ -\bm{e} & \frac{1}{2} \bI & \frac{1}{2} \bI & \bm{0} & \frac{1}{2} \bI \\ -\bm{e} & \frac{1}{2} \bI & \frac{1}{2} \bI & \bm{0} & -\frac{1}{2} \bI \\ \hline 1 & \bm{0}^T & \bm{0}^T & 0 & \bm{0}^T \\ \bm{0} & \frac{1}{2} \bI & -\frac{1}{2} \bI & \bm{0} & \bm{0} \end{array}\right], \label{eq_min_frob_example_invF}
\end{align}
and so $\|\hat{\bF}^{-1}\|_{\infty} = 4n+1$.
Applying \lemref{lem_min_frob_bounded_hessian}, we get the bound $\kappa_H = \bigO(pn) = \bigO(n^2)$.
Combining with $\|\hat{\bM}^{\dagger}\|_{\infty}=1$ \revision{(which can be verified by direct computation)}, \thmref{thm_min_frob_fully_linear} tells us that the model is fully linear with $\kappamf,\kappamg = \bigO(\sqrt{n}\: \kappa_H) = \bigO(n^{5/2})$.
Applying these results to \corref{cor_wcc_dfo}, we get a worst-case complexity for \algref{alg_basic_tr_dfo} of $\bigO(n^7 \epsilon^{-2})$ iterations and $\bigO(n^8 \epsilon^{-2})$ objective evaluations.
This is significantly worse than linear and even (fully) quadratic interpolation!

However, luckily for us it turns out that the bound $\kappa_H = \bigO(n^2)$ from \lemref{lem_min_frob_bounded_hessian} is overly conservative in this case, because only the first row of $\hat{\bF}^{-1}$ has an absolute row sum of size $\bigO(n)$, and this row is only used to calculate the model coefficient $\lambda_1$.
Moreover, $\lambda_1$ only contributes to the model Hessian via $\lambda_1 (\by_1-\bx)(\by_1-\bx)^T$, so since $\by_1=\bx$, it does not affect the Hessian.
Indeed, using the explicit form of $\hat{\bF}^{-1}$ \eqref{eq_min_frob_example_invF}, we compute that $\lambda_{i+1}=\lambda_{n+i+1} = \Delta^{-4}[-f(\bx) + \frac{1}{2}f(\bx+\Delta \be_i) + \frac{1}{2}f(\bx-\Delta \be_i)]$ for $i=1,\ldots,n$, and so the model Hessian is
\begin{align}
    \bH = \sum_{i=1}^{n} (\lambda_{i+1}+\lambda_{n+i+1}) \Delta^2 \be_i \be_i^T = \operatorname{diag}\left(\left\{\frac{f(\bx+\Delta \be_i) - 2 f(\bx) + f(\bx-\Delta \be_i)}{\Delta^2} : i=1,\ldots,n \right\}\right). \label{eq_min_frob_example_hess_bound}
\end{align}
As per the proof of \lemref{lem_min_frob_bounded_hessian}, we may without loss of generality assume $|f(\by_i)| \leq \frac{\Lgrad}{2}\beta^2 \Delta^2$, and so we get the improved bound $\|\bH\| = \max_{i=1,\ldots,n} |H_{i,i}| \leq 2 \Lgrad \beta^2$.
Hence we may take $\kappa_H = \bigO(1)$, a significant improvement over \lemref{lem_min_frob_bounded_hessian}.

Taking $\kappamf,\kappamg=\bigO(\sqrt{n}\: \kappa_H)$ from \thmref{thm_min_frob_fully_linear} with our tighter bound $\kappa_H = \bigO(1)$ in \corref{cor_wcc_dfo}, we arrive at the same worst-case complexity as linear interpolation.

\begin{corollary} \label{cor_wcc_dfo_min_frob}
    Under the assumptions of \thmref{thm_wcc_dfo}, if the local quadratic model at each iteration of in \algref{alg_basic_tr_dfo} is generated by minimum Frobenius norm quadratic interpolation to the points $\{\bx_k, \bx_k + \Delta_k \bm{e}_1, \ldots, \bx_k + \Delta_k \bm{e}_n, \bx_k - \Delta_k \be_1, \ldots, \bx_k - \Delta_k \be_n\}$, then the iterates of \algref{alg_basic_tr_dfo} achieve $\|\grad f(\bx_k)\| < \epsilon$ for the first time after at most $\bigO(n \epsilon^{-2})$ iterations and $\bigO(n^2 \epsilon^{-2})$ objective evaluations.
\end{corollary}

\revision{Again, by scaling the perturbations in terms of the problem dimension, we can get improved dependence on dimension in the worst-case complexity bounds.

\begin{corollary} \label{cor_wcc_dfo_min_frob_improved}
    Under the assumptions of \thmref{thm_wcc_dfo}, and if $\Delta_k \leq \Delta_{\max}$ for all $k$, if the local quadratic model at each iteration of in \algref{alg_basic_tr_dfo} is generated by minimum Frobenius norm quadratic interpolation to the points $\{\bx_k, \bx_k + \Delta_k \bm{e}_1 / n^{1/4}, \ldots, \bx_k + \Delta_k \bm{e}_n / n^{1/4}, \bx_k - \Delta_k \be_1 / n^{1/4}, \ldots, \bx_k - \Delta_k \be_n / n^{1/4}\}$, then the iterates of \algref{alg_basic_tr_dfo} achieve $\|\grad f(\bx_k)\| < \epsilon$ for the first time after at most $\bigO( \epsilon^{-2})$ iterations and $\bigO(n \epsilon^{-2})$ objective evaluations.
\end{corollary}
\begin{proof}
    This construction gives the same $\bg_k$ as in \eqref{eq_fully_quadratic_example_points_improved}, and so $\|\bg_k-\grad f(\bx_k)\| \leq \frac{\Lgrad}{2}\Delta_k^2 \leq \frac{\Lgrad \Delta_{\max}}{2} \Delta_k$.
    By the same reasoning as above, we have $\kappa_H=\bigO(1)$, and so $\kappamf,\kappamg=\bigO(1)$.
    The remainder of the proof is identical to the linear interpolation case (\corref{cor_fully_linear_wcc_improved}).
\end{proof}
}

\begin{remark}
    A simple modification of this approach, which can be beneficial in practice, is to look at the minimum change in the Hessian between successive iterations of our main algorithm.
    That is, given some old Hessian approximation $\Hprev$, we choose $\bH = \Hprev + \Delta \bH$ where $\Delta \bH$ solves
    \begin{align}
        \min_{c,\bg,\Delta \bH} \frac{1}{4}\|\Delta \bH\|_F^2, \qquad \text{s.t.} \quad \text{$m(\by_i-\bx) = f(\by_i)$ for all $i=1,\ldots,p$} \quad \text{and} \quad \Delta \bH = \Delta \bH^T. \label{eq_min_chg_frob_problem}
    \end{align}
    The corresponding linear system is the same as \eqref{eq_min_frob_interp_system} but where the entries in the right-hand side are replaced with $f(\by_i) - \frac{1}{2}(\by_i-\bx)^T \Hprev (\by_i-\bx)$, and where $\bH = \Hprev + \sum_{i=1}^{p} \lambda_i (\by_i-\bx) (\by_i-\bx)^T$.
    Following the proof of \lemref{lem_min_frob_bounded_hessian}, we get $\|\Delta \bH\| \leq \frac{\Lgrad + \|\Hprev\|}{2} p \beta^4 \|\hat{\bF}^{-1}\|_{\infty}$ and so our model Hessian can potentially grow rapidly,
    \begin{align}
        \|\bH\| \leq \|\Hprev\| + \|\Delta \bH\| \leq \left(1 + \frac{1}{2} p \beta^4 \|\hat{\bF}^{-1}\|_{\infty}\right)\|\Hprev\| + \frac{\Lgrad}{2} p \beta^4 \|\hat{\bF}^{-1}\|_{\infty}.
    \end{align}
    Unfortunately this rate of growth in the model Hessian is sufficient that global convergence results do not apply.\footnote{Trust-region methods can converge in the case of unbounded model Hessians, but the growth cannot exceed $\|\bH_k\| = \bigO(k)$; see \cite{Toint1988,Diouane2024} for more details.}
\end{remark}

\revision{\subsection{Practical Model Construction} \label{sec_interp_practical_ibcdfo}
We have seen that choosing interpolation sets based on structured perturbations around $\bx_k$ can yield fully linear/\minorrev{fully} quadratic models.
Although these can yield very good worst-case complexity bounds (e.g.~\corref{cor_fully_linear_wcc_improved}), they do not allow us to reuse any existing objective evaluations during model construction.
Given we often turn to \renaming{MBDFO} methods when objective evaluations are expensive, a pragmatic approach would allow us to construct models based on existing objective evaluations, while controlling the size of the relevant matrix norms.

To do this, given $\bx_k$, $\Delta_k$ and a collection of existing evaluations $\mathcal{Y}\subset\R^n$ that are sufficiently close to $\bx_k$, we first try to find a subset of $\mathcal{Y}$ that (together with $\bx_k$) is suitable for linear interpolation.
If we cannot find $n$ points in $\mathcal{Y}$ that yield a sufficiently good linear interpolation set, we augment this set based on the null space of selected perturbations around $\bx_k$ (i.e.~$\by-\bx_k$ for selected points $\by\in\mathcal{Y}$).
By maintaining a QR factorization of the selected perturbations, we can control $\|\hat{\bM}^{-1}\|$ \eqref{eq_linear_interp_system_scaled} directly based on the diagonal entries of the $R$ matrix \cite[Lemma 2.4]{Wild2013}.
This procedure is formalized in \cite[Figure 4.2]{Wild2008}, and extended to bound-constrained problems in \cite{Curtis2024}.

We can then optionally choose to extend the selected interpolation set by adding more points (that are sufficiently close to $\bx_k$) while ensuring that $\|\hat{\bF}^{-1}\|$ \eqref{eq_min_frob_interp_system_scaled} is not too large.
Two approaches for doing this are:
\begin{itemize}
    \item Ensuring the new point added does not significantly increase $\|\hat{\bP}^{-1}\|$ by using a QR factorization of $\hat{\bM}^T$ and a Cholesky-like factorization of $\hat{\bP}$, as described in \cite[Algorithm 4.2]{Wild2008a};
    \item Ensuring the new point added does not significantly increase $\|\hat{\bF}^{-1}\|$ by controlling the Schur complement of the resulting $\hat{\bF}$ after the new point is added \cite[Algorithm 5.2]{Roberts2024}.
\end{itemize}
The first of these approaches is used in the IBCDFO software collection (see \secref{sec_conclusion}).

In both these approaches, points are \emph{appended} to the interpolation set from the set of existing evaluations, with the model quality controlled via the norm of the (inverse) interpolation matrix as per the error bounds derived in this section. 
In \secref{sec_self_correcting}, we will instead consider geometric approaches for incrementally improving an interpolation set, where we take an existing interpolation set and determine how to \emph{change} a small number of points to improve its quality.
}

\subsection{Composite Models} \label{sec_composite_models}
Another important construction is for models where the objective function $f$ has some known structure.
For example, this may be where $f$ results from composing a known function with a function for which derivatives are unavailable, namely
\begin{align}
    f(\bx) = \minorrev{H(\br(\bx))}, \label{eq_composite_generic}
\end{align}
where $\br:\R^n\to\R^q$ (for some $q$) is smooth but has only a zeroth order oracle, and the smooth function $\minorrev{H}:\R^q\to\R$ and its derivatives are completely known.
Perhaps the most common example of \eqref{eq_composite_generic} is nonlinear least-squares minimization, where $\minorrev{H(\br)} = \frac{1}{2}\|\br\|^2$.

In such settings, because we receive more information with every zeroth order oracle query (the full vector $\br(\bx)$ instead of the scalar $f(\bx)$), it can often be sufficient to build linear interpolation models for $\br(\bx)$, and extend these to a quadratic model for $f$ using \eqref{eq_composite_generic}.
For example, in the nonlinear least-squares case, if we sample $\br(\by_i)$ for $i=1,\ldots,n+1$ and perform linear interpolation on each component of $\br(\bx)$, we get a model
\begin{align}
    \br(\by) \approx \bem(\by) := \bc + \bJ(\by-\bx),
\end{align}
for some $\bc\in\R^q$ and $\bJ\in\R^{q\times n}$.
We then naturally get a quadratic model for $f$ via
\begin{align}
    f(\by) \approx m(\by) := \frac{1}{2}\|\bem(\by)\|^2 = \frac{1}{2}\|\bc\|^2 + (\bJ^T\bc)^T (\by-\bx) + \frac{1}{2}(\by-\bx)^T (\bJ^T \bJ) (\by-\bx). \label{eq_gn_model}
\end{align}
In effect, by using the structure of the problem \eqref{eq_composite_generic} we get an approximate quadratic model for the cost of a linear interpolation model.
If each component of $\bem(\by)$ is a fully linear model for the corresponding component of $\br(\by)$, then it can be shown that the quadratic model \eqref{eq_gn_model} is a fully linear model for $f$ \cite[Lemma 3.3]{Cartis2019}.
Of course, if a quadratic model for $\br$ is available, then a more complex quadratic (or potentially up to quartic) model for $f$ can be derived in a similar way \cite{Zhang2010,Wild2017}.

This composite approach, where we receive more oracle information from $\br(\bx)$ and directly exploit some known structure, can significantly improve the practical performance of \renaming{MBDFO} methods.
These ideas also extend to the case where $H$ is nonsmooth (e.g.~$\minorrev{H(\br)}=\|\br\|_1$ or $\minorrev{H(\br)}=\max(\br)$), although solving the resulting trust-region step calculation and convergence theory is more complex \cite{Larson2021,Larson2024}.

\subsubsection*{Notes and References}
{\small The first fully linear/\minorrev{fully} quadratic error bounds for linear and quadratic interpolation were given in \cite{Conn2007}, and for linear regression and underdetermined quadratic interpolation in \cite{Conn2008}, although minimum Frobenius norm interpolation had been proposed earlier in \cite{Powell2004}. These results are collected in the reference book \cite{Conn2009} and summarized in  \cite{Schwertner2024}. \revision{The proofs presented here generalize to constrained interpolation (\secref{sec_convex_constraints}), and are based on the theory for linear interpolation and minimum Frobenius norm interpolation given in \cite{Hough2022,Roberts2024}.}


The interpolation constructions here can also be thought of as a way to approximate the gradient or Hessian of a function just with function evaluations. Such approximations, known as simplex gradients/Hessians, have been studied in \cite{Regis2015,Hare2020,Hare2024}.
\revision{Using perturbations of size $\bigO(\Delta_k/\sqrt{n})$ to improve the dimension dependence of the worst-case complexity bound for linear models (as in \corref{cor_fully_linear_wcc_improved}) was originally proposed for quadratic regularization methods in \cite{Grapiglia2023}, and was adapted to MBDFO in \cite{Chaudhry2025,Davar2025}. The extension to fully quadratic models is based on the constructions in \cite{Cao2024,Doikov2025}. In particular, \cite{Doikov2025} achieves a better bound of $\bigO(n^{3/2})$ objective evaluations by only re-sampling the Hessian every $n$ iterations. Although \cite{Doikov2025} uses a cubic regularization method, it is likely a similar improvement would be possible for a trust-region method.}

{\small How to best construct models in the underdetermined quadratic case is still actively being studied. Early alternatives to \eqref{eq_min_frob_problem} or \eqref{eq_min_chg_frob_problem} were proposed in \cite{Conn1996,Zhang2014}, but recent years have seen several alternatives proposed \cite{Xie2025,Xie2025a,Xie2025b}. 
The method in \cite{Xie2025b}, for example, attempts to optimize the parameters of a specific model construction problem over the course of the algorithm.}

Composite models for nonlinear least-squares problems were first introduced in \cite{Zhang2010} with the resulting worst-case complexity studied in \cite{Cartis2019}.
Similar ideas can also be used when $h$ is nonsmooth \cite{Grapiglia2016,Larson2021,Larson2024,Liu2024}, although the resulting trust-region subproblem becomes significantly more difficult to solve.
However, \renaming{MBDFO} for general nonsmooth problems is a relatively under-studied topic (see e.g.~\cite{Judice2015,Audet2020}), in contrast to the relatively established approaches for direct search DFO methods \cite{Audet2006a,Audet2022}.

We conclude by noting that other forms of interpolation model have been used in \renaming{MBDFO}.
The most notable example of non-polynomial models is the use of radial basis function (RBF) models, which take the form
\begin{align}
    f(\by) \approx m(\by) := c + \bg^T (\by-\bx) + \sum_{i=1}^{p} \lambda_i \psi(\|\by-\by_i\|),
\end{align}
where $\psi:[0,\infty)\to\R$ is some predetermined scalar function, such as the Gaussian function $\psi(r) = e^{-cr^2}$ for some $c>0$ \cite{Buhmann2000}. 
Because of the inclusion of a linear term in the model, such models can also be made fully linear through proper selection of the interpolation points \cite{Wild2008,Wild2013}.
RBF models are primarily designed to globally approximate a function (rather than locally as here), and so are more commonly used in global optimization algorithms~\cite{Bjorkmann2000,Gutmann2001}.
Such models have many similarities with statistically motivated global function approximations such as Gaussian Processes \cite{Rasmussen2006}, which are used in Bayesian (global) optimization algorithms~\cite{Shahriari2016} but again have also been used in \renaming{MBDFO} \cite{Augustin2017}.
}

\section{Incremental Geometry Improvement} \label{sec_self_correcting}

In \secref{sec_model_construction}, we described how to generate a set of interpolation points in order to be fully linear/\minorrev{fully} quadratic in $B(\bx,\Delta)$.
There, we had full control over where to place the interpolation points (e.g.~using coordinate perturbations around $\bx$) to ensure the model was sufficiently accurate.
\revision{This approach is of theoretical interest, and does motivate procedures to find a collection of suitable interpolation points from a database (and determine any extra evaluations required to get a good model).
However, in many practically successful \renaming{MBDFO} algorithms, an incremental approach is used, where interpolation points for the next iteration $k+1$ are taken to be the points from the current iteration $k$, with minimal changes to the set, such as adding the new iterate $\bx_{k+1}$.
This ensures we only ask for a very small number of new (possibly expensive) objective evaluations at each iteration.
Here, we describe how to assess the quality of an interpolation set using geometric conditions, which determine how to minimally alter a given interpolation set to ensure its quality.
}



The primary tool we will use to assess the quality of an interpolation set, and determine which points to move (and where to move them to), are \emph{Lagrange polynomials}.
For a set of points $\by_1,\ldots,\by_p\in\R^n$, the associated Lagrange polynomials are a collection of $p$ polynomials $\ell_1,\ldots,\ell_p:\R^n\to\R$ satisfying the conditions $\ell_i(\by_j) = \delta_{i,j}$, where $\delta_{i,j}$ is the Kronecker delta (i.e.~$\delta_{i,j}=1$ if $i=j$ and 0 otherwise).
The degree of the Lagrange polynomials will match the degree of the interpolation model we are trying to assess.

Lagrange polynomials are an important concept from polynomial approximation theory.
For example, suppose we have a continuous function $f:[-1,1]\to\R$ which we will approximate by a degree-$d$ polynomial $p_d$ via interpolation to points $x_1,\ldots,x_{d+1}\in[-1,1]$.
The associated Lagrange polynomials determine the \emph{Lebesgue constant} for the interpolation set,  $\Lambda_1 := \max_{x\in[-1,1]} \sum_{i=1}^{d+1} |\ell_i(x)|$.
We have the following result, which says that the interpolant $p_d$ has error within $\bigO(\Lambda_1)$ of the best possible error.

\begin{proposition}[Theorem 15.1, \cite{Trefethen2020}] \label{prop_lebesgue_constant}
    The error between the true (continuous) function $f:[-1,1]\to\R$ and degree-$d$ polynomial interpolant $p_d$ satisfies
    \begin{align}
        \|f - p_d\|_{\infty} \leq (\Lambda_1 + 1) \|f - p_d^*\|_{\infty},
    \end{align}
    where $\|f\|_{\infty}:=\max_{x\in[-1,1]} |f(x)|$ is the supremum norm and $p_d^* \minorrev{\in} \argmin_{p_d} \|f-p_d\|_{\infty}$ is \minorrev{an} optimal degree-$d$ polynomial interpolant to $f$.
\end{proposition}

That is, choosing interpolation sets which make the associated Lagrange polynomials small in magnitude are associated with better polynomial approximations.
We will see that a similar property holds in the multi-dimensional case.

We begin by re-deriving the fully linear/\minorrev{fully} quadratic constants from \secref{sec_model_construction}, replacing all matrix norms with quantities related to the magnitude of the relevant Lagrange polynomials. 
We will then use this to build a \renaming{MBDFO} algorithm with self-correcting interpolation sets.

\subsection{Interpolation Error Bounds}

To begin, consider the case of linear interpolation to build a model \eqref{eq_linear_model_generic} by solving the system \eqref{eq_linear_interp_system} or \eqref{eq_linear_interp_system_scaled}.
In this case, given the same base point $\bx\in\R^n$ as the model \eqref{eq_linear_model_generic}, the Lagrange polynomials are the linear functions
\begin{align}
    \ell_i(\by) := c_i + \bg_i^T (\by-\bx), \qquad \forall i=1,\ldots,p, \label{eq_lagrange_linear_generic}
\end{align}
where $p:=n+1$, defined by the interpolation conditions $\ell_i(\by_j) = \delta_{i,j}$ for $i,j=1,\ldots,p$.
That is, from \eqref{eq_linear_interp_system} or \eqref{eq_linear_interp_system_scaled} we can construct $\ell_i$ by solving
\begin{align}
    \bM \begin{bmatrix} c_i \\ \bg_i \end{bmatrix} = \be_i, \qquad \text{or} \qquad \hat{\bM} \begin{bmatrix} c_i \\ \Delta \: \bg_i \end{bmatrix} = \be_i, \label{eq_lagrange_linear_system}
\end{align}
provided that $\bM$ and $\hat{\bM}$ are invertible.
Comparing to \eqref{eq_linear_interp_system} or \eqref{eq_linear_interp_system_scaled}, one important consequence of this is that the interpolation model can be written as a linear combination of Lagrange polynomials,
\begin{align}
    \begin{bmatrix} c \\ \bg \end{bmatrix} = \sum_{i=1}^{p} f(\by_i) \begin{bmatrix} c_i \\ \bg_i \end{bmatrix}, \qquad \text{giving} \qquad m(\by) = \sum_{i=1}^{p} f(\by_i) \ell_i(\by). \label{eq_lagrange_interpolant}
\end{align}
If $f$ is constant, $f(\by)=1$ for all $\by\in\R^n$, then since the linear model $m$ is unique, the model must perfectly match the original function $m(\by)=f(\by)$.
Applying \eqref{eq_lagrange_interpolant}, we conclude that Lagrange polynomials form a partition of unity,
\begin{align}
    \sum_{i=1}^{p} \ell_i(\by) = 1, \qquad \forall \by\in\R^n. \label{eq_lagrange_partition_unity}
\end{align}
Separately, from \eqref{eq_lagrange_linear_generic} and \eqref{eq_lagrange_linear_system} we observe that
\begin{align}
    \ell_i(\by) = \begin{bmatrix} 1 \\ \by-\bx \end{bmatrix}^T \begin{bmatrix} c_i \\ \bg_i \end{bmatrix} = \be_i^T \bM^{-T} \begin{bmatrix} 1 \\ \by-\bx \end{bmatrix} = \be_i^T \hat{\bM}^{-T} \begin{bmatrix} 1 \\ (\by-\bx)/\Delta \end{bmatrix}.
\end{align}
This means that we can evaluate all Lagrange polynomials at a single point by solving a single linear system, without calculating any $c_i$ or $\bg_i$:
\begin{align}
    \blambda(\by) := \begin{bmatrix} \ell_1(\by) \\ \vdots \\ \ell_p(\by) \end{bmatrix} = \bM^{-T} \begin{bmatrix} 1 \\ \by-\bx \end{bmatrix} = \hat{\bM}^{-T} \begin{bmatrix} 1 \\ (\by-\bx)/\Delta \end{bmatrix}. \label{eq_lagrange_system_linear}
\end{align}
We have seen a similar relationship before, in \eqref{eq_linear_interp_error_tmp3} in the proof of \thmref{thm_fully_linear}.
Specifically, the vector $\bv$ in \eqref{eq_linear_interp_error_tmp3} is equal to $\blambda(\by)$, and so \eqref{eq_linear_interp_error_tmp2} can be written as
\begin{align}
    |m(\by) - f(\bx) - \grad f(\bx)^T (\by-\bx)| \leq \frac{\Lgrad}{2} \beta^2 \Delta^2 \|\blambda(\by)\|_1, \label{eq_linear_interp_lagrange_basic}
\end{align}
for all $\by\in B(\bx,\Delta)$.

The bound \eqref{eq_linear_interp_lagrange_basic} tells us that, just as for higher-order polynomial interpolation in the scalar case (\propref{prop_lebesgue_constant}), the Lebesgue constant $\max_{\by\in B(\bx,\Delta)} \|\blambda(\by)\|_1$ is a useful way to measure the quality of an interpolation model.
\revision{However, we will not be able to use the Lebesgue constant to determine how to incrementally improve an interpolation set.
Instead, we consider $\|\blambda(\by)\|_{\infty}$ in the following definition, }
which will apply equally to quadratic interpolation as well as linear.

\begin{definition} \label{def_lebesgue_lambda_poised}
    Suppose we \minorrev{have} an interpolation set $\mathcal{Y} := \{\by_1,\ldots,\by_p\}\subset\R^n$ such that its Lagrange polynomials exist. 
    Given $\bx\in\R^n$ and $\Delta>0$, we say that $\mathcal{Y}$ is $\Lambda$-poised in $B(\bx,\Delta)$ for some $\Lambda>0$ if
    \begin{align}
        \revision{\max_{\by\in B(\bx,\Delta)} \|\blambda(\by)\|_{\infty} \leq \Lambda.}
    \end{align}
\end{definition}

There are \minorrev{two} noteworthy observations that can be made at this point:
\begin{itemize}
    \item We can construct the Lagrange polynomials (i.e.~$\Lambda$ exists) if and only if the relevant interpolation linear system is invertible;
    \item Since the Lagrange polynomials always form a partition of unity \eqref{eq_lagrange_partition_unity} we have $\|\blambda(\by)\|_1 \geq \sum_{i=1}^{p} \ell_i(\by) = 1$, and so \minorrev{$\Lambda \geq 1/p$}.
    Moreover, if any $\by_i\in B(\bx,\Delta)$, which is usually the case, the condition $\ell_i(\by_i)=1$ implies \minorrev{$\Lambda \geq 1$}.
\end{itemize}

\minorrev{For example, the set} $\mathcal{Y} = \{\bx,\bx+\Delta\be_1, \ldots, \bx+\Delta \be_n\}$ from \corref{cor_fully_linear_wcc} \minorrev{has} Lagrange polynomials 
\begin{align}
    \ell_1(\by) = 1 - \be^T (\by-\bx)/\Delta, \qquad \text{and} \qquad \ell_{i+1}(\by) = \be_i^T (\by-\bx)/\Delta, \quad \forall i=1,\ldots,n.
\end{align}
The maximizers of $|\ell_i(\by)|$ over $\by\in B(\bx,\Delta)$ are $\by=\bx - \Delta \frac{\be}{\|\be\|}$ for $i=1$ and $\by=\bx\pm \Delta \be_{i-1}$ for $i=2,\ldots,n+1$, yielding $\Lambda=\sqrt{n}+1$.

As a direct consequence of \eqref{eq_linear_interp_lagrange_basic}, we get the following alternative to \thmref{thm_fully_linear}.

\begin{theorem} \label{thm_fully_linear_lagrange}
    Suppose $f$ satisfies \minorrev{\assref{ass_smoothness_1}\ref{ass_smoothness_1_smooth}} and we construct a linear model \eqref{eq_linear_model_generic} for $f$ by solving \eqref{eq_linear_interp_system_scaled}, where we assume $\hat{\bM}$ is invertible.
    If $\|\by_i-\bx\| \leq \beta \Delta$ for some $\beta>0$ and all $i=1,\ldots,p$ (where $p=n+1$ is the number of interpolation points) and the interpolation set is $\Lambda$-poised in $B(\bx,\Delta)$, then the model is fully linear in $B(\bx,\Delta)$ with constants
    \begin{align}
        \kappamf = \frac{\Lgrad}{2} \beta^2 \minorrev{p \Lambda} + \frac{\Lgrad}{2}, \quad \text{and} \quad \kappamg = 2\kappamf.
    \end{align}
\end{theorem}
\begin{proof}
    Follows immediately from \eqref{eq_linear_interp_lagrange_basic}, \minorrev{$\|\blambda(\by)\|_1 \leq p \Lambda$} and \lemref{lem_fully_linear_quadratic_from_taylor_approx}(a) with $\kappa_H=0$.
\end{proof}

We get similar results in the case of (fully) quadratic interpolation \eqref{eq_quadratic_model_generic}.
If we have interpolation points $\by_1,\ldots,\by_p$ with $p=(n+1)(n+2)/2$, then the associated Lagrange polynomials are
\begin{align}
    \ell_i(\by) := c_i + \bg_i^T (\by-\bx) + \frac{1}{2} (\by-\bx)^T \bH_i (\by-\bx), \qquad \forall i=1,\ldots,p, \label{eq_lagrange_quadratic}
\end{align}
which we can construct by solving the same systems as \eqref{eq_quadratic_interp_system} and \eqref{eq_quadratic_interp_system_scaled}, namely
\begin{align}
    \bQ \begin{bmatrix} c_i \\ \bg_i \\ \upper(\bH_i) \end{bmatrix} = \be_i, \qquad \text{or} \qquad \hat{\bQ} \begin{bmatrix} c_i \\ \Delta \: \bg_i \\ \Delta^2 \upper(\bH_i) \end{bmatrix} = \be_i, \qquad \forall i=1,\ldots,p. \label{eq_lagrange_fully_quad_system}
\end{align}
This again gives us the relationship \eqref{eq_lagrange_interpolant} and similar reasoning also gives \eqref{eq_lagrange_partition_unity}.
Recalling the natural quadratic basis \eqref{eq_natural_quadratic_basis} and \eqref{eq_quad_model_natural_basis}, we also have
\begin{align}
    \ell_i(\by) = \bphi(\by-\bx)^T \begin{bmatrix} c_i \\ \bg_i \\ \upper(\bH_i) \end{bmatrix} = \bphi(\hat{\bs})^T \begin{bmatrix} c_i \\ \Delta \: \bg_i \\ \Delta^2 \upper(\bH_i) \end{bmatrix},
\end{align}
where $\hat{\bs} := (\by-\bx)/\Delta$, and combining with \eqref{eq_lagrange_fully_quad_system} we get the analog of \eqref{eq_lagrange_system_linear}, namely
\begin{align}
    \blambda(\by) = \bQ^{-T} \bphi(\by-\bx) = \hat{\bQ}^{-T} \bphi(\hat{\bs}). \label{eq_lagrange_system_quadratic}
\end{align}
Comparing with the proof of \thmref{thm_fully_quadratic}, we see that $\bv$ in \eqref{eq_fully_quadratic_tmp2} is equal to $\blambda(\by)$, and so \eqref{eq_fully_quadratic_tmp3} can be replaced with
\begin{align}
    \left|m(\by) - f(\bx) - \grad f(\bx)^T (\by-\bx) - \frac{1}{2}(\by-\bx)^T \grad^2 f(\bx) (\by-\bx)\right| \leq \frac{\Lhess}{6}  \beta^3 \|\blambda(\by)\|_1 \Delta^3, \label{eq_fully_quadratic_tmp4}
\end{align}
for all $\by\in B(\bx,\Delta)$.

Hence we get the following version of \thmref{thm_fully_quadratic}.

\begin{theorem} \label{thm_fully_quadratic_lagrange}
    Suppose $f$ satisfies \minorrev{\assref{ass_smoothness_2}\ref{ass_smoothness_2_smooth}} and we construct a quadratic model $m$ \eqref{eq_quadratic_model_generic} for $f$ by solving \eqref{eq_quadratic_interp_system_scaled}, where we assume $\hat{\bQ}$ is invertible.
    If $\|\by_i-\bx\| \leq \beta \Delta$ for some $\beta>0$ and all $i=1,\ldots,p$ (where $p=(n+1)(n+2)/2$ is the number of interpolation points) and the interpolation set is \minorrev{$\Lambda$}-poised in $B(\bx,\Delta)$, then the model is fully quadratic in $B(\bx,\Delta)$ with constants 
    \begin{align}
        \kappamf = \frac{1}{6} \Lhess \beta^3 \minorrev{p \Lambda} + \frac{\Lhess}{6}, \quad \kappamg = \frac{17}{3} \Lhess  \beta^3 \minorrev{p \Lambda} + \frac{\Lhess}{2}, \quad \text{and} \quad \kappamh = 4 \Lhess  \beta^3 \minorrev{p \Lambda} + \Lhess.
    \end{align}
\end{theorem}
\begin{proof}
    Combine \eqref{eq_fully_quadratic_tmp4} with \minorrev{$\|\blambda(\by)\|_1 \leq p \Lambda$} and \lemref{lem_fully_linear_quadratic_from_taylor_approx}(b).
\end{proof}

\minorrev{Considering again} the structured interpolation points \eqref{eq_fully_quadratic_example_points} for fully quadratic interpolation, in \appref{app_fully_quadratic_lebesgue} we show that, for this set, \minorrev{$\Lambda=\bigO(n)$}.


Lastly, consider the case of minimum Frobenius norm quadratic interpolation, for an interpolation set $\by_1,\ldots,\by_p$ with $p\in\{n+2, \ldots, (n+1)(n+2)/2-1\}$.
Here, the associated Lagrange polynomials are again quadratic, \eqref{eq_lagrange_quadratic}, but are constructed by solving (c.f.~\eqref{eq_min_frob_problem})
\begin{align}
    \min_{c_i,\bg_i,\bH_i} \frac{1}{4}\|\bH_i\|_F^2, \qquad \text{s.t.} \quad \text{$\ell_i(\by_j) = \delta_{i,j}$ for all $j=1,\ldots,p$} \quad \text{and} \quad \bH_i = \bH_i^T,
\end{align}
which gives us (c.f.~\eqref{eq_min_frob_interp_system} and \eqref{eq_min_frob_interp_system_scaled})
\begin{align}
   \bF \left[\begin{array}{c} \lambda_{i,1} \\ \vdots \\ \lambda_{i,p} \\ \hline c_i \\ \bg_i \end{array}\right] = \be_i, \qquad \text{or} \qquad \hat{\bF} \left[\begin{array}{c} \Delta^4 \lambda_{i,1} \\ \vdots \\ \Delta^4 \lambda_{i,p} \\ \hline c_i \\ \Delta\: \bg_i \end{array}\right] = \be_i, \qquad \forall i=1,\ldots,p,
\end{align}
with associated Hessians $\bH_i = \sum_{j=1}^{p} \lambda_{i,j} (\by_j-\bx) (\by_j-\bx)^T$.
Once again, we get \eqref{eq_lagrange_interpolant} and \eqref{eq_lagrange_partition_unity}.

We now define $\bvarphi(\by)\in\R^{p+n+1}$ to be the polynomial vector
\begin{align}
    \bvarphi(\by) = \begin{bmatrix}\frac{1}{2}[(\by_1-\bx)^T (\by-\bx)]^2 & \cdots & \frac{1}{2}[(\by_p-\bx)^T (\by-\bx)]^2 & 1 & (\by-\bx)^T \end{bmatrix} ^T, \label{eq_min_frob_cols}
\end{align}
defined so that the $i$-th row/column of $\bF$ is $\bvarphi(\by_i)$, i.e.~$\bvarphi(\by_i) = \bF\be_i$ for $i=1,\ldots,p$.
We also define $\hat{\bvarphi}(\hat{\bs})$ analogously, replacing $\by_i-\bx$ with $\hat{\bs}_i = (\by_i-\bx)/\Delta$ and $\by-\bx$ with $\hat{\bs}=(\by-\bx)/\Delta$, so $\hat{\bvarphi}(\hat{\bs}_i) = \hat{\bF} \be_i$.
We can then observe that
\begin{align}
    \ell_i(\by) = c_i + \bg_i^T (\by-\bx) + \frac{1}{2} \sum_{j=1}^{p} \lambda_{i,j} [(\by_j-\bx)^T (\by-\bx)]^2 = \bvarphi(\by)^T \begin{bmatrix} \begin{array}{c} \lambda_{i,1} \\ \vdots \\ \lambda_{i,p} \\ \hline c_i \\ \bg_i \end{array} \end{bmatrix} = \bvarphi(\by)^T \bF^{-1} \be_i,
\end{align}
or $\ell_i(\by) = \hat{\bvarphi}(\hat{\bs})^T \hat{\bF}^{-1} \be_i$, and, recalling that $\bF$ and $\hat{\bF}$ are symmetric, once again we can evaluate all Lagrange polynomials at a single point without constructing any $\ell_i$ explicitly, via
\begin{align}
    \blambda(\by) = [\bF^{-1} \bvarphi(\by)]_{1,\ldots,p} = [\hat{\bF}^{-1} \hat{\bvarphi}(\hat{\bs})]_{1,\ldots,p}, \label{eq_lagrange_system_min_frob}
\end{align}
where $[\cdot]_{1,\ldots,p}$ refers to the first $p$ entries of the vector.

We now get our new version of \lemref{lem_min_frob_bounded_hessian}.

\begin{lemma} \label{lem_min_frob_bounded_hessian_lagrange}
    Suppose $f$ satisfies \minorrev{\assref{ass_smoothness_1}\ref{ass_smoothness_1_smooth}} and we construct a quadratic model $m$ \eqref{eq_quadratic_model_generic} for $f$ by solving \eqref{eq_min_frob_interp_system_scaled}, where we assume $\hat{\bF}$ is invertible.
    If $\|\by_i-\bx\| \leq \beta \Delta$ for some $\beta>0$ and all $i=1,\ldots,p$ and the interpolation set is $\minorrev{\Lambda}$-poised in $B(\bx,\Delta)$, then the model Hessian satisfies
    \begin{align}
        \|\bH\| \leq \kappa_H := 12 \Lgrad p \beta^2 \minorrev{\Lambda}. \label{eq_min_frob_hess_bound_poised}
    \end{align}
\end{lemma}
\begin{proof}
    Just as in the proof of \lemref{lem_min_frob_bounded_hessian}, we note that we get the same Hessian if we interpolate to $\tilde{f}(\by) := f(\by) - f(\bx) - \grad f(\bx)^T (\by-\bx)$, with $|\tilde{f}(\by_i)| \leq \frac{\Lgrad}{2} \beta^2 \Delta^2$.
    Hence from \eqref{eq_lagrange_interpolant} we have
    \begin{align}
        \|\bH\| \leq \sum_{i=1}^{p} |\tilde{f}(\by_i)| \cdot \|\bH_i\| \leq \frac{\Lgrad}{2} \beta^2 \Delta^2 \sum_{i=1}^{p} \|\bH_i\|,
    \end{align}
    where $\bH_i$ is the Hessian of $\ell_i(\by)$.
    Since $|\ell_i(\by)| \leq \minorrev{\Lambda}$ for all $\by\in B(\bx,\Delta)$, we can apply \lemref{lem_difference_of_quadratic} to compare $\ell_i$ with the zero function and conclude $\|\bH_i\| \leq \frac{24 \minorrev{\Lambda}}{ \Delta^2}$, from which the result follows.
\end{proof}

Our new version of  \thmref{thm_min_frob_fully_linear} is the following.

\begin{theorem} \label{thm_min_frob_fully_linear_lagrange}
    Suppose $f$ satisfies \minorrev{\assref{ass_smoothness_1}\ref{ass_smoothness_1_smooth}} and we construct a quadratic model $m$ \eqref{eq_quadratic_model_generic} for $f$ by solving \eqref{eq_min_frob_interp_system_scaled}, where we assume $\hat{\bF}$ is invertible.
    If $\|\by_i-\bx\| \leq \beta \Delta$ for some $\beta>0$ and all $i=1,\ldots,p$ and the interpolation set is \minorrev{$\Lambda$}-poised in $B(\bx,\Delta)$, then the model is fully linear with constants
    \begin{align}
        \kappamf = \frac{\Lgrad + \kappa_H}{2} \beta^2 \minorrev{p \Lambda} + \frac{\Lgrad + \kappa_H}{2}, \qquad \text{and} \qquad \kappamg = 2\kappamf + 2\kappa_H,
    \end{align}
    using the value of $\kappa_H$ defined in \lemref{lem_min_frob_bounded_hessian_lagrange}.
\end{theorem}
\begin{proof}
    By following the same argument as used to derive \eqref{eq_min_frob_fully_linear_tmp1} in the proof of \thmref{thm_min_frob_fully_linear}, we get
    \begin{align}
        |c + \bg^T (\by-\bx) - f(\bx) - \grad f(\bx)^T (\by-\bx)| \leq \frac{\Lgrad + \kappa_H}{2} \beta^2 \|\bv\|_1 \Delta^2, \label{eq_min_frob_lagrange_tmp1}
    \end{align}
    for any $\by\in B(\bx,\Delta)$, where $\bv\in\R^p$ is any vector satisfying
    \begin{align}
        \bM^T \bv = \begin{bmatrix} 1 \\ \by-\bx \end{bmatrix}. \label{eq_min_frob_v_needed}
    \end{align}
    We now show that $\blambda(\by)$ satisfies  \eqref{eq_min_frob_v_needed}.
    Since $\bH=\bm{0}$ is a global minimizer of \eqref{eq_min_frob_problem}, if $f$ is linear then we have exact interpolation, $m(\by)=f(\by)$ for all $\by$.
    Applying this to $f(\by)=(\by-\bx)^T \be_i$ for $i=1,\ldots,n$ and writing the resulting models in the form \eqref{eq_lagrange_interpolant}, we see
    \begin{align}
        (\by-\bx)^T \be_i = \sum_{i=1}^{p} (\by_i-\bx)^T \be_i \: \ell_i(\by), \label{eq_min_frob_fully_linear_lagrange_tmp1}
    \end{align}
    which, together with \eqref{eq_lagrange_partition_unity}, is equivalent to $\blambda(\by)$ satisfying \eqref{eq_min_frob_v_needed}.
    Hence, \eqref{eq_min_frob_lagrange_tmp1} becomes
    \begin{align}
        |c + \bg^T (\by-\bx) - f(\bx) - \grad f(\bx)^T (\by-\bx)| \leq \frac{\Lgrad + \kappa_H}{2} \beta^2 \|\blambda(\by)\|_1 \Delta^2, \label{eq_min_frob_fully_linear_lagrange_tmp2}
    \end{align}
    and so
    \begin{align}
        |m(\by) - f(\bx) - \grad f(\bx)^T (\by-\bx)| &\leq \frac{\Lgrad + \kappa_H}{2} \beta^2 \|\blambda(\by)\|_1 \Delta^2 + \frac{1}{2}|(\by-\bx)^T \bH (\by-\bx)|, \\
        &\leq \left(\frac{\Lgrad + \kappa_H}{2} \beta^2 \|\blambda(\by)\|_1  + \frac{\kappa_H}{2}\right) \Delta^2. \label{eq_min_frob_fully_linear_lagrange_tmp3}
    \end{align}
    The result then follows from \lemref{lem_fully_linear_quadratic_from_taylor_approx}(a) \minorrev{and $\|\blambda(\by)\|_1 \leq p \Lambda$}.
\end{proof}

\minorrev{Now consider} the set $\mathcal{Y} = \{\bx, \bx+\Delta\be_1, \ldots, \bx+\Delta \be_n, \bx-\Delta \be_1, \ldots, \bx-\Delta \be_n\}$, as in \corref{cor_wcc_dfo_min_frob}.
The associated Lagrange polynomials are:
\begin{subequations} \label{eq_min_frob_example_lag_poly}
\begin{align}
    \ell_1(\by) &= 1 - \frac{1}{\Delta^2} \|\by-\bx\|^2, \\
    \ell_{i+1}(\by) &= \frac{(y_i-x_i)^2}{2\Delta^2} + \frac{y_i-x_i}{2\Delta}, \\
    \ell_{n+i+i}(\by) &= \frac{(y_i-x_i)^2}{2\Delta^2} - \frac{y_i-x_i}{2\Delta}, 
\end{align}
\end{subequations}
for $i=1,\ldots,n$.
Each of these takes values in $[-\frac{1}{8},1]$ for all $\by\in B(\bx,\Delta)$ and so \minorrev{$\Lambda=1$}.

\begin{remark}
    For the example interpolation sets considered here, applying the value of $\Lambda$ in the fully linear/\minorrev{fully} quadratic constants gives a worse bound (in terms of $n$) than \secref{sec_model_construction}.
    These can be improved by calculating the \minorrev{constants} in terms of the Lebesgue constant, but this does not allow for the procedures in \secref{sec_model_improvement} to be used.
\end{remark}

\subsection{Model Improvement} \label{sec_model_improvement}
As mentioned previously, the main reason for using Lagrange polynomials in our model error bounds is that it allows us \revision{both check if a model is fully linear, and to decide how to make small changes to the interpolation set to make it fully linear.} 
Essentially, we wish to find an interpolation set for which \minorrev{$\Lambda$} is sufficiently small, below some given value.

Firstly, we can determine \minorrev{$\Lambda$} by maximizing each Lagrange polynomial and its negative within the trust-region.
Since the Lagrange polynomials are themselves linear or quadratic, this reduces to solving $2p$ trust-region subproblems \eqref{eq_trs_generic}.
In theory, this requires finding the global trust-region solution for each, but in practice solving these trust-region subproblems approximately can give a sufficiently good estimate \cite{Powell2002}.

Now, suppose we have an interpolation set which is \minorrev{$\Lambda$}-poised, but we want to change some points to reduce its value of \minorrev{$\Lambda$} below some threshold \minorrev{$\Lambda^*>1$} (recalling that \minorrev{$\Lambda \geq 1$} if any $\by_i \in B(\bx,\Delta)$).
The procedure for changing the interpolation set to improve its poisedness constant is quite simple, and given in \algref{alg_geom_improvement}.
At each iteration of \algref{alg_geom_improvement} we find $\by\in B(\bx,\Delta)$ and $i$ such that $|\lambda_i(\by)|=\minorrev{\Lambda}$ and replace $\by_i$ with $\by$.
We do not specify a particular interpolation model type; \algref{alg_geom_improvement} applies equally to linear, quadratic or minimum Frobenius norm quadratic interpolation without modification.

\begin{algorithm}[tb]
\begin{algorithmic}[1]
\Require Interpolation set $\mathcal{Y}$ which is \minorrev{$\Lambda$}-poised in $B(\bx,\Delta)$ for some \minorrev{$\Lambda>0$} and all $\|\by_i-\bx\| \leq \beta \Delta$ for some $\beta\geq 1$, desired \minorrev{$\Lambda$}-poisedness constant $\minorrev{\Lambda^*} > 1$.
\While{\minorrev{$\Lambda > \Lambda^*$}}
    \State Find $\by\in B(\bx,\Delta)$ and $i\in\{1,\ldots,p\}$ such that $|\ell_i(\by)| > \minorrev{\Lambda^*}$.
    \State Update $\mathcal{Y}$ by replacing $\by_i$ with $\by$. \label{ln_geom_improvement_update}
    \State Recompute the poisedness constant \minorrev{$\Lambda$} of the new $\mathcal{Y}$.
\EndWhile
\State \Return $\mathcal{Y}$
\end{algorithmic}
\caption{Interpolation set improvement using Lagrange polynomials.}
\label{alg_geom_improvement}
\end{algorithm}

To analyze \algref{alg_geom_improvement}, we consider the impact on the determinant of the relevant interpolation matrix (i.e.~$\bM$ \eqref{eq_linear_interp_system}, $\bQ$ \eqref{eq_quadratic_interp_system} or $\bF$ \eqref{eq_min_frob_interp_system} depending on the type of interpolation model) from changing a single point in the interpolation set.

The case for linear or quadratic interpolation is quite straightforward.

\begin{lemma} \label{lem_det_updating_lin_quad}
    Suppose the set $\mathcal{Y}$ is used for either linear or quadratic interpolation, with interpolation matrix $\bA\in\R^{p\times p}$ (i.e.~$\bM$ or $\bQ$), where $\mathcal{Y}=\{\by_1,\ldots,\by_p\}$. 
    If we replace $\by_i$ in $\mathcal{Y}$ with another point $\by$, then the resulting interpolation matrix $\bA_{\textnormal{new}}$ satisfies
    \begin{align}
        |\det(\bA_{\textnormal{new}})| \geq |\ell_i(\by)| \: |\det(\bA)|,
    \end{align}
    where the Lagrange polynomial $\ell_i$ is given by the original interpolation set $\mathcal{Y}$ (before the replacement is made).
\end{lemma}
\begin{proof}
    If $\bA$ is singular, the result holds trivially, so assume that $\bA$ is invertible.
    Replacing $\by_i$ with $\by$ changes the $i$th row of $\bA$ from $\bphi(\by_i)^T$ to $\bphi(\by)^T$ in $\bA_{\textnormal{new}}$, where $\bphi(\by) = \begin{bmatrix} 1 \\ \by-\bx \end{bmatrix}$ in the linear interpolation case or the natural quadratic basis \eqref{eq_natural_quadratic_basis} for quadratic interpolation.
    Equivalently, $\bA_{\textnormal{new}}^T$ is the same as $\bA^T$, but with its $i$th column changed from $\bphi(\by_i)$ to $\bphi(\by)$.
    However, from either \eqref{eq_lagrange_system_linear} or \eqref{eq_lagrange_system_quadratic}, we have that $\ell_i(\by)$ is the $i$th entry of $\bA^{-T} \bphi(\by)$.
    By Cramer's rule for linear systems, this means  $\ell_i(\by) = \det(\bA_{\textnormal{new}}^T) / \det(\bA^T) = \det(\bA_{\textnormal{new}}) / \det(\bA)$, giving the result.
\end{proof}

For minimum Frobenius norm updating, a similar result holds, but it is more complicated to prove because of the more complicated structure of the interpolation matrix $\bF$ \eqref{eq_min_frob_interp_system}.
In particular, replacing a single interpolation point changes both a row and column of $\bF$, so Cramer's rule cannot be used as in the proof of \lemref{lem_det_updating_lin_quad}. 

\begin{lemma} \label{lem_det_updating_min_frob}
    Suppose the set $\mathcal{Y}$ is used for minimum Frobenius norm interpolation, with interpolation matrix $\bF\in\R^{(p+n+1)\times(p+n+1)}$, where $\mathcal{Y} = \{\by_1,\ldots,\by_p\}$. 
    If we replace $\by_i$ in $\mathcal{Y}$ with another point $\by$, then the resulting interpolation matrix $\bF_{\textnormal{new}}$ satisfies
    \begin{align}
        |\det(\bF_{\textnormal{new}})| \geq \ell_i(\by)^2 \: |\det(\bF)|,
    \end{align}
    where the Lagrange polynomial $\ell_i$ is given by the original interpolation set $\mathcal{Y}$ (before the replacement is made).
\end{lemma}
\begin{proof}
    See \cite[Theorem 5.2]{Roberts2024}.
\end{proof}

We are now ready to show that \algref{alg_geom_improvement} achieves the desired outcome.

\begin{theorem} \label{thm_geom_improvement}
    \algref{alg_geom_improvement} terminates in finite time, and the resulting interpolation set $\mathcal{Y}$ is \minorrev{$\Lambda^*$}-poised in $B(\bx,\Delta)$, and we have $\|\by_i-\bx\| \leq \beta \Delta$ for all $\by_i\in \mathcal{Y}$.
\end{theorem}
\begin{proof}
    The termination condition for \algref{alg_geom_improvement} implies that $\mathcal{Y}$ is \minorrev{$\Lambda^*$}-poised in $B(\bx,\Delta)$ if it terminates.
    All points in the resulting $\mathcal{Y}$ are either in $B(\bx,\beta\Delta)$ (if they were in $\mathcal{Y}$ initially) or are in $B(\bx,\Delta)$ (if they come from an update in line~\ref{ln_geom_improvement_update} of \algref{alg_geom_improvement}).
    Hence all points satisfy $\|\by_i-\bx\| \leq \beta \Delta$ since we assume $\beta \geq 1$.
    It remains to show that \algref{alg_geom_improvement} terminates in finite time.

    Since the initial set $\mathcal{Y}$ is \minorrev{$\Lambda$}-poised, the corresponding interpolation matrix, say $\bA$ (either $\bM$, $\bQ$ or $\bF$ depending on the type of interpolation) is invertible, so $|\det(\bA)| > 0$ initially.
    Call this value $d_0$.
    However, since all points from the initial $\mathcal{Y}$ and any new points generated are all in the compact set $B(\bx,\beta \Delta)$, there is a finite maximum value, say $d_{\max}$, of $|\det(\bA)|$ over all interpolation sets contained in $B(\bx,\beta \Delta)$, 
    From either \lemref{lem_det_updating_lin_quad} or \lemref{lem_det_updating_min_frob}, whenever we make the update in line~\ref{ln_geom_improvement_update} of \algref{alg_geom_improvement}, we increase $|\det(\bA)|$ by a factor of at least $\min(|\ell_i(\by)|, \ell_i(\by)^2) > \minorrev{\Lambda^*}$ (since $|\ell_i(\by)| > \minorrev{\Lambda^*} > 1$).
    Hence \algref{alg_geom_improvement} must terminate after at most $\lceil \log_{\minorrev{\Lambda}^*}(d_{\max}/d_0)\rceil$ iterations.
\end{proof}

\begin{figure}[tb]
  \centering
  \begin{subfigure}[b]{0.3\textwidth}
    \includegraphics[width=\textwidth]{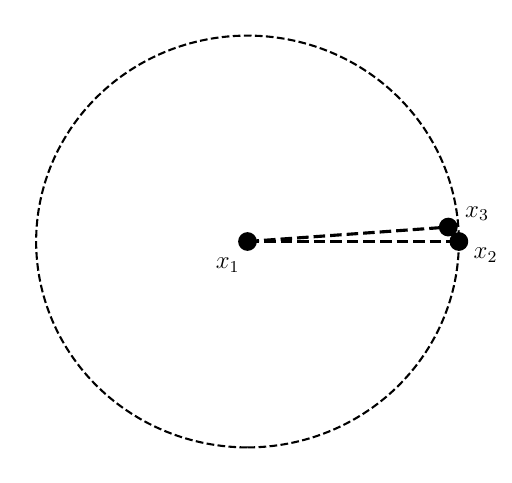}
    \caption{Initial set ($\minorrev{\Lambda} \approx 14.3$)}
    \label{fig_poisedness_demo_1}
  \end{subfigure}
  ~
  \begin{subfigure}[b]{0.3\textwidth}
    \includegraphics[width=\textwidth]{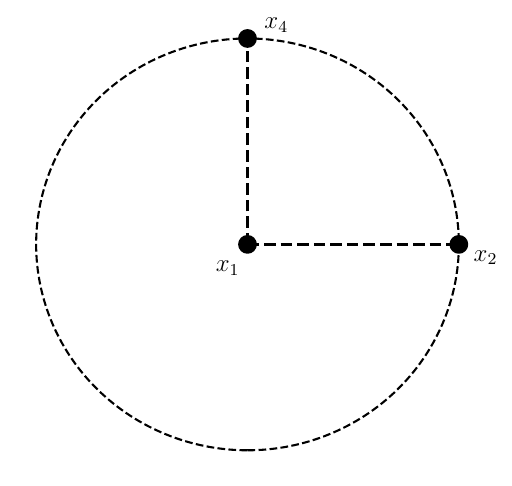}
    \caption{After 1 iteration ($\minorrev{\Lambda} \approx 2.41$)}
    \label{fig_poisedness_demo_2}
  \end{subfigure}
  ~
  \begin{subfigure}[b]{0.3\textwidth}
    \includegraphics[width=\textwidth]{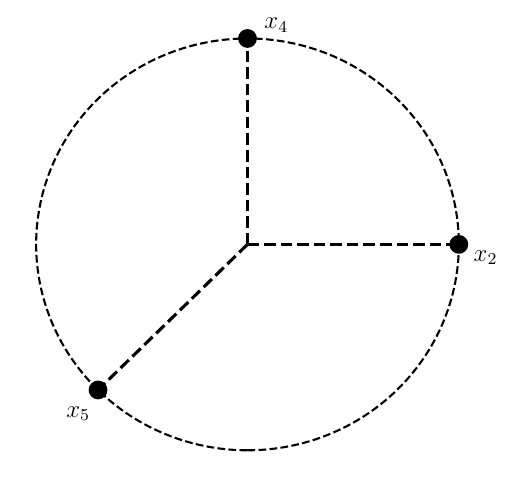}
    \caption{After 2 iterations ($\minorrev{\Lambda} \approx 1.06$)}
    \label{fig_poisedness_demo_3}
  \end{subfigure}
  \caption{Demonstration of \algref{alg_geom_improvement} for linear interpolation in $B(\bm{0},1)\subset\R^2$. (a) Initial interpolation set is $\bx_1=(0,0)$, $\bx_2=(1,0)$ and $\bx_3=(0.95,0.07)$. The value $\minorrev{\Lambda}\approx 14.3$ comes from Lagrange polynomial associated with $\bx_3$ at point $\bx_4=(0,1)$. (b) New set is $\bx_1$, $\bx_2$ and $\bx_4$. The value $\minorrev{\Lambda}\approx 2.41$ comes from Lagrange polynomial associated with $\bx_1$ at point $\bx_5 \approx (-0.707, -0.707)$. (c) New set is $\bx_2$, $\bx_4$ and $\bx_5$.}
  \label{fig_poisedness_demo}
\end{figure}

Three iterations of \algref{alg_geom_improvement} for linear interpolation in two dimensions are illustrated in \figref{fig_poisedness_demo}. 
Initially, two interpolation points $\bx_2$ and $\bx_3$ are very close, yielding a large value of $\minorrev{\Lambda}$.\footnote{The value of $\minorrev{\Lambda}$ is still not unreasonably large here, simply so the locations of $\bx_2$ and $\bx_3$ can be illustrated clearly.}
After two iterations the value of $\minorrev{\Lambda}$ has been decreased by more than an order of magnitude, and is very close to 1.
Note that the center of the region $\bx_1$ is no longer chosen as an interpolation point, so $m(\bx) \neq f(\bx)$ is possible for this final interpolation set.
However, after one iteration of \algref{alg_geom_improvement} we already have a small value of $\minorrev{\Lambda}$ while retaining $m(\bx) = f(\bx)$, which may be desirable in practice.

\subsection{Incremental Geometry Improvement Algorithm} \label{sec_geom_incr_algo}
\revision{We conclude this section by describing a variant of \algref{alg_basic_tr_dfo} which avoids the requirement that every model $m_k$ is fully linear, and instead handles non-fully linear models by incrementally improving the geometry of the current interpolation set. This core principle---only worry about the geometry of the interpolation set if the algorithm is not progressing, and only make minimal changes to the set if so---is used in many state-of-the-art MBDFO codes, such as Powell's methods (see \secref{sec_conclusion}).}


At each iteration $k$ of this algorithm, we have an interpolation set $\mathcal{Y}_k = \{\by_{k,1}, \ldots, \by_{k,p}\}$ with $\bx_k \in \mathcal{Y}_k$, which can be used for linear, quadratic or minimum Frobenius norm quadratic interpolation.
The algorithm ensures that the interpolation linear system is invertible at every iteration.
\revision{After the model is constructed and a tentative step determined, a slightly different the trust-region mechanism to \algref{alg_basic_tr_dfo} is used: if a step is not accepted and the interpolation set is not guaranteed to be fully linear, we perform a single iteration of the geometry-improving \algref{alg_geom_improvement} (and do not reduce $\Delta_k$) instead of having an unsuccessful step (where $\Delta_k$ is reduced).}

\revision{Formally, the `replace distant point' and `replace bad point' steps in \algref{alg_self_correcting} have the minimal requirements that the new geometry-improving point $\by_{\textnormal{new}}$ must satisfy, in practice a good choice---motivated by Lemma~\ref{lem_det_updating_lin_quad} or \ref{lem_det_updating_min_frob}---is to find this point by maximizing the (magnitude of the) relevant Lagrange polynomial inside the current trust region $B(\bx_k,\Delta_k)$. More details are given in  \remref{rem_incr_update_ideas} below.}

\begin{algorithm}[tb]
\begin{algorithmic}[1]
\Require Starting point $\bx_0\in\R^n$ and trust-region radius $\Delta_0>0$. Algorithm parameters: scaling factors $0 < \gammadec < 1 < \gammainc$, acceptance thresholds $0 < \eta_U \leq \eta_S < 1$, criticality threshold $\mu_c > 0$, number of interpolation points $p$, distance threshold $\beta>1$ and poisedness threshold $\minorrev{\Lambda}>1$.
\State Select an initial interpolation set $\mathcal{Y}_0 \subset \R^n$ of size $p$ with $\bx_0 \in \mathcal{Y}_0$ such that the resulting interpolation linear system is invertible.
\For{$k=0,1,2,\ldots$}
    \State Solve the relevant interpolation problem to obtain the model $m_k$ \eqref{eq_tr_model} and Lagrange polynomials $\ell_{k,1},\ldots,\ell_{k,p}$.
    \State Solve the trust-region subproblem \eqref{eq_trs} to get a step $\bs_k$ satisfying \assref{ass_cauchy_decrease}. 
    \State Evaluate $f(\bx_k+\bs_k)$ and calculate the ratio $\rho_k$ \eqref{eq_ratio_generic}.
    \If{$\rho_k \geq \eta_U$ and $\|\bg_k\| \geq \mu_c \Delta_k$} 
        \State \textit{((Very) successful iteration)} Set $\bx_{k+1}=\bx_k+\bs_k$ and either $\Delta_{k+1}=\gammainc\Delta_k$ if $\rho_k \geq \eta_S$ or $\Delta_{k+1}=\Delta_k$ if $\eta_U \leq \rho_k < \eta_S$.
        \State Set $\mathcal{Y}_{k+1}$ to be any interpolation set of size $p$ with $\bx_{k+1}\in \mathcal{Y}_{k+1}$ so that the resulting interpolation linear system is invertible. \label{ln_generic_update}
    \ElsIf{$\max_{i=1,\ldots,p} \|\by_{k,i}-\bx_k\| > \beta \Delta_k$}
        \State \textit{(Replace distant point)} Set $\bx_{k+1}=\bx_k$ and $\Delta_{k+1}=\Delta_k$.
        \State Set $\mathcal{Y}_{k+1} = \mathcal{Y}_k \setminus \{\by_{k,i_k}\} \cup \{ \by_{\textnormal{new}}\}$, where $i_k\in\{1,\ldots,p\}$ satisfies $\|\by_{k,i_k}-\bx_k\| > \beta \Delta_k$ and $\by_{\textnormal{new}} \in B(\bx_k,\Delta_k)$ satisfies $\ell_{k,i_k}(\by_{\textnormal{new}}) \neq 0$.
    \ElsIf{$\max_{\by \in B(\bx_k,\Delta_k)} \max_{i=1,\ldots,p \: \text{s.t.} \by_i \neq \bx_k} |\ell_{k,i}(\by)| > \minorrev{\Lambda}$}
        \State \textit{(Replace bad point)} Set $\bx_{k+1}=\bx_k$ and $\Delta_{k+1}=\Delta_k$.
        \State Set $\mathcal{Y}_{k+1} = \mathcal{Y}_k \setminus \{\by_{k,i_k}\} \cup \{ \by_{\textnormal{new}} \}$, where $i_k\in\{1,\ldots,p\}$ and $\by_{\textnormal{new}}\in B(\bx_k,\Delta_k)$ satisfy $|\ell_{k,i_k}(\by_{\textnormal{new}})| > \minorrev{\Lambda}$ and $\by_{k,i_k} \neq \bx_k$.
    \Else
        \State \textit{(Unsuccessful iteration)} Set $\bx_{k+1}=\bx_k$, $\Delta_{k+1}=\gammadec\Delta_k$ and $\mathcal{Y}_{k+1} = \mathcal{Y}_k$. 
    \EndIf 
\EndFor
\end{algorithmic}
\caption{\revision{Incremental geometry improving MBDFO algorithm} for \eqref{eq_problem}.}
\label{alg_self_correcting}
\end{algorithm}

The full algorithm is given in \algref{alg_self_correcting}.
\revision{It is written in a generic way, without specifying a specific interpolation type, which can be linear, fully quadratic\footnote{Note that a fully quadratic model with $\Delta \leq \Delta_{\max}$ is fully linear with constants $\kappamf \Delta_{\max}$ and $\kappamg \Delta_{\max}$, or use minimum Frobenius interpolation with $p=(n+1)(n+2)/2$; see \remref{rem_fully_quad_is_min_frob}.} or minimum Frobenius quadratic.}


We begin with some basic results.

\begin{lemma} \label{lem_basic_self_correcting}
    \minorrev{Regarding \algref{alg_self_correcting}, we have:}
    \begin{enumerate}[label=(\alph*)]
        \item The interpolation linear system is invertible at every iteration (and hence Lagrange polynomials exist at every iteration);
        \item If $\max_{i=1,\ldots,p} \|\by_{k,i}-\bx_k\| \leq \beta \Delta_k$ and $\max_{\by\in B(\bx,\Delta)} |\ell_{k,i}(\by)| \leq \minorrev{\Lambda}$ for all $i=1,\ldots,p$ except if $\by_i=\bx_k$, then  the model $m_k$ is fully linear in $B(\bx_k.\Delta_k)$ with constants $\kappamf,\kappamg>0$ possibly depending on $p$, $\beta$ and $\minorrev{\Lambda}$; \label{lem_fully_linear_self_correcting} 
        \item We can only have finitely many \minorrev{consecutive} iterations of types `replace distant point' or `replace bad point' before a (very) successful or unsuccessful iteration must occur; and \label{lem_basic_self_correcting_finite_geom}
    \end{enumerate}
\end{lemma}
\begin{proof}
    First, we show (a) by induction.
    The case $k=0$ holds by construction of $\mathcal{Y}_0$, so suppose it is invertible for some iteration $k$.
    If iteration $k$ is (very) successful, then the next linear system is invertible by definition of $\mathcal{Y}_{k+1}$.
    If iteration $k$ replaces a distant or bad interpolation point, then $\mathcal{Y}_{k+1} = \mathcal{Y}_k \setminus \{\by_{k,i_k}\} \cup \{ \by_{\textnormal{new}} \}$ for $\ell_{k,i_k}(\by_{\textnormal{new}}) \neq 0$, so the system is invertible by \minorrev{either Lemma~\ref{lem_det_updating_lin_quad} or \ref{lem_det_updating_min_frob}}. 
    Lastly, if iteration $k$ is unsuccessful, then the linear system is unchanged.

    To show (b), denote $i^*$ to be the index such that $\by_{k,i^*}=\bx_k$.
    Since the Lagrange polynomials always form a partition of unity, we have, for any $\by\in B(\bx,\Delta)$,
    \begin{align}
        |\ell_{k,i^*}(\by)| \leq 1 + \sum_{\substack{i=1 \\ i \neq i^*}}^{p} |\ell_{k,i}(\by)| \leq 1 + (p-1)\minorrev{\Lambda}.
    \end{align}
    Hence $\max_{\by\in B(\bx,\Delta)} \|\blambda(\by)\|_{\infty} \leq 1+(p-1)\minorrev{\Lambda}$, and the model $m_k$ is fully linear in $B(\bx_k,\Delta_k)$ with suitable constants by \minorrev{Theorem~\ref{thm_fully_linear_lagrange} or \ref{thm_min_frob_fully_linear_lagrange} (which also covers fully quadratic interpolation, per \remref{rem_fully_quad_is_min_frob})}. 
    \minorrev{Finally, (c) is an immediate consequence of \thmref{thm_geom_improvement}.}
\end{proof}

\revision{The properties in \lemref{lem_basic_self_correcting} are sufficient to prove global convergence for \algref{alg_self_correcting} (i.e.~$\liminf_{k\to\infty} \|\grad f(\bx_k)\| = 0$), but to get a worst-case complexity bound, we need to slightly strengthen \lemref{lem_basic_self_correcting}\ref{lem_basic_self_correcting_finite_geom} to a uniform bound.

\begin{assumption} \label{ass_bounded_geom_steps}
    The number of consecutive iterations of types `replace distant point' or `replace bad point' cannot exceed some $G_{\max}\in\N$, independent of the starting interpolation set.
\end{assumption}

It was recently shown that \assref{ass_bounded_geom_steps} holds with $G_{\max} = \bigO(p \log n)$ for linear and fully quadratic interpolation \cite[Theorem 5.1]{Chaudhry2025}.
No such bound is known for minimum Frobenius quadratic interpolation, but it is likely that one exists.
}


\begin{lemma} \label{lem_very_successful_self_correcting}
    Suppose Assumptions~\minorrev{\ref{ass_smoothness_1}\ref{ass_smoothness_1_smooth}},  \minorrev{\ref{ass_cauchy_decrease}, and \ref{ass_model_dfo}\ref{ass_model_dfo_H}} hold.
    \minorrev{If, on iteration $k$ of \algref{alg_self_correcting},} $\max_{i=1,\ldots,p} \|\by_{k,i}-\bx_k\| \leq \beta \Delta_k$ and $\max_{\by\in B(\bx,\Delta)} |\ell_{k,i}(\by)| \leq \minorrev{\Lambda}$ for all $i=1,\ldots,p$ except if $\by_{k,i}=\bx_k$, and $\bg_k\neq\bm{0}$ and
    \begin{align}
        \Delta_k \leq \min\left(\frac{\kappa_s (1-\eta_S)}{2\kappamf}, \frac{1}{\kappa_H}, \frac{1}{\mu_c}\right) \|\bg_k\|,
    \end{align}
    where $\kappamf$ is from \lemref{lem_basic_self_correcting}\ref{lem_fully_linear_self_correcting}, then $\rho_k \geq \eta_S$ and $\|\bg_k\| \geq \mu_c \Delta_k$ (i.e.~iteration $k$ is very successful).
\end{lemma}
\begin{proof}
    By \lemref{lem_basic_self_correcting}\ref{lem_fully_linear_self_correcting}, the model $m_k$ is fully linear in $B(\bx_k,\Delta_k)$ and so the proof of \lemref{lem_very_successful_dfo} holds.
\end{proof}

\begin{lemma} \label{lem_delta_min_self_correcting}
    Suppose Assumptions~\minorrev{\ref{ass_smoothness_1}\ref{ass_smoothness_1_smooth}}, \minorrev{\ref{ass_cauchy_decrease}, and \ref{ass_model_dfo}\ref{ass_model_dfo_H}} hold and we run \algref{alg_self_correcting}.
    If $\|\grad f(\bx_k)\| \geq \epsilon$ for all $k=0,\ldots,K-1$, then $\Delta_k \geq \Delta_{\min}(\epsilon)$ for all $k=0,\ldots,K$, where $\Delta_{\min}(\epsilon)$ is defined in \lemref{lem_delta_min_dfo}.
\end{lemma}
\begin{proof}
    The proof is identical to that of \lemref{lem_delta_min_dfo}, except noting that if $\Delta_{k+1} < \Delta_k$ then we must have an unsuccessful iteration, and so $\max_{i=1,\ldots,p} \|\by_{k,i}-\bx_k\| \leq \beta \Delta_k$ and $\max_{\by\in B(\bx,\Delta)} |\ell_{k,i}(\by)| \leq \minorrev{\Lambda}$ for all $i=1,\ldots,p$ except if $\by_{k,i}=\bx_k$.
    This allows us to use \lemref{lem_very_successful_self_correcting} to achieve the contradiction.
\end{proof}

\revision{Our final worst-case complexity result is the following.

\begin{theorem} \label{thm_wcc_self_correcting}
    Suppose Assumptions~\ref{ass_smoothness_1}, \minorrev{\ref{ass_cauchy_decrease},  \ref{ass_model_dfo}\ref{ass_model_dfo_H}} and \ref{ass_bounded_geom_steps} hold. 
    If $k_{\epsilon}$ is the first iteration of \algref{alg_self_correcting} such that $\|\grad f(\bx_{k_\epsilon})\| < \epsilon$, then $k_{\epsilon} = \bigO( G_{\max} \kappa_H (\kappam + \kappa_H)^2 \epsilon^{-2})$, also requiring a total of $\bigO( G_{\max} \kappa_H (\kappam + \kappa_H)^2 \epsilon^{-2})$ objective evaluations.
    Hence $\liminf_{k\to\infty} \|\grad f(\bx_k)\| = 0$.
\end{theorem}
\begin{proof}
    The bounds on the number of (very) successful \eqref{eq_wcc_dfo_tmp1} and unsuccessful \eqref{eq_wcc_tmp2} iterations from the proof of \thmref{thm_wcc_dfo} hold here, where we only need to note that $\Delta_k$ is unchanged for `replace distant point' and `replace bad point' iterations to recover \eqref{eq_wcc_tmp2}.
    Lastly, the total number of `replace distant point' and `replace bad point' iterations is at most $G_{\max}(|\mathcal{S}| + |\mathcal{U}|)$ by \assref{ass_bounded_geom_steps}, so the total number of iterations of all types is at most $(G_{\max}+1)(|\mathcal{S}| + |\mathcal{U}|)$, i.e.~a factor of $\bigO(G_{\max})$ worse than \thmref{thm_wcc_dfo}.
    The objective evaluation bound comes from noting that \algref{alg_self_correcting} requires at most 2 evaluations per iteration (once the first interpolation set is chosen). 
\end{proof}

If we could construct linear interpolation models to give $\kappam$ independent of $n$ as in \corref{cor_fully_linear_wcc_improved}, \thmref{thm_wcc_self_correcting} with $G_{\max} = \bigO(n\log n)$ gives an evaluation complexity bound that is slightly worse than \corref{cor_fully_linear_wcc_improved} by a factor of $\bigO(\log n)$.
We can both make $\kappam$ independent of $n$ and remove the $\bigO(\log n)$ factor from the complexity bound by \emph{excluding} the Lagrange polynomial associated with the current iterate $\bx_k$ in the $\Lambda$-poisedness requirement (i.e.~by modifying \defref{def_lebesgue_lambda_poised}) \cite{Chaudhry2025}.
In general, if the interpolation set is updated well on (very) successful iterations---see \remref{rem_incr_update_ideas}---then $\kappam$ and $G_{\max}$ are usually small and the approach of \algref{alg_self_correcting} is efficient in practice.
}

\begin{remark} \label{rem_incr_update_ideas}
On (very) successful iterations, we have a flexible choice in how to update the \minorrev{interpolation set} (line~\ref{ln_generic_update} of \algref{alg_self_correcting}).
A simple way to do this would be $\mathcal{Y}_{k+1} = \mathcal{Y}_k \setminus \{\by_{k,i_k}\} \cup \{\bx_{k+1}\}$, where $i_k$ is any index such that $\ell_{k,i_k}(\bx_{k+1}) \neq 0$, which must exist by \eqref{eq_lagrange_partition_unity}.
Usually, all Lagrange polynomials will be nonzero at $\bx_{k+1}$, and so, motivated by the geometry improving procedure, a good choice might be to replace a point that is far from $\bx_{k+1}$ and/or improves the poisedness of the interpolation set (i.e.~$|\ell_{k,i_k}(\bx_{k+1})|$ is large, as in \algref{alg_geom_improvement}).
\revision{For example, in \cite{Powell2002}, $i_k$ is selected as}
\begin{align}
    i_k \minorrev{\in} \argmax_{i=1,\ldots,p} \left\{\max\left(\frac{\|\by_{k,i}-\bx_{k+1}\|^3}{\Delta_k^3}, 1\right) \: |\ell_{k,i}(\bx_{k+1})|\right\}.
\end{align}
For `replace distant point' and `replace bad point' iterations, typically the point removed from the interpolation set is the maximizer of the relevant quantity (distance from $\bx_k$ or poisedness), and is replaced with a point that maximizes the magnitude of the Lagrange polynomial of the point being removed, as in \algref{alg_geom_improvement}.
\end{remark}

\revision{\paragraph{Powell's Methods} The framework of \algref{alg_self_correcting} broadly aligns with the approach in Powell's highly regarded \renaming{MBDFO} software (see \secref{sec_conclusion}). The most notable distinction is that two trust-region radii are used in a complex way that partially decouples the step size constraint $\|\bs_k\| \leq \Delta_k$ from the radius used to assess $\Lambda$-poisedness. It also uses a simplified interpolation set management procedure:
\begin{itemize}
    \item On (very) successful steps, update the interpolation set as per \remref{rem_incr_update_ideas};
    \item Otherwise, check that all interpolation points are sufficiently close to $\bx_k$. If yes, decrease both trust-region radii. If no, replace the furthest interpolation point with the approximate maximizer of (the magnitude of) the associated Lagrange polynomial, and possibly decrease one trust-region radius. 
\end{itemize}
Additionally, the implementation pays careful attention to efficient linear algebra and subproblem solvers. This approach has proven extremely successful in practice, despite the simplified interpolation set management having limited convergence guarantees (e.g.~\cite{Han2004}).\footnote{\revision{The use of two trust-region radii on its own still allows for convergence theory, e.g.~\cite{Zhang2010}.}} An accessible overview of these algorithms can be found in \cite{Ragonneau2024}.}

\subsubsection*{Notes and References}
{\small More information on Lebesgue constants in approximation theory may be found in resources such as \cite{Trefethen2020}. The theory of Lagrange polynomials to construct fully linear/\minorrev{fully} quadratic interpolation models was originally developed in \cite{Conn2007} for linear/quadratic interpolation and \cite{Conn2008} for minimum Frobenius norm interpolation, with these results collected in \cite{Conn2009}. We again note that minimum Frobenius norm quadratic models were originally studied in \cite{Powell2004}. 

A more concise alternative to \eqref{eq_linear_interp_lagrange_basic} and \eqref{eq_fully_quadratic_tmp4} are the bounds \cite{Powell2001,Cao2023}
\begin{align}
    |m(\by) - f(\by)| \leq \frac{\Lgrad}{2} \sum_{i=1}^{p} |\ell_i(\by)| \: \|\by-\by_i\|^2, \quad \text{and} \quad 
    |m(\by) - f(\by)| \leq \frac{\Lhess}{6} \sum_{i=1}^{p} |\ell_i(\by)| \: \|\by-\by_i\|^3, \label{eq_linear_error_lagrange_simple}
\end{align}
for all $\by\in\R^n$, for linear and (fully) quadratic interpolation under Assumptions~\minorrev{\ref{ass_smoothness_1}\ref{ass_smoothness_1_smooth}} and \minorrev{\ref{ass_smoothness_2}\ref{ass_smoothness_2_smooth}} respectively.

The \minorrev{incremental geometry improving} algorithm (\algref{alg_self_correcting}) is \minorrev{based on} the approach from \cite{Scheinberg2010}, which also includes an example where simply updating $\mathcal{Y}_k$ by removing the furthest point from $\bx_k$ and replacing it with $\bx_k+\bs_k$ does not converge, demonstrating that consideration of interpolation set geometry is necessary for a convergent \renaming{MBDFO} method.
\revision{The presentation in \secref{sec_geom_incr_algo} draws on the complexity bounds from the more recent analysis \cite{Chaudhry2025}, which also proves that \assref{ass_bounded_geom_steps} can be satisfied for polynomial interpolation over an arbitrary basis.}


{\small Recent work on Lagrange polynomials in \renaming{MBDFO} have included refining the error bounds from Lagrange polynomials \cite{Cao2023}, extending to the constrained optimization case \cite{Hough2022,Roberts2024} (see \secref{sec_convex_constraints}), and drawing connections with methods for outlier detection \cite{Zhang2024}.


\section{Constrained Optimization} \label{sec_constraints}

In this section we consider how the above algorithmic ideas and approximation theory can be adapted to constrained optimization problems.
We first consider simple convex constraints (not involving any derivative-free functions), such as bounds, and then consider the case of general, potentially derivative-free, constraints.

\subsection{Simple Convex Constraints} \label{sec_convex_constraints}

We first consider the case of \renaming{MBDFO} in the presence of simple constraints on the decision variables.
We will focus on the case where the feasible set is an easy-to-describe, convex set.
Specifically, we will extend the above ideas to solve problems of the form
\begin{align}
    \min_{\bx\in\C} f(\bx), \label{eq_convex_cons_generic}
\end{align}
where $f:\R^n\to\R$ is, as usual, smooth and nonconvex but with only a zeroth order oracle available, and the feasible set $\C$ satisfies the following.

\begin{assumption} \label{ass_convex_cons}
    The feasible set $\C\subseteq \R^n$ in \eqref{eq_convex_cons_generic} is closed, convex and has nonempty interior.
\end{assumption}

This covers some of the most important constraint types, such as lower/upper bounds and linear inequality constraints. 
To allow us to apply our method to many choices of feasible sets, the only information about $\C$ available to our algorithm will be its Euclidean projection operator,
\begin{align}
    \proj_{\C}(\bx) := \argmin_{\by\in \C} \|\by-\bx\|,
\end{align}
which returns the closest point in $\C$ to $\bx$.
This function is well-defined (i.e.~there is a unique minimizer) whenever $\C$ satisfies \assref{ass_convex_cons}~\cite[Theorem 6.25]{Beck2017}.
Since $\proj_{\C}(\bx)=\bx$ if and only if $\bx\in\C$, the projection operator also gives us a test for feasibility.
There are many simple sets $\C$ for which $\proj_{\C}$ is easy to compute (see~\cite[Table 6.1]{Beck2017} for others), for example:
\begin{itemize}
    \item For bound constraints $\C=\{\bx : \ba \leq \bx \leq \bb\}$ (including unbounded constraints, $a_i=-\infty$ and/or $b_i=\infty$), we have $[\proj_{\C}(\bx)]_i = \min(\max(a_i, x_i), b_i)$ for $i=1,\ldots,n$; and
    \item If $\C=\{\bx : \ba^T \bx \leq b\}$ for $\ba\neq \bm{0}$, then $\proj_{\C}(\bx) = \bx - \frac{\max(\ba^T\bx-b,0)}{\|\ba\|^2}\ba$.
\end{itemize}
If we have many sets, each of which has their own projection operator, the projection onto the intersection of these sets can be computed using Dykstra's algorithm \cite{Dykstra1986}. 

A good criticality measure for \eqref{eq_convex_cons_generic}, to assess how close we are to a solution, is the following\footnote{Other criticality measures can also be used, with another common choice being $\pi(\bx) = \|\proj_{\C}(\bx-\grad f(\bx)) - \bx\|$.}:
\begin{subequations}
\begin{align}
    \pi(\bx) := \minorrev{-\min_{\bd\in\R^n}} &\: \grad f(\bx)^T \bd, \qquad \forall \bx\in\C. \\
    \minorrev{\text{s.t.}} &\: \minorrev{\bx+\bd\in \C}, \\
    &\: \minorrev{\|\bd\| \leq 1}, \label{eq_convex_criticality}
\end{align}
\end{subequations}
In the unconstrained case, $\C=\R^n$, this reduces to  $\pi(\bx) = \|\grad f(\bx)\|$, as we might hope.
The suitability of $\pi(\bx)$ as a criticality measure comes from the following.

\begin{proposition}[Theorem 12.1.6, \cite{Conn2000}]
    If $f$ is twice continuously differentiable, then $\pi$ is continuous in $\bx$, $\pi(\bx) \geq 0$ for all $\bx\in\C$, and $\pi(\bx)=0$ if and only if $\bx$ is a first-order critical point of \eqref{eq_convex_cons_generic}.\footnote{A point $\bx$ is first-order critical for \eqref{eq_convex_cons_generic} if $-\grad f(\bx) \in \mathcal{N}(\bx)$, where $\mathcal{N}(\bx)$ is the normal cone for $\C$ at $\bx$. This is a first-order necessary condition for a constrained local minimizer \cite[Theorem 12.8]{Nocedal2006}.} 
\end{proposition}

Our goal in this section is to develop a \emph{strictly feasible} trust-region method for solving \eqref{eq_convex_cons_generic}; that is, where we can only evaluate $f$ at feasible points.
This may be necessary, for example, if we have a $\sqrt{x}$ term in $f$ with bounds $x\geq 0$.

If at iteration $k$ we have a quadratic model $m_k$ \eqref{eq_tr_model_dfo} for the objective, we naturally get an estimate of the criticality measure $\pi(\bx_k)$ \eqref{eq_convex_criticality}, namely
\begin{subequations}
\begin{align}
    \pi^m_k := \minorrev{-\min_{\bd\in\R^n}} &\: \bg_k^T \bd, . \\
    \minorrev{\text{s.t.}} &\: \minorrev{\bx_k+\bd\in \C}, \\
    &\: \minorrev{\|\bd\| \leq 1}, \label{eq_convex_criticality_model}
\end{align}
\end{subequations}
Again, if $\C=\R^n$ then $\pi^m_k = \|\bg_k\|$.
Assuming strict feasibility, our trust-region subproblem now looks like: \minorrev{set $\bs_k$ to be an approximate minimizer of}
\begin{align}
    \min_{\bs\in\R^n} m_k(\bx_k+\bs) \qquad \text{s.t.} \qquad \bx_k+\bs\in\C, \quad \text{and} \quad \|\bs\| \leq \Delta_k. \label{eq_trs_convex_cons}
\end{align}
The presence of the feasibility condition $\bx_k+\bs\in\C$ means that the theory and algorithms in \secref{sec_tr} for solving the trust-region subproblem no longer apply.
In this setting, the analog of the Cauchy point is to consider the \emph{projected gradient} path, $\bs_k(t) := \proj_{\C}(\bx_k-t\bg_k) - \bx_k$.
A suitable linesearch applied to $\bs_k(t)$ can yield \cite[Theorem 12.2.2]{Conn2000} a step $\bs_k = \bs_k(t^*)$ satisfying the following assumption, the constrained version of \assref{ass_cauchy_decrease}.
Importantly, instead of the sufficient decrease condition depending on $\|\bg_k\|$, it now depends on $\pi^m_k$. 

\begin{assumption} \label{ass_cauchy_decrease_convex_cons}
    For all $k$, the computed step $\bs_k$ \eqref{eq_trs_convex_cons} satisfies $\bx_k+\bs_k\in\C$,  $\|\bs_k\| \leq \Delta_k$ and
    \begin{align}
        m_k(\bx_k) - m_k(\bx_k+\bs_k) \geq \kappa_s \pi^m_k \min\left(\Delta_k, \frac{\pi^m_k}{\|\bH_k\|+1}, 1\right),
    \end{align}
    for some $\kappa_s\in(0,\frac{1}{2})$.
\end{assumption}

\assref{ass_cauchy_decrease_convex_cons} is almost identical to \assref{ass_cauchy_decrease}, with the unconstrained criticality measure $\|\bg_k\|$  replaced with $\pi^m_k$.
Hence, as we shall see, this new assumption will not have a significant impact on the convergence of our algorithm.

The biggest change comes from the model accuracy requirements.
Our fully linear definition (\defref{def_fully_linear}) is based on comparing the model with the objective at all points inside the trust region.
However, our strict feasibility requirement means that our interpolation model can only be constructed using points in $\C$.
This limits our ability to have accurate models outside $\C$, which is required for \defref{def_fully_linear} to be satisfied. 

\begin{example} \label{ex_convex_poised_demo}
Suppose we have problem dimension $n=2$ and our feasible region is defined by the simple bound constraints $\C = \{\bx\in\R^2 : |x_2| \leq \delta\}$ for some small $\delta\in(0,1)$.
Now suppose we have the base point $\bx=\bm{0}$ and we wish to perform linear interpolation to build a fully linear model in $B(\bm{0},1)$ (i.e.~$\Delta=1$).
Following the approach in \secref{sec_linear_interp}, a natural choice of \emph{feasible} interpolation points is $\{\bm{0}, \be_1, \delta\, \be_2\}$.
\revision{The resulting interpolation model can be shown to be fully linear with $\kappamf,\kappamg=\bigO(1/\delta)$},
and so our fully linear constants can be made arbitrarily large just by changing the feasible region.
So, if we continued to use \defref{def_fully_linear} as our measure of model accuracy, the worst-case complexity of our algorithm could grow like $\bigO(\delta^{-2} \epsilon^{-2})$ for accuracy level $\epsilon$. 
\end{example}

However this definition is stronger than we need.
Requiring strictly feasible interpolation points only limits our ability to approximate $f$ outside the feasible region, but that is exactly the region which is not relevant to our optimization.
So, we now introduce a generalized definition of fully linear interpolation models, capturing the notion that we only care about model accuracy within the feasible region.

\begin{definition} \label{def_fully_linear_convex}
    Suppose we have $\bx\in\C$ and $\Delta>0$.
    A local model $m:\R^n\to\R$ approximating $f:\R^n\to\R$ is $\C$-feasible fully linear in $B(\bx,\Delta)$ if there exist constants $\kappamf,\kappamg>0$, independent of $m$, $\bx$ and $\Delta$, such that
    \begin{subequations} \label{eq_fully_linear_convex}
    \begin{align}
        \max_{\substack{\by \in \C \\ \by\in B(\bx,\Delta)}} |m(\by) - f(\by)| &\leq \kappamf \Delta^2, \label{eq_fully_linear_convex_f} \\
        \max_{\substack{\bx+\bd \in \C \\ \|\bd\| \leq 1}} |(\grad m(\bx) - \grad f(\bx))^T \bd| &\leq \kappamg \Delta. \label{eq_fully_linear_convex_g}
    \end{align}
    \end{subequations}
    Again, we will sometimes use $\kappam := \max(\kappamf,\kappamg)$ for notational convenience.
\end{definition}

In particular, we note the distinction $\|\by-\bx\|\leq\Delta$ in \eqref{eq_fully_linear_convex_f} but $\|\bd\| \leq 1$ in \eqref{eq_fully_linear_convex_g} (as in \eqref{eq_convex_criticality_model}).
Also, compared to \defref{def_fully_linear}, the gradient accuracy condition \eqref{eq_fully_linear_convex_g} only considers the accuracy of $\grad m(\bx) \approx \grad f(\bx)$, not at any other points in the trust region.\footnote{Indeed, our convergence analysis in \secref{sec_dfotr} only requires \eqref{eq_fully_linear_g} to hold at $\by=\bx$.}

We are now ready to state our main \renaming{MBDFO} algorithm for solving \eqref{eq_convex_cons_generic}, given in \algref{alg_basic_tr_dfo_convex}.
It is essentially the same as \algref{alg_basic_tr_dfo}, but replacing $\|\bg_k\|$ with $\pi^m_k$, and using our new sufficient decrease and fully linear conditions.
Specifically, our new model assumptions are given below.

\begin{assumption} \label{ass_model_dfo_convex}
    At each iteration $k$ of \algref{alg_basic_tr_dfo_convex}, the model $m_k$ \eqref{eq_tr_model_dfo} satisfies:
    \begin{enumerate}[label=(\alph*)]
        \item $m_k$ is $\C$-feasible fully linear in $B(\bx_k,\Delta_k)$ with constants $\kappamf,\kappamg>0$ independent of $k$; \label{ass_model_dfo_convex_g}
        \item $\|\bH_k\| \leq \kappa_H - 1$ for some $\kappa_H\geq 1$ (independent of $k$). \label{ass_model_dfo_convex_H}
    \end{enumerate}
\end{assumption}

\begin{algorithm}[tb]
\begin{algorithmic}[1]
\Require Starting point $\bx_0\in\C$ and trust-region radius $\Delta_0>0$. Algorithm parameters: scaling factors $0 < \gammadec < 1 < \gammainc$, acceptance thresholds $0 < \eta_U \leq \eta_S < 1$, and criticality threshold $\mu_c > 0$.
\For{$k=0,1,2,\ldots$}
    \State Build a local quadratic model $m_k$ \eqref{eq_tr_model_dfo} satisfying \assref{ass_model_dfo_convex}.
    \State Solve the trust-region subproblem \eqref{eq_trs_convex_cons} to get a step $\bs_k$ satisfying \assref{ass_cauchy_decrease_convex_cons}.
    \State Evaluate $f(\bx_k+\bs_k)$ and calculate the ratio $\rho_k$ \eqref{eq_ratio_generic}.
    \If{$\rho_k \geq \eta_S$ and $\pi^m_k \geq \mu_c \Delta_k$} 
        \State \textit{(Very successful iteration)} Set $\bx_{k+1}=\bx_k+\bs_k$ and $\Delta_{k+1}=\gammainc\Delta_k$. 
    \ElsIf{$\eta_U \leq \rho_k < \eta_S$ and $\pi^m_k \geq \mu_c \Delta_k$}
        \State \textit{(Successful iteration)} Set $\bx_{k+1}=\bx_k+\bs_k$ and $\Delta_{k+1}=\Delta_k$.
    \Else
        \State \textit{(Unsuccessful iteration)} Set $\bx_{k+1}=\bx_k$ and $\Delta_{k+1}=\gammadec\Delta_k$. 
    \EndIf 
\EndFor
\end{algorithmic}
\caption{Constrained \renaming{MBDFO} trust-region method for solving \eqref{eq_convex_cons_generic}.}
\label{alg_basic_tr_dfo_convex}
\end{algorithm}

We begin by stating the obvious result that all iterates of \algref{alg_basic_tr_dfo_convex} are feasible.

\begin{lemma}
    Suppose \assref{ass_cauchy_decrease_convex_cons} holds and we run \algref{alg_basic_tr_dfo_convex}.
    Then $\bx_k\in\C$ for all $k$.
\end{lemma}
\begin{proof}
    This holds by induction on $k$ because $\bx_0\in\C$ by definition, and $\bx_{k+1}\in\{\bx_k,\bx_k+\bs_k\}$, with $\bx_k+\bs_k\in\C$ whenever $\bx_k\in\C$ from \assref{ass_cauchy_decrease_convex_cons}.
\end{proof}

The $\C$-feasible fully linear requirement \eqref{eq_fully_linear_convex_g} is required principally to ensure the following, which essentially replaces \eqref{eq_fully_linear_g} in measuring the error in the criticality measure.

\begin{lemma}[\revision{Lemma 2.4, \cite{Hough2022}}] \label{lem_convex_crit_measure}
    Suppose Assumptions~\minorrev{\ref{ass_smoothness_1}\ref{ass_smoothness_1_smooth}} and \ref{ass_model_dfo_convex} hold.
    Then $|\pi(\bx_k) - \pi^m_k| \leq \kappamg \Delta_k$ for all iterations $k$ \minorrev{of \algref{alg_basic_tr_dfo_convex}}.
\end{lemma}


Now, the worst-case complexity of \algref{alg_basic_tr_dfo_convex} can be proven using essentially identical arguments to the complexity of \algref{alg_basic_tr_dfo}, replacing $\|\grad f(\bx_k)\|$ and  $\|\bg_k\|$ with $\pi(\bx_k)$ and $\pi^m_k$ respectively.

\begin{theorem} \label{thm_wcc_dfo_convex}
    Suppose Assumptions~\ref{ass_smoothness_1}, \minorrev{\ref{ass_cauchy_decrease_convex_cons} and \ref{ass_model_dfo_convex}} hold and we run \algref{alg_basic_tr_dfo_convex}.
    \revision{If $k_{\epsilon}$ is the first iteration of \algref{alg_basic_tr_dfo_convex} such that $\pi(\bx_k) < \epsilon$, then $k_{\epsilon} = \bigO(\kappa_H (\kappam + \kappa_H)^2 \epsilon^{-2})$.
    Hence $\liminf_{k\to\infty} \pi(\bx_k) = 0$.}
\end{theorem}
\begin{proof}
    \revision{See \cite[Theorem 3.14]{Hough2022} for details.}
\end{proof}



We now have a worst-case complexity bound for \algref{alg_basic_tr_dfo_convex} which matches the unconstrained version, \corref{cor_wcc_dfo} for \algref{alg_basic_tr_dfo}.
What remains is to specify how to construct $\C$-feasible fully linear models using only feasible points.

\subsubsection{Constructing Feasible Interpolation Models}

We wish to construct interpolation models that are $\C$-feasible fully linear (\defref{def_fully_linear_convex}) using only feasible interpolation points.
To do this, we will use the notions of Lebesgue measure and $\Lambda$-poisedness (\defref{def_lebesgue_lambda_poised}), since these more naturally extend to the constrained case than the matrix norms in \secref{sec_model_construction}.

Since we care about the accuracy of our model \minorrev{only} at feasible points, the natural generalization of \defref{def_lebesgue_lambda_poised} to the constrained case is to measure the size of Lagrange polynomials \minorrev{only} within the feasible region.

\begin{definition} \label{def_lebesgue_lambda_poised_convex}
    Suppose we an interpolation set $\mathcal{Y} := \{\by_1,\ldots,\by_p\}\subset\C$ such that its Lagrange polynomials exist. 
    Given $\bx\in\C$ and $\Delta>0$, we say that $\mathcal{Y}$ is $\minorrev{\Lambda}$-poised in $B(\bx,\Delta)\cap \C$ for some $\minorrev{\Lambda}>0$ if
    \begin{align}
        \max_{\by\in B(\bx,\min(\Delta,1))\cap\C} \|\blambda(\by)\|_{\infty} \leq \minorrev{\Lambda}.
    \end{align}
\end{definition}

Aside from only considering $\by\in \C$, the main difference in \defref{def_lebesgue_lambda_poised_convex} compared to \defref{def_lebesgue_lambda_poised} is that we also only consider $\by\in B(\bx,\min(\Delta,1))$ rather than $B(\bx,\Delta)$.
This is required to handle the differing constraints $\|\by-\bx\|\leq \Delta$ and $\|\bd\|\leq 1$ in \eqref{eq_fully_linear_convex}.

\paragraph{Linear Interpolation Models}
First, suppose we construct linear interpolation models using the interpolation set $\mathcal{Y} = \{\by_1,\ldots,\by_p\}\in \C$, for $p=n+1$ by solving \eqref{eq_linear_interp_system} or \eqref{eq_linear_interp_system_scaled}, as usual.
We first note that we have an alternative version of \eqref{eq_linear_interp_lagrange_basic}, namely
\begin{align}
    |m(\by) - f(\bx) - \grad f(\bx)^T (\by-\bx)| \leq \frac{\Lgrad}{2} \beta^2 \|\blambda(\by)\|_1 \min(\Delta,1)^2, \label{eq_linear_interp_lagrange_convex}
\end{align}
for all $\by\in B(\bx,\min(\Delta,1))\cap \C$, where the proof is identical to the unconstrained case, just replacing $B(\bx,\Delta)$ with $B(\bx,\min(\Delta,1))$.
\revision{If $\mathcal{Y}$ is $\Lambda$-poised, then we have $\|\blambda(\by)\|_1 \leq p\Lambda$. By treating the cases $\Delta\leq 1$ and $\Delta>1$ separately, we get a constrained version of \thmref{thm_fully_linear_lagrange}.}


\begin{theorem} \label{thm_fully_linear_lagrange_convex}
    Suppose $f$ satisfies \minorrev{\assref{ass_smoothness_1}\ref{ass_smoothness_1_smooth}} and we construct a linear model \eqref{eq_linear_model_generic} for $f$ by solving \eqref{eq_linear_interp_system_scaled}, where we assume $\hat{\bM}$ is invertible.
    If $\bx\in\C$, with $\by_i\in\C$ and $\|\by_i-\bx\| \leq \beta \Delta$ for some $\beta>0$ and all $i=1,\ldots,p$ (where $p=n+1$ is the number of interpolation points), and the interpolation set is $\minorrev{\Lambda}$-poised in $B(\bx,\Delta)\cap \C$, then the model is $\C$-feasible fully linear in $B(\bx,\Delta)$ with constants
    \begin{align}
        \kappamf = \Lgrad \beta^2 \minorrev{p \Lambda} + \frac{\Lgrad}{2}, \quad \text{and} \quad \kappamg = \Lgrad \beta^2 \minorrev{p \Lambda}. 
    \end{align}
\end{theorem}
\begin{proof}
    \revision{See \cite[Theorem 4.4]{Hough2022}.}
\end{proof}

The fully linear constants in \thmref{thm_fully_linear_lagrange_convex} are of size $\kappamf,\kappamg=\bigO(\Lgrad \minorrev{p \Lambda})$, which matches the values in the unconstrained case $\C=\R^n$ from \thmref{thm_fully_linear_lagrange}.

\begin{example}[\exref{ex_convex_poised_demo} revisited]
Consider again the case $n=2$ with $\C=\{\bx\in\R^2 : |x_2| \leq \delta\}$, with linear interpolation in $B(\bm{0},1)$ using $\{\bm{0}, \be_1, \delta \: \be_2\}$.
By considering the enlarged trust region $[-1,1]^2$ instead of $B(\bm{0},1)$, we see that the interpolation set is $\minorrev{\Lambda}$-poised in $B(\bm{0},1)\cap \C$ with $\minorrev{\Lambda} \leq 3$ and hence $\kappamf$ and $\kappamg$ are independent of $\delta$ from \thmref{thm_fully_linear_lagrange_convex}, a significant improvement over $\kappamf,\kappamg=\bigO(1/\delta)$ by applying the unconstrained approach from \defref{def_fully_linear} and \thmref{thm_fully_linear}.
\end{example}

\paragraph{Minimum Frobenius Norm Models}
A similar approach works if we wish to use minimum Frobenius norm quadratic interpolation \eqref{eq_min_frob_problem}, with the same definition of Lebesgue measure and $\Lambda$-poisedness \defref{def_lebesgue_lambda_poised_convex}.

To begin, we state a bound on the size of the model Hessian, similar to \eqref{eq_min_frob_hess_bound_poised} in the unconstrained case.
Unfortunately we cannot bound $\|\bH\|$ entirely; we instead bound certain Rayleigh quotient-type expressions.

\begin{lemma} \label{lem_min_frob_hess_bound_convex}
    Suppose $f$ satisfies \minorrev{\assref{ass_smoothness_1}\ref{ass_smoothness_1_smooth}} and we construct a quadratic model $m$ \eqref{eq_quadratic_model_generic} for $f$ by solving \eqref{eq_min_frob_interp_system_scaled}, where we assume $\hat{\bF}$ is invertible.
    If $\bx\in\C$, with $\by_i\in\C$ and $\|\by_i-\bx\| \leq \beta \min(\Delta,1)$ for some $\beta>0$ and all $i=1,\ldots,p$, and $\{\by_1,\ldots,\by_p\}$ is $\minorrev{\Lambda}$-poised in $B(\bx,\Delta)\cap \C$ then the model Hessian $\bH$ satisfies
    \begin{align}
	\max_{i,j=1,\ldots,p} \frac{|(\by_i-\bx)^T \bH (\by_j-\bx)|}{\beta^2 \min(\Delta,1)^2} \leq \kappa^{\C}_H := \Lgrad p (8\minorrev{\Lambda}\beta^2 + 36\minorrev{\Lambda}\beta + 58\minorrev{\Lambda} + 6), \label{eq_min_frob_hess_bound_convex1}
    \end{align}
    and
    \begin{align}
        |(\by-\bx)^T \bH (\by-\bx)| \leq \kappa^{\C}_H \beta^2 \minorrev{p \Lambda^2} \min(\Delta,1)^2, \label{eq_min_frob_hess_bound_convex2}
    \end{align}
    for all $\by \in B(\bx,\min(\Delta,1))\cap \C$.
\end{lemma}
\begin{proof}
    The proof of \eqref{eq_min_frob_hess_bound_convex1} uses similar ideas to the proof of \eqref{eq_min_frob_hess_bound_poised} but requires a more lengthy calculation which we omit.
    See \cite[Lemma 4.4]{Roberts2024} for details.
    To show \eqref{eq_min_frob_hess_bound_convex2}, fix $\by \in B(\bx,\min(\Delta,1))\cap \C$. 
    We recall from \eqref{eq_min_frob_fully_linear_lagrange_tmp1} that we may write $\by-\bx = \sum_{i=1}^{p} \ell_i(\by) (\by_i-\bx)$.
    So, from \eqref{eq_min_frob_hess_bound_convex1} we get
    \begin{align}
        |(\by-\bx)^T \bH (\by-\bx)| &\leq \sum_{i,j=1}^{p} |\ell_i(\by)| \: |\ell_j(\by)| \: |(\by_i-\bx)^T \bH (\by_j-\bx)|, \\
        &\leq \kappa^{\C}_H \beta^2 \min(\Delta,1)^2 \|\blambda(\by)\|_1 \max_{i=1,\ldots,p} |\ell_i(\by)|,
    \end{align}
    and are done, after noting $\|\blambda(\by)\|_1 \leq \minorrev{p \Lambda}$ and $\max_{i=1,\ldots,p} |\ell_i(\by)| \leq \minorrev{\Lambda}$.
\end{proof}

\begin{remark}
    If we have $\|\bH\| \leq \kappa_H-1$ (as in \assref{ass_model_dfo_convex}\ref{ass_model_dfo_convex_H}), then we can replace \eqref{eq_min_frob_hess_bound_convex2} with $|(\by-\bx)^T \bH (\by-\bx)| \leq \kappa_H \min(\Delta,1)^2$.
\end{remark}

The constrained version of \thmref{thm_min_frob_fully_linear_lagrange} is the following.

\begin{theorem} \label{thm_min_frob_fully_linear_lagrange_convex}
    Suppose $f$ satisfies \minorrev{\assref{ass_smoothness_1}\ref{ass_smoothness_1_smooth}} and we construct a quadratic model $m$ \eqref{eq_quadratic_model_generic} for $f$ by solving \eqref{eq_min_frob_interp_system_scaled}, where we assume $\hat{\bF}$ is invertible.
    If $\bx\in\C$, with $\by_i\in\C$ and $\|\by_i-\bx\| \leq \beta \min(\Delta,1)$ for some $\beta>0$, and the interpolation set is $\minorrev{\Lambda}$-poised in $B(\bx,\Delta)\cap \C$, then the model is $\C$-feasible fully linear in $B(\bx,\Delta)$ with constants
    \begin{align}
        \kappamf = \frac{3(\Lgrad + \kappa^{\C}_H)}{2} \beta^2 \minorrev{p \Lambda} + \frac{1}{2}\kappa^{\C}_H \beta^2 \minorrev{p \Lambda^2}  + \frac{\Lgrad}{2}, \qquad \text{and} \qquad \kappamg = (\Lgrad + \kappa^{\C}_H) \beta^2 \minorrev{p \Lambda},
    \end{align}
    with $\kappa^{\C}_H$ defined in \lemref{lem_min_frob_hess_bound_convex}.
\end{theorem}
\begin{proof}
    First, we follow the proof of \thmref{thm_min_frob_fully_linear_lagrange}, except using $|(\by_i-\bx)^T \bH (\by_j-\bx)| \leq \kappa^{\C}_H \beta^2 \min(\Delta,1)^2$ from \eqref{eq_min_frob_hess_bound_convex1} in the derivation of \eqref{eq_min_frob_fully_linear_lagrange_tmp2} to get
    \begin{align}
        |c + \bg^T(\by-\bx) - f(\bx) - \grad f(\bx)^T (\by-\bx)| \leq \frac{\Lgrad + \kappa^{\C}_H}{2} \beta^2 \|\blambda(\by)\|_1 \min(\Delta,1)^2, \label{eq_min_frob_fully_linear_convex_tmp1}
    \end{align}
    for all $\by\in B(\bx,\min(\Delta,1))\cap \C$.
    Otherwise, the proof is largely similar to that of \thmref{thm_fully_linear_lagrange_convex}; see \cite[Theorem 4.7]{Roberts2024} for more details.
\end{proof} 

\paragraph{Model Improvement}
As in \secref{sec_model_improvement}, we can use our new notion of $\Lambda$-poisedness (\defref{def_lebesgue_lambda_poised_convex}) to build algorithms to improve the geometry of an interpolation set.
Just like the unconstrained case, this works equally well in the linear and minimum Frobenius norm quadratic interpolation cases.
Indeed, our definitions ensure that \algref{alg_geom_improvement} extends trivially to the constrained case, as shown in \algref{alg_geom_improvement_convex}.
The only changes are that the initial interpolation set must all be feasible, and our new point $\by$ must also be feasible.

\begin{algorithm}[tb]
\begin{algorithmic}[1]
\Require Interpolation set $\mathcal{Y}$ which is $\minorrev{\Lambda}$-poised in $B(\bx,\Delta)\cap\C$, $\bx\in\C$, with $\by_i\in\C$ and $\|\by_i-\bx\| \leq \beta \min(\Delta,1)$ for some $\beta\geq 1$, desired poisedness constant $\minorrev{\Lambda}^* > 1$.
\While{$\minorrev{\Lambda} > \minorrev{\Lambda}^*$}
    \State Find $\by\in B(\bx,\min(\Delta,1))\cap \C$ and $i\in\{1,\ldots,p\}$ such that $|\ell_i(\by)| > \minorrev{\Lambda}^*$.
    \State Update $\mathcal{Y}$ by replacing $\by_i$ with $\by$. 
    \State Recompute the poisedness constant $\minorrev{\Lambda}$ of the new $\mathcal{Y}$.
\EndWhile
\State \Return $\mathcal{Y}$.
\end{algorithmic}
\caption{Interpolation set improvement (convex constrained case).}
\label{alg_geom_improvement_convex}
\end{algorithm}

\begin{theorem} \label{thm_geom_improvement_convex}
    \algref{alg_geom_improvement_convex} terminates in finite time, and the resulting interpolation set $\mathcal{Y}$ is $\minorrev{\Lambda}^*$-poised in $B(\bx,\Delta)$, and we have $\by_i\in\C$ and $\|\by_i-\bx\| \leq \beta \min(\Delta,1)$ for all $\by_i\in \mathcal{Y}$.
\end{theorem}
\begin{proof}
    The proof is identical to that of \thmref{thm_geom_improvement}, since Lemmas~\ref{lem_det_updating_lin_quad} and \ref{lem_det_updating_min_frob} still hold, and replacing $\Delta$ with $\min(\Delta,1)$.
\end{proof}

However, the assumption in \algref{alg_geom_improvement_convex} that the initial interpolation set is $\minorrev{\Lambda}$-poised in $B(\bx,\Delta)\cap \C$ for some $\minorrev{\Lambda}$ (with only feasible points) is not so easy to ensure.
In the unconstrained case, we saw example sets in \secref{sec_model_construction} that ensure the relevant interpolation linear system is invertible, but these do not necessarily work here as we cannot guarantee they are feasible.
Instead, we can use \algref{alg_interp_init_convex} to construct a strictly feasible set with invertible interpolation linear system.

\begin{algorithm}[tb]
\begin{algorithmic}[1]
\Require Interpolation region $B(\bx,\min(\Delta,1))$.
\State Construct an initial interpolation set $\mathcal{Y} \subset B(\bx,\min(\Delta,1))$ which has an invertible interpolation linear system (e.g.~as in \secref{sec_model_construction}).
\While{there exists $\by_i \in \mathcal{Y}$ with $\by_i \notin \C$}
    \State Find $\by\in B(\bx,\min(\Delta,1))\cap \C$ such that $\ell_i(\by) \neq 0$.
    \State Update $\mathcal{Y}$ by replacing $\by_i$ with $\by$.
\EndWhile
\State \Return $\mathcal{Y}$
\end{algorithmic}
\caption{Interpolation set initialization (convex constrained case).}
\label{alg_interp_init_convex}
\end{algorithm}

\begin{theorem} 
    \algref{alg_interp_init_convex} terminates in finite time, and the resulting interpolation set $\mathcal{Y}$ is contained in $B(\bx,\min(\Delta,1))\cap \C$ and has an invertible interpolation linear system.
\end{theorem}
\begin{proof}
    All initial points are in $B(\bx,\min(\Delta,1))$, and any initial point $\by_i \notin \C$ is replaced by a new point $\by\in B(\bx,\min(\Delta,1))\cap \C$, and so \algref{alg_interp_init_convex} terminates after at most $p$ iterations.
    The interpolation linear system remains invertible after each iteration because it starts invertible, and we apply Lemma~\ref{lem_det_updating_lin_quad} or \ref{lem_det_updating_min_frob} together with $\ell_i(\by) \neq 0$.
\end{proof}

\subsection{General Constraints}

Since \renaming{MBDFO} differs most from derivative-based trust-region methods in its model construction, our focus in this section has been on building accurate models using strictly feasible points, and hence we have restricted ourselves to problems with simple constraints \eqref{eq_convex_cons_generic}, rather than problems with general nonlinear constraints (for which we also may only have zeroth order oracles).
To conclude this section, we briefly outline one approach for solving general constrained problems of the form \eqref{eq_generic_cons}.
This form is often augmented with explicit bound constraints, although we do not do this here for simplicity. 
We assume that the objective $f$ and all constraint functions $c_i$ for $i\in\mathcal{E}\cup\mathcal{I}$ only have zeroth order oracles, and so their derivatives must be approximated.

The algorithm we outline will be a sequential quadratic programming (SQP) \renaming{MBDFO} method, based on the COBYQA (Constrained Optimization BY Quadratic Approximation) algorithm available in recent versions of Python's widely used SciPy optimization library. 
At each iteration of an SQP method, we construct a local quadratic model for the objective $f$, 
\begin{align}
    f(\by) \approx m_k(\by) := c_k + \bg_k^T (\by-\bx_k) + \frac{1}{2}(\by-\bx_k)^T \bH_k (\by-\bx_k), \label{eq_cobyqa_f_model}
\end{align}
and local linear models for all constraints $c_i$,
\begin{align}
    c_i(\by) \approx m_{k,i}(\by) := c_{k,i} + \bg_{k,i}^T (\by-\bx_k), \qquad \forall i\in\mathcal{E}\cup\mathcal{I}. \label{eq_cobyqa_c_model}
\end{align}
The main difference in these approximations is that we want to choose the model Hessian in $m_k$ to achieve $\bH_k \approx \grad^2_{\bx} L(\bx_k,\blambda_k)$ (recalling the Lagrangian \eqref{eq_lagrangian}) rather than $\bH_k\approx \grad^2 f(\bx_k)$, and so we actually need to approximate Hessians for each $c_i$ and have suitable Lagrange multiplier estimates.
With these approximations, the trust-region subproblem becomes
\begin{subequations} \label{eq_cobyqa_trs}
\begin{align}
    \min_{\bs\in\R^n} &\: m_k(\bs), \\
    \text{s.t.} &\: m_{k,i}(\bx_k+\bs) = 0, \qquad \forall i\in\mathcal{E}, \\
    &\: m_{k,i}(\bx_k+\bs) \leq 0, \qquad \forall i\in\mathcal{I}, \\
    &\: \|\bs\| \leq \Delta_k.
\end{align}
\end{subequations}
In the unconstrained case, we decided if a step was good by measuring the decrease in the objective \eqref{eq_ratio_generic}.
Here, we measure the quality of a step $\bs_k$ by using a \emph{merit function}, which combines the objective value and size of any constraint violations into a single scalar.
For example, the $\ell_2$ merit function used in COBYQA is defined as
\begin{align}
    \phi(\bx, \gamma) := f(\bx) + \gamma \Phi(\bx), \qquad \text{where} \quad \Phi(\bx) := \sqrt{\sum_{i\in\mathcal{E}} c_i(\bx)^2 + \sum_{i\in\mathcal{I}} \max(c_i(\bx),0)^2},
\end{align}
where the penalty parameter $\gamma>0$ controls the relative weight applied to reductions in $f$ compared to improvements in feasibility.\footnote{An alternative but common choice of merit function uses the $\ell_1$ norm of the constraint violation.}
Note that $\phi(\bx,\gamma) = f(\bx)$ if $\bx$ is feasible.
Given our models \eqref{eq_cobyqa_f_model} and \eqref{eq_cobyqa_c_model}, we can derive our approximate merit function
\begin{align}
    \phi_k(\bx,\gamma) := m_k(\bx) + \gamma \Phi_k(\bx), \qquad \text{where} \quad \Phi_k(\bx) := \sqrt{\sum_{i\in\mathcal{E}} m_{k,i}(\bx)^2 + \sum_{i\in\mathcal{I}} \max(m_{k,i}(\bx), 0)^2}.
\end{align}
The choice of accepting/rejecting a step and updating the trust-region radius is similar to the unconstrained case.

A prototypical algorithm is given in \algref{alg_cobyqa}.
We note that increase the merit penalty parameter $\gamma_k$ as $k\to\infty$ ensures we encourage $\bx_k$ to gradually become feasible.
The condition $\gamma_k \geq \|\blambda_k\|$, where $\blambda_k$ is a Lagrange multiplier estimate is to ensure that the merit function is \emph{exact}; that is, the true solution $\bx^*$ is also a minimizer of $\phi(\cdot,\gamma)$ for all $\gamma$ sufficiently large.
In the case of the $\ell_2$ merit function, `sufficiently large' is $\gamma \geq \|\blambda^*\|$ \cite[Theorem 4.2.1]{Ragonneau2022}.
Motivated by \eqref{eq_kkt}, the Lagrange multiplier estimates $\blambda_k$ are computed by solving the (bound-constrained) linear least-squares problem
\begin{subequations} \label{eq_cobyqa_multipler_est}
\begin{align}
    \blambda_k \minorrev{\in} \argmin_{\blambda\in\R^{|\mathcal{E}|+|\mathcal{I}|}} &\: \left\|\bg_k + \sum_{i\in\mathcal{E}\cup\mathcal{I}} \lambda_i \bg_{k,i}\right\|^2, \\
    \text{s.t.} &\: \lambda_i = 0, \qquad \forall i\in\{j\in\mathcal{I} : c_{k,i} < 0\}, \\
    &\: \lambda_i \geq 0, \qquad \forall i\in\{j\in \mathcal{I} : c_{k,i} \geq 0\}.
\end{align}
\end{subequations}

\begin{algorithm}[tb]
\begin{algorithmic}[1]
\Require Starting point $\bx_0\in\R^n$, initial Lagrange multiplier estimates $\blambda_0\in\R^{|\mathcal{E}|+|\mathcal{I}|}$ and trust-region radius $\Delta_0>0$. Algorithm parameters: scaling factors $0 < \gammadec < 1 < \gammainc$, and acceptance thresholds $0 < \eta_U \leq \eta_S < 1$.
\For{$k=0,1,2,\ldots$}
    \State Build \minorrev{quadratic} models $m_k$ \eqref{eq_cobyqa_f_model} for $f$ and $m_{k,i}$ \eqref{eq_cobyqa_c_model} for each $c_i$.
    \State Approximately solve the subproblems \eqref{eq_cobyqa_trs_normal} and \eqref{eq_cobyqa_trs_tangential} to get a step $\bs_k = \bn_k + \bt_k$.
    \State Evaluate $f$ and all constraints $c_i$ at $\bx_k+\bs_k$ and calculate the ratio 
    \begin{align}
        \rho_k = \frac{\text{actual reduction in merit function}}{\text{predicted reduction in merit function}} = \frac{\phi(\bx_k,\gamma_k) - \phi(\bx_k+\bs_k, \gamma_k)}{\phi_k(\bx_k,\gamma_k) - \phi_k(\bx_k+\bs_k,\gamma_k)},
    \end{align}
    for some value $\gamma_k \geq \max(\gamma_{k-1}, \|\blambda_k\|)$, chosen such that $\phi_k(\bx_k+\bs_k,\gamma_k) < \phi_k(\bx_k,\gamma_k)$.
    \If{$\rho_k \geq \eta_S$} 
        \State \textit{(Very successful iteration)} Set $\bx_{k+1}=\bx_k+\bs_k$ and $\Delta_{k+1}=\gammainc\Delta_k$. 
    \ElsIf{$\eta_U \leq \rho_k < \eta_S$}
        \State \textit{(Successful iteration)} Set $\bx_{k+1}=\bx_k+\bs_k$ and $\Delta_{k+1}=\Delta_k$.
    \Else
        \State \textit{(Unsuccessful iteration)} Set $\bx_{k+1}=\bx_k$ and $\Delta_{k+1}=\gammadec\Delta_k$. 
    \EndIf 
    \State Estimate the Lagrange multipliers $\blambda_{k+1}$ associated with $\bx_{k+1}$ by solving \eqref{eq_cobyqa_multipler_est}.
\EndFor
\end{algorithmic}
\caption{Example \renaming{MBDFO} trust-region SQP method for solving \eqref{eq_generic_cons}.}
\label{alg_cobyqa}
\end{algorithm}

Aside from the management of an interpolation set, which can use the techniques described in previous sections (but the interpolation problem is solved for $f$ and all constraints $c_i$), the main practical consideration is the calculation of the step \eqref{eq_cobyqa_trs}.
The largest difficulty is that the linearized constraints may be infeasible.
A common solution to this issue is to decompose the computed step into a normal and tangential component, $\bs_k = \bn_k + \bt_k$.
The normal step $\bn_k$ aims to reduce constraint violation and the tangential step aims to decrease the objective without worsening the constraint violation.
There are several ways to do this, but COBYQA uses the following approach.
First, calculate the normal step \minorrev{$\bn_k$ by (approximately) solving}
\begin{subequations} \label{eq_cobyqa_trs_normal}
\begin{align}
    \min_{\bn\in\R^n} &\: \Phi_k(\bx_k+\bn)^2, \\
    \text{s.t.} &\: \|\bn\| \leq \theta \Delta_k,
\end{align}
\end{subequations}
for some scalar $\theta\in(0,1)$ (e.g.~$\theta=0.8$), and then calculate the tangential step \minorrev{$\bt_k$ by approximately solving}
\begin{subequations} \label{eq_cobyqa_trs_tangential}
\begin{align}
    \min_{\bt\in\R^n} &\: m_k(\bx_k+\bn_k + \bt), \\
    \text{s.t.} &\: \bg_{k,i}^T \bt = 0, \qquad \forall i\in\mathcal{E}, \label{eq_cobyqa_trs_tangential_cons1} \\
    &\: \bg_{k,i}^T \bt \leq \max(-m_{k,i}(\bx_k+\bn_k), 0), \qquad \forall i\in\mathcal{I}, \label{eq_cobyqa_trs_tangential_cons2} \\
    &\: \|\bn_k + \bt\| \leq \Delta_k.
\end{align}
\end{subequations}
Both of these subproblems \eqref{eq_cobyqa_trs_normal} and \eqref{eq_cobyqa_trs_tangential} can be reformulated to be trust-region-like subproblems with linear constraints which can be solved using algorithms similar to the conjugate gradient-based method developed in \cite{Powell2015}; see \cite[Chapter 6]{Ragonneau2022} for details.

We \minorrev{also note} that the requirement that $\gamma_k$ is chosen to ensure $\phi_k(\bx_k+\bs_k,\gamma_k) < \phi_k(\bx_k,\gamma_k)$ can be satisfied relatively easily, under mild assumptions on the step components $\bn_k$ and $\bt_k$.

\begin{lemma}
    \minorrev{If, on iteration $k$ of \algref{alg_cobyqa},} $\bs_k = \bn_k + \bt_k$ where $\bn_k$ \eqref{eq_cobyqa_trs_normal} either satisfies $\Phi_k(\bx_k+\bn_k) < \Phi(\bx_k)$ or $\bn_k = \bm{0}$,
    and $\bt_k$ \eqref{eq_cobyqa_trs_tangential} satisfies $m_k(\bx_k+\bt_k) < m_k(\bx_k)$, then $\phi_k(\bx_k+\bs_k,\gamma) < \phi_k(\bx_k,\gamma)$ for all $\gamma>0$ sufficiently large.
\end{lemma}
\begin{proof}
    The constraints \eqref{eq_cobyqa_trs_tangential_cons1} and \eqref{eq_cobyqa_trs_tangential_cons2} ensure that $\Phi_k(\bx_k+\bs_k) = \Phi_k(\bx_k+\bn_k+\bt_k) \leq \Phi_k(\bx_k+\bn_k)$.
    
    First, if $\Phi_k(\bx_k+\bn_k) < \Phi_k(\bx_k)$, then we have $\Phi_k(\bx_k+\bs_k) < \Phi_k(\bx_k)$, and the result follows by noting that $\phi_k(\bx_k,\gamma) - \phi_k(\bx_k+\bs,\gamma) = m_k(\bx_k) - m_k(\bx_k+\bs) + \gamma (\Phi_k(\bx_k) - \Phi_k(\bx_k+\bs))$ is linear and increasing in $\gamma>0$.
    Instead, if $\bn_k=\bm{0}$ then $\bs_k=\bt_k$ so $\Phi_k(\bx_k+\bs_k) \leq \Phi_k(\bx_k)$ and $m_k(\bx_k+\bs_k) = m_k(\bx_k+\bt_k) < m_k(\bx_k)$, so $\phi_k(\bx_k+\bs_k,\gamma) < \phi_k(\bx_k,\gamma)$ for all $\gamma>0$.
\end{proof}

\revision{\algref{alg_cobyqa} does not have any convergence guarantees as written, and was developed purely as a practical method. Other \renaming{MBDFO} SQP methods such as \cite{Troltzsch2016,Hannanu2024} also do not have convergence guarantees.
However, stochastic SQP methods (extending the ideas from \secref{sec_stochastic_noise} to the constrained setting) with convergence guarantees \cite{Fang2024,Fang2026complexity,Fang2026} can be used in the deterministic setting, and there does not appear to be any fundamental issue preventing the convergence theory for derivative-based SQP trust-region methods \cite[Chapter 15]{Conn2000} being adapted to the \renaming{MBDFO} setting.
}

\subsubsection{Evaluation Failures and Hidden Constraints}

\revision{In some situations, we may also have to concern ourselves with objective evaluations failing (i.e.~our procedure for evaluating $f(\bx)$ fails or crashes unexpectedly). 
For example, the objective evaluation for a helicopter rotor blade design problem in \cite{Booker1999} failed to compute over 60\% of the time. 
More generally, large-scale computations (such as many objective evaluations in DFO contexts) can be exposed to hardware failure; for one high-performance computing system, hardware failures occurred on average once every 7.5 days \cite{Snir2014}.

This can be thought of as an example of \emph{hidden constraints} \cite{Chen2016,LeDigabel2024}, where our problem takes the form
\begin{align}
    \min_{\bx\in\mathcal{D}} \: f(\bx),
\end{align}
for some feasible set $\mathcal{D}\subseteq\R^n$ which is unknown to the algorithm (and where $\bx\in\mathcal{D}$ represents `evaluating $f(\bx)$ succeeds').
This can be rigorously handled with unconstrained algorithms if we let $f:\R^n\to\R\cup\{+\infty\}$ and define $f(\bx) := +\infty$ whenever $\bx\notin\mathcal{D}$.
Here, we must prove convergence to Clarke stationary points, although to the author's knowledge this has not been applied to MBDFO, only direct search \cite{Chen2016}.
In practice, this can instead be handled by setting $f(\bx)$ to some large (finite) value if $\bx\notin\mathcal{D}$ \cite{Ragonneau2024}.
}

\subsubsection*{Notes and References}
{\small 
\secref{sec_convex_constraints} is based primarily on \cite{Hough2022,Roberts2024}.
These works draw on the derivative-based trust-region theory outlined in \cite[Chapter 12]{Conn2000}, the complexity analysis for convex-constrained cubic regularization \cite{Cartis2012b} and \minorrev{convergence} theory for convex-constrained \renaming{MBDFO} using the original notion (\defref{def_fully_linear}) of fully linear models \cite{Conejo2013}.
A good resource for the theory of convex sets and convex optimization is \cite{Beck2017}.

The description of COBYQA (\algref{alg_cobyqa}) for general constraints is taken from \cite{Ragonneau2022}.
More detail about derivative-based SQP methods, including convergence theory, alternative approaches for calculating $\bs_k$ and alternative globalization mechanisms such as filters (instead of a merit function), see \cite[Chapter 15]{Conn2000}.
The specific step calculation given by \eqref{eq_cobyqa_trs_normal} and \eqref{eq_cobyqa_trs_tangential} is a variant of the Byrd--Omojokun method \cite[Chapters 15.4.2 \& 15.4.4]{Conn2000}, originally from \cite{Omojokun1989} (a thesis by Omojokun, supervised by Byrd).
There are many other works which describe \renaming{MBDFO} methods for handling general constraints \eqref{eq_generic_cons} described in the survey \cite[Section 7]{Larson2019}.
However, we note \minorrev{in particular \cite{Troltzsch2016,Hannanu2024}, which propose} SQP methods very similar to COBYQA. 
See \cite{Nocedal2006} for an overview of theory and algorithms for general constrained optimization, including quadratic programming.}


{\small \renaming{MBDFO} for nonlinearly constrained problems is not as well-developed as for unconstrained problems. There are few implementations outside of COBYQA, and to \minorrev{the author's} knowledge none with both good practical performance and theoretical guarantees.
\revision{By contrast, constrained optimization within the Mesh Adaptive Direct Search framework is well-studied in the sense of proving convergence to Clarke stationary points (e.g.~\cite{Audet2004,Audet2006a,Audet2009,Audet2010,Audet2015}).}}

\section{Optimization for Noisy Problems} \label{sec_noisy}

A key use case of \renaming{MBDFO} methods is optimization for problems with some degree of noise present.
In general, this refers to situations where an exact zeroth order oracle is not available (i.e.~we can only evaluate $f$, but not exactly).
This comprises situations, for example, where the calculation of the objective $f$ requires performing a Monte Carlo simulation, or involves a physical/real-world experiment subject to some inherent noise or uncertainty.

In this section we discuss the cases of both \emph{deterministic} and \emph{stochastic} noise.
By `deterministic noise', we mean the following: if we evaluate the noisy function repeatedly at the same point, we get the same answer (e.g.~roundoff errors, finite termination of an iterative process).
Stochastic noise refers to the case where repeated evaluations at the same point can produce different value, and so the noisy value of $f(\bx)$ can be treated as a random variable (e.g.~Monte Carlo simulation, experimental errors).

\revision{In both cases, we outline realistic model accuracy assumptions (in place of `fully linear at every iteration'), and analyze the resulting variants of \algref{alg_basic_tr_dfo}.}
In practice, a simple heuristic \revision{(with theoretical justification for deterministic noise, see \remref{rem_gammadec_noisy})} for handling noisy problems is to use a generic algorithm such as \algref{alg_basic_tr_dfo}, but set $\gammadec$ very close to 1, for example $\gammadec=0.98$, rather than more typical values such as $\gammadec=0.5$ \cite{Cartis2019a}.

\subsection{Deterministic Noise}
In the case of deterministic noise, we essentially have a zeroth order oracle that reliably returns an incorrect value.
We will assume that the noise is uniformly bounded.

\begin{assumption} \label{ass_noisy_oracle}
    Given objective $f:\R^n\to\R$ \eqref{eq_problem}, we only have access to the noisy zeroth order oracle $\tilde{f}:\R^n\to\R$, satisfying $|\tilde{f}(\bx) - f(\bx)| \leq \epsilon_f$ for all $\bx\in\R^n$, for some $\epsilon_f >0$.
\end{assumption}

We illustrate the difficulty that the presence of noise introduces to optimization algorithms by first considering the case of finite differencing.
The below result says that, without due care, finite differencing a noisy zeroth order oracle can give very inaccurate derivative estimates.

\begin{lemma} \label{lem_noisy_fin_diff}
    Suppose $f:\R^n\to\R$ is twice continuously differentiable and $\|\grad^2 f(\by)\| \leq M$ for all $\by\in B(\bx_k,h)$, for some $h>0$.
    If we have a noisy oracle $\tilde{f}$ satisfying \assref{ass_noisy_oracle} and we compute (c.f.~\eqref{eq_fin_diff})
    \begin{align}
        [\bg_k]_i := \frac{\tilde{f}(\bx_k+h\be_i) - \tilde{f}(\bx_k)}{h}, \qquad \forall i=1,\ldots,n, \label{eq_forward_fin_diff_noisy}
    \end{align}
    then $\|\bg_k - \grad f(\bx_k)\|_{\infty} \leq \frac{Mh}{2} + \frac{2\epsilon_f}{h}$.
\end{lemma}
\begin{proof}
    Let $\hat{\bg}_k$ be the noise-free finite difference estimate from \eqref{eq_fin_diff}.
    From standard finite differencing theory (e.g.~\cite[Chapter 8.1]{Nocedal2006}) we have $\|\hat{\bg}_k - \grad f(\bx_k)\|_{\infty} \leq \frac{Mh}{2}$.
    We also have
    \begin{align}
        |[\bg_k - \hat{\bg}_k]_i| \leq \frac{|\tilde{f}(\bx_k+h\be_i) - f(\bx_k+h\be_i)| + |\tilde{f}(\bx_k) - f(\bx_k)|}{h} \leq \frac{2\epsilon_f}{h}, \qquad \forall i=1,\ldots,n,
    \end{align}
    \minorrev{as required}.
\end{proof}

This demonstrates that, in the noisy case, the error in finite differencing can increase significantly as $h\to 0^{+}$.
In fact, to minimize the error bound in \lemref{lem_noisy_fin_diff} we should take $h = 2\sqrt{\epsilon_f/M}$ to get gradient error $\|\bg_k-\grad f(\bx_k)\|_{\infty} \leq 2\sqrt{M\epsilon_f}$.\footnote{For this reason, when estimating the gradient of an objective using \eqref{eq_forward_fin_diff_noisy}, both Python's \texttt{scipy.optimize.minimize} and MATLAB's \texttt{fminunc} use the default value $h = \sqrt{\epsilon_{\text{machine}}}$, where $\epsilon_{\text{machine}}$ is the machine epsilon.}
The advantage of \renaming{MBDFO} is that interpolation models are constructed from well-spaced points, typically of distance $\bigO(\Delta_k)$ from each other, and so we avoid the $h\to 0^{+}$ case until $\Delta_k$ is very small. 

For example, suppose we construct our model \eqref{eq_tr_model_dfo} using linear interpolation to points $\by_1,\ldots,\by_p$ for $p=n+1$ \eqref{eq_linear_interp_system}.
With a noisy zeroth order oracle, we now solve
\begin{align}
    \bM\begin{bmatrix} c \\ \bg \end{bmatrix} = \begin{bmatrix} \tilde{f}(\by_1) \\ \vdots \\ \tilde{f}(\by_p) \end{bmatrix} = \begin{bmatrix} f(\by_1) \\ \vdots \\ f(\by_p) \end{bmatrix} + \begin{bmatrix} \tilde{f}(\by_1)-f(\by_1) \\ \vdots \\ \tilde{f}(\by_p)-f(\by_p) \end{bmatrix}. \label{eq_linear_interp_system_noisy}
\end{align}
This gives us the following version of \thmref{thm_fully_linear}.

\begin{theorem} \label{thm_fully_linear_noisy}
    Suppose $f$ satisfies \minorrev{\assref{ass_smoothness_1}\ref{ass_smoothness_1_smooth}} and $\tilde{f}$ satisfies \assref{ass_noisy_oracle} and we construct a linear model \eqref{eq_linear_model_generic} for $f$ by solving \eqref{eq_linear_interp_system_noisy}, where we assume $\hat{\bM}$ is invertible.
    If $\|\by_i-\bx\| \leq \beta \Delta$ for some $\beta>0$ and all $i=1,\ldots,p$ (where $p=n+1$ is the number of interpolation points), then the model satisfies
    \begin{subequations} \label{eq_thm_fully_linear_noisy_constants}
    \begin{align}
        |m(\by) - f(\by)| &\leq \kappamf \Delta^2 + \tkappamf \: \epsilon_f, \\
        \|\grad m(\by) - \grad f(\by)\| &\leq \kappamg \Delta + \tkappamg \frac{\epsilon_f}{\Delta}, \label{eq_thm_fully_linear_noisy_constants_g}
    \end{align}
    \end{subequations}
    for all $\by\in B(\bx,\Delta)$, where $\kappamf$ and $\kappamg$ are the same as in \thmref{thm_fully_linear} (i.e.~\eqref{eq_thm_fully_linear_constants}), and 
    \begin{align}
        \tkappamf = (1+\sqrt{n})\|\hat{\bM}^{-1}\|_{\infty}, \qquad \text{and} \qquad \tkappamg = 2\tkappamf.
    \end{align}
\end{theorem}
\begin{proof}
    Noting $m(\by_i) = \tilde{f}(\by_i)$, the noisy version of \eqref{eq_linear_interp_error_tmp1} is
    \begin{align}
        |m(\by_i) - f(\bx) - \grad f(\bx)^T (\by_i-\bx)| \leq \frac{\Lgrad}{2} \beta^2 \Delta^2 + \epsilon_f,
    \end{align}
    and so the analog of \eqref{eq_linear_interp_error_tmp4} is 
    \begin{align}
        |m(\by) - f(\bx) - \grad f(\bx)^T (\by-\bx)| \leq \frac{\Lgrad}{2} (1+\sqrt{n}) \beta^2 \|\hat{\bM}^{-1}\|_{\infty} \Delta^2 + (1+\sqrt{n}) \|\hat{\bM}^{-1}\|_{\infty} \epsilon_f,
    \end{align}
    for all $\by \in B(\bx,\Delta)$.
    The result then follows from \lemref{lem_fully_linear_quadratic_from_taylor_approx}(a) with $\kappa_H=0$.
\end{proof}

A similar result holds for minimum Frobenius norm quadratic models.

\begin{theorem} \label{thm_min_frob_fully_linear_noisy}
    Suppose $f$ satisfies \minorrev{\assref{ass_smoothness_1}\ref{ass_smoothness_1_smooth}} and $\tilde{f}$ satisfies \assref{ass_noisy_oracle}, and we construct a quadratic model $m$ \eqref{eq_quadratic_model_generic} for $f$ by solving \eqref{eq_min_frob_interp_system_scaled} with function estimates $\tilde{f}(\by_i)$, where we assume $\hat{\bF}$ is invertible.
    If $\|\by_i-\bx\| \leq \beta \Delta$ for some $\beta>0$ and all $i=1,\ldots,p$, then the model satisfies \eqref{eq_thm_fully_linear_noisy_constants} for all $\by\in B(\bx,\Delta)$ with constants
    \begin{align}
        \kappamf &= \frac{\Lgrad + \kappa_H}{2} (1+\sqrt{n}) \beta^2  \|\hat{\bM}^{\dagger}\|_{\infty} + \frac{\Lgrad + \kappa_H}{2}, \qquad \text{and} \qquad \tkappamf = (1+\sqrt{n})\|\hat{\bM}^{\dagger}\|_{\infty},
    \end{align}
    with $\kappamg=2\kappamf+2\kappa_H$ and $\tkappamg=2\tkappamg$, and where
    \begin{align}
        \|\bH\| \leq \kappa_H := \frac{\Lgrad}{2}  p \beta^4 \|\hat{\bF}^{-1}\|_{\infty} + \frac{\epsilon_f}{\Delta^2} p \beta^2 \|\hat{\bF}^{-1}\|_{\infty}.
    \end{align}
\end{theorem}
\begin{proof}
    Noting $m(\by_i)=\tilde{f}(\by_i)$, the noisy equivalent of \eqref{eq_min_frob_fully_linear_tmp2} in the proof of \thmref{thm_min_frob_fully_linear} is
    \begin{align}
        |c + \bg^T (\by-\bx) - f(\bx) - \grad f(\bx)^T (\by-\bx)| &\leq \frac{\Lgrad + \|\bH\|}{2} \beta^2 (1+\sqrt{n}) \|\hat{\bM}^{\dagger}\| \Delta^2 + \epsilon_f (1+\sqrt{n}) \|\hat{\bM}^{\dagger}\|,
    \end{align}
    The argument from \thmref{thm_min_frob_fully_linear} then gives the values of the constants in \eqref{eq_thm_fully_linear_noisy_constants}.
    To get the bound on $\|\bH\|$, following the proof of \lemref{lem_min_frob_bounded_hessian} we get
    \begin{align}
        |\hat{\lambda}_i| \leq \left(\frac{\Lgrad}{2} \beta^2 \Delta^2 + \epsilon_f\right) \|\hat{\bF}^{-1}\|_{\infty}, 
    \end{align}
    from which the bound follows.
\end{proof}

\revision{Theorems~\ref{thm_fully_linear_noisy} and \ref{thm_min_frob_fully_linear_noisy} motivate the following assumption on our interpolation model, replacing \assref{ass_model_dfo}.

\begin{assumption} \label{ass_model_dfo_storm_deterministic}
    At each iteration $k$ of \algref{alg_storm_deterministic}, the model $m_k$ \eqref{eq_tr_model_dfo} satisfies:
    \begin{enumerate}[label=(\alph*)]
        \item There exist constants $\kappamf,\tkappamf,\kappamg,\tkappamg>0$ (independent of $k$) such that $m_k$ satisfies \eqref{eq_thm_fully_linear_noisy_constants} for all $\by\in B(\bx_k,\Delta_k)$;
        \item $\|\bH_k\| \leq \kappa_H - 1$ for some $\kappa_H\geq 1$ (independent of $k$). \label{ass_model_dfo_storm_deterministic_H}
    \end{enumerate}
    We again denote $\kappam := \max(\kappamf,\kappamg)$ and $\tkappam:=\max(\tkappamf,\tkappamg)$ for notational convenience.
\end{assumption}

We now consider the adaptation of \algref{alg_basic_tr_dfo} to the deterministic noise setting.
The only other change to the stated algorithm is that the numerator in the ratio test \eqref{eq_ratio_test_noisy} is now the estimated predicted decrease, based on querying $\t{f}$ (satisfying \assref{ass_noisy_oracle}) instead of $f$.
Our new algorithm is given in \algref{alg_storm_deterministic}.
}

\begin{algorithm}[tb]
\begin{algorithmic}[1]
\Require Starting point $\bx_0\in\R^n$ and trust-region radius $\Delta_0>0$. Algorithm parameters: scaling factors $0 < \gammadec < 1 < \gammainc$, acceptance threshold $0 < \eta_S < 1$, \minorrev{and} criticality threshold $\mu_c > 0$.
\For{$k=0,1,2,\ldots$}
    \State Build a local quadratic model $m_k$ \eqref{eq_tr_model_dfo} satisfying \assref{ass_model_dfo_storm_deterministic}.
    \State Solve the trust-region subproblem \eqref{eq_trs} to get a step $\bs_k$ satisfying \assref{ass_cauchy_decrease}.
    \State Calculate the ratio 
    \begin{align}
        \rho_k = \frac{\revision{\text{estimated actual decrease}}}{\text{predicted decrease}} := \frac{\revision{\t{f}(\bx_k) - \t{f}(\bx_k+\bs_k)}}{m_k(\bx_k) - m_k(\bx_k+\bs_k)}. \label{eq_ratio_test_noisy}
    \end{align}
    \If{$\rho_k \geq \eta_S$ and $\|\bg_k\| \geq \mu_c \Delta_k$} 
        \State \textit{(Successful iteration)} Set $\bx_{k+1}=\bx_k+\bs_k$ and $\Delta_{k+1}=\gammainc\Delta_k$.
    \Else
        \State \textit{(Unsuccessful iteration)} Set $\bx_{k+1}=\bx_k$ and $\Delta_{k+1}=\gammadec\Delta_k$. 
    \EndIf 
\EndFor
\end{algorithmic}
\caption{\renaming{MBDFO} trust-region method for problems with deterministic noise.}
\label{alg_storm_deterministic}
\end{algorithm}



The worst-case complexity analysis of \algref{alg_storm_deterministic} is broadly similar to the \revision{noise-free} case, but more complicated.
\revision{The iteration complexity bound matches the noise-free case (\corref{cor_wcc_dfo}) but only when the desired first-order accuracy is sufficiently large, $\bigO(\sqrt{\epsilon_f})$.}

\begin{theorem} \label{thm_deterministic_noise_wcc}
    Suppose Assumptions~\ref{ass_smoothness_1}, \ref{ass_cauchy_decrease}, \ref{ass_noisy_oracle} and \ref{ass_model_dfo_storm_deterministic} hold.
    For any
    \begin{align}
        \epsilon > \epsilon_{\min} = \revision{\bigO\left(\sqrt{(\kappam+\kappa_H) \tkappam \epsilon_f}\right)},
    \end{align}
    \revision{provided $\Delta_0$ is sufficiently large}, \algref{alg_storm_deterministic} achieves $\|\grad f(\bx_k)\| \leq \epsilon$ for the first time after at most \revision{$\bigO(\kappa_H (\kappam + \kappa_H)^2 \epsilon^{-2})$} iterations.
\end{theorem}
\begin{proof}
    \revision{See \appref{sec_storm_deterministic_proofs}.}
\end{proof}



\subsection{Stochastic Noise} \label{sec_stochastic_noise}
If the noise in our objective is stochastic, then, unlike the deterministic noise case, we can decrease the noise level by averaging many noisy values of $f(\bx)$ to achieve, in principle, arbitrarily good estimates of the true objective.
This means that, with care, we can achieve any desired level of optimality, $\|\grad f(\bx_k)\| < \epsilon$, rather than being limited to optimality levels $\epsilon=\bigO(\sqrt{\epsilon_f})$ as in the deterministic case.

Formally, we will assume that our stochastic noise \minorrev{is} unbiased, and so optimizing a function which has stochastic noise in its evaluations can be viewed as the stochastic optimization problem
\begin{align}
    \min_{\bx\in\R^n} \: f(\bx) := \E_{\omega}[f(\bx,\omega)], \label{eq_general_stochastic}
\end{align}
where $\omega\in \Omega$ is a source of random noise.
Here, we assume that we only have zeroth order \emph{stochastic oracles} $\bx \mapsto f(\bx,\omega)$, for realizations of the noise $\omega$.
To make the problem tractable, we will assume that the randomness has bounded variance.

\begin{assumption} \label{ass_noise_bounded_variance}
    There exists $\sigma>0$ such that $\Var_{\omega}(f(\bx,\omega)) \leq \sigma^2$ for all $\bx\in\R^n$.
\end{assumption}

The most important consequence of \assref{ass_noise_bounded_variance} is that we can easily construct estimates of $f(\bx)$ \minorrev{in \eqref{eq_general_stochastic} which} are arbitrarily good, with high probability.

\begin{proposition}[Chebyshev's inequality, e.g.~Chapter 2.1 of \cite{Boucheron2013}]
    Suppose \assref{ass_noise_bounded_variance} holds.
    For any $\bx\in\R^n$, suppose we generate i.i.d.~realizations $\omega_1,\ldots,\omega_N$ from $\Omega$ and calculate the sample average, $\overline{f}_N(\bx,\omega) = \frac{1}{N}\sum_{i=1}^{N} f(\bx,\omega_i)$.
    Then $\P\left[|\overline{f}_N(\bx,\omega) - f(\bx)| > t \right] \leq \frac{\sigma^2}{N t^2}$ for any $t>0$.
\end{proposition}

Roughly speaking, the Central Limit Theorem says that the typical error between the sample mean from $N$ samples and the true mean is of size $\bigO(\frac{\sigma}{\sqrt{N}})$ as $N\to\infty$.
Chebyshev's inequality gives a `finite $N$' version of this idea.
Specifically, motivated by the fully linear assumption, the below result confirms that to get a sample error of the desirable size $\bigO(\Delta^2)$ requires averaging $N=\bigO(\Delta^{-4})$ samples.

\begin{corollary} \label{cor_chebyshev}
    Suppose \assref{ass_noise_bounded_variance} holds.
    For any $\bx\in\R^n$, and $\epsilon_f,\Delta>0$ and $\alpha_f\in(0,1)$, if $N \geq \frac{\sigma^2}{\epsilon_f^2 (1-\alpha_f) \Delta^4}$ then $\P\left[|\overline{f}_N(\bx,\omega) - f(\bx)| \leq \epsilon_f \Delta^2 \right] \geq \alpha_f$.
\end{corollary}

Another consequence of Chebyshev's inequality is that, again by choosing the number of samples sufficiently large, we can guarantee that an interpolation model is fully linear with high probability.

\begin{theorem} \label{thm_sample_average_fully_linear}
    Suppose $f$ satisfies \minorrev{\assref{ass_smoothness_1}\ref{ass_smoothness_1_smooth}}, \assref{ass_noise_bounded_variance} holds, and we construct a linear model \eqref{eq_linear_model_generic} for $f$ by solving \eqref{eq_linear_interp_system_noisy} with noisy estimates $\tilde{f}(\by_i) = \overline{f}_N(\by_i,\omega)$ (i.i.d.~for each $i$), where we assume $\hat{\bM}$ is invertible.
    Assume also that $\|\by_i-\bx\| \leq \beta \Delta$ for some $\beta>0$ and all $i=1,\ldots,p$ (where $p=n+1$ is the number of interpolation points).
    If $N \geq \frac{\sigma^2}{\epsilon_f^2 (1-\alpha_m^{1/p})\Delta^4}$ for some $\epsilon_f>0$ and $\alpha_m\in(0,1)$, then with probability at least $\alpha_m$ the model is fully linear with constants
    \begin{align}
        \kappamf = \frac{\Lgrad}{2} (1+\sqrt{n}) \beta^2 \|\hat{\bM}^{-1}\|_{\infty} + \frac{\Lgrad}{2} + (1+\sqrt{n})\|\hat{\bM}^{-1}\|_{\infty} \epsilon_f, \quad \text{and} \quad \kappamg = 2\kappamf.
    \end{align}
\end{theorem}
\begin{proof}
    By assumption on $N$, \corref{cor_chebyshev} gives
    \begin{align}
        \P\left[|\overline{f}_N(\by_i,\omega) - f(\by_i)| \leq \epsilon_f \Delta^2 \right] \geq \alpha_m^{1/p},
    \end{align}
    for each $i=1,\ldots,p$, and since these events are all independent we have
    \begin{align}
        \P\left[\max_{i=1,\ldots,p}|\overline{f}_N(\by_i,\omega) - f(\by_i)| \leq \epsilon_f \Delta^2 \right] \geq \alpha_m. \label{eq_all_rhs_accurate}
    \end{align}
    That is, with probability at least $\alpha_m$ we have $|\overline{f}_N(\by_i,\omega) - f(\by_i)| \leq \epsilon_f \Delta^2$ for all $i=1,\ldots,p$.
    The result then follows from  \thmref{thm_fully_linear_noisy}.
\end{proof}

An analogous result can be derived for minimum Frobenius norm models, by the same probabilistic argument as above and the use of \thmref{thm_min_frob_fully_linear_noisy}, where we again need $N\geq \frac{\sigma^2}{\epsilon_f (1-\alpha_m^{1/p}) \Delta^4}$ samples for each objective estimate.

The above results motivate \algref{alg_storm}, a  stochastic version of \algref{alg_basic_tr_dfo}.
\algref{alg_storm} has very similar algorithmic structure to \algref{alg_basic_tr_dfo}, but with some simplifications to aid the analysis.
Specifically,
\begin{itemize}
    \item We always increase $\Delta_k$ on successful iterations, i.e.~setting $\eta_U=\eta_S$ in \algref{alg_basic_tr_dfo};
    \item We maintain the same ratio for increasing and decreasing $\Delta_k$, i.e.~setting $\gammadec = \gammainc^{-1}$ in \algref{alg_basic_tr_dfo};
    \item We impose a maximum trust-region radius, $\Delta_k \leq \Delta_{\max}$ for all $k$, similar to \algref{alg_basic_tr_dfo_2}. We assume $\Delta_{\max} = \gammainc^{j_{\max}} \Delta_0$ for some $j_{\max}\in\N$ again for simplicity.
\end{itemize}
More importantly, in  \algref{alg_storm} we now assume that the model $m_k$ is fully linear with high probability, and the approximate evaluations of $f(\bx_k)$ and $f(\bx_k+\bs_k)$ in the calculation of $\rho_k$ are accurate with high probability.
In particular, we care about the likelihood of the events
\begin{subequations}
\begin{align}
    I_k &:= \ind\{\text{$m_k$ is fully linear in $B(\bx_k,\Delta_k)$ with constants $\kappamf,\kappamg>0$}\}, \qquad \text{and} \\
    J_k &:= \ind\{\text{$|f_k^0 - f(\bx_k)| \leq \epsilon_f \Delta_k^2$ and $|f_k^s - f(\bx_k+\bs_k)| \leq \epsilon_f \Delta_k^2$ for some $\epsilon_f > 0$}\},
\end{align}
\end{subequations}
where $\ind$ is the indicator function of an event (i.e.~taking values 1 if the event occurs and 0 otherwise).
Formally, all the randomness in \algref{alg_storm} comes from the models $m_k$ and the estimates $f_k^0,f_k^s$.
So, we define the filtration $\mathcal{F}_{k-1}$ to be the $\sigma$-algebra generated by $\{m_0,f_0^0,f_0^s, \ldots, m_{k-1}, f_{k-1}^0, f_{k-1}^s\}$, representing all randomness up to the start of iteration $k$. 
Similarly, we define $\mathcal{F}_{k-1/2}$ to be generated by $\{m_0,f_0^0,f_0^s, \ldots, m_{k-1}, f_{k-1}^0, f_{k-1}^s, m_k\}$, representing all randomness up to the calculation of $\bs_k$ (i.e.~from $\mathcal{F}_{k-1}$ and $m_k$).

Given this filtration structure, we make the following assumptions about the model.

\begin{assumption} \label{ass_model_dfo_storm}
    At each iteration $k$ of \algref{alg_storm}, the model $m_k$ \eqref{eq_tr_model} satisfies:
    \begin{enumerate}[label=(\alph*)]
        \item $\P[I_k=1 | \mathcal{F}_{k-1}] \geq \alpha_m$ for some $\alpha_m\in(\frac{1}{2},1]$; \label{ass_model_dfo_storm_g}
        \item $\|\bH_k\| \leq \kappa_H - 1$ for some $\kappa_H\geq 1$ (independent of $k$). \label{ass_model_dfo_storm_H}
    \end{enumerate}
\end{assumption}

We also require the following assumption about the function value estimates used to calculate $\rho_k$.

\begin{assumption} \label{ass_evals_storm}
    At each iteration $k$ of \algref{alg_storm}, the estimates $f_k^0$ and $f_k^s$ satisfy $\P[J_k=1 | \mathcal{F}_{k-1/2}] \geq \alpha_f$ for some $\alpha_f\in(\frac{1}{2},1]$ and some $\epsilon_f > 0$.
\end{assumption}

\revision{Under \assref{ass_noise_bounded_variance}, \corref{cor_chebyshev} and \thmref{thm_sample_average_fully_linear}} show that Assumptions~\ref{ass_model_dfo_storm} and \ref{ass_evals_storm} can be satisfied by using sample averages $\overline{f}_N(\bx,\omega)$ for $N\geq \frac{\sigma^2}{\epsilon_f (1-\alpha_f^{1/2}) \Delta^4}$ samples.\footnote{We need $\alpha_f^{1/2}$ rather than $\alpha_f$ as in \corref{cor_chebyshev} because the condition $J_k=1$ requires we have accurate estimates for \emph{both} $f(\bx_k)$ and $f(\bx_k+\bs_k)$.} \revision{However, Assumptions~\ref{ass_model_dfo_storm} and \ref{ass_evals_storm} are more general and may be satisfied under other assumptions on the randomness in the objective.}

\begin{algorithm}[tb]
\begin{algorithmic}[1]
\Require Starting point $\bx_0\in\R^n$ and trust-region radius $\Delta_0>0$. Algorithm parameters: scaling factor $\gammainc > 1$, acceptance threshold $0 < \eta_S < 1$, criticality threshold $\mu_c > 0$, and maximum trust-region radius $\Delta_{\max} = \gammainc^{j_{\max}} \Delta_0$ for some $j_{\max}\in\N$.
\For{$k=0,1,2,\ldots$}
    \State Build a local quadratic model $m_k$ \eqref{eq_tr_model_dfo} satisfying \assref{ass_model_dfo_storm}.
    \State Solve the trust-region subproblem \eqref{eq_trs} to get a step $\bs_k$ satisfying \assref{ass_cauchy_decrease}.
    \State Calculate estimates $f_k^0 \approx f(\bx_k)$ and $f_k^s \approx f(\bx_k+\bs_k)$ satisfying \assref{ass_evals_storm}, and calculate the ratio 
    \begin{align}
        \rho_k = \frac{\text{estimated actual decrease}}{\text{predicted decrease}} := \frac{f_k^0 - f_k^s}{m_k(\bx_k) - m_k(\bx_k+\bs_k)}.
    \end{align}
    \If{$\rho_k \geq \eta_S$ and $\|\bg_k\| \geq \mu_c \Delta_k$} 
        \State \textit{(Successful iteration)} Set $\bx_{k+1}=\bx_k+\bs_k$ and $\Delta_{k+1}=\min(\gammainc\Delta_k, \Delta_{\max})$. 
    \Else
        \State \textit{(Unsuccessful iteration)} Set $\bx_{k+1}=\bx_k$ and $\Delta_{k+1}=\gammainc^{-1}\Delta_k$. 
    \EndIf 
\EndFor
\end{algorithmic}
\caption{Stochastic \renaming{MBDFO} trust-region method for solving \eqref{eq_general_stochastic}.}
\label{alg_storm}
\end{algorithm}

The worst-case complexity of \algref{alg_storm} is derived from the following general probabilistic result.
\revision{This considers a non-negative random process $\phi_k$ which (before stopping) typically decreases in proportion to another non-negative random process $\Delta_k$. It says that, provided $\Delta_k$ is unlikely to get small before stopping (specifically, increases in $\Delta_k$ occur when a biased event $\P[W_k=1] > 1/2$ occurs), we can bound how long it takes for stopping to occur.}

\begin{proposition}[Theorem 2, \cite{Blanchet2019}] \label{prop_supermartingale}
    Suppose we have random processes $\{(\phi_k,\Delta_k,W_k)\}_{k=0}^{\infty}$ with $\phi_k, \Delta_k \geq 0$, and 
    \begin{align}
        \P[W_{k} = 1 | \mathcal{F}_{k-1}] = q, \qquad \text{and} \qquad \P[W_{k} = -1 | \mathcal{F}_{k-1}] = 1-q, \label{eq_storm_W_defn}
    \end{align}
    for some $q>\frac{1}{2}$, where $\mathcal{F}_{k-1}$ is the $\sigma$-algebra generated by $\{(\phi_j,\Delta_j)\}_{j=0}^{k}$ and $\{W_j\}_{j=0}^{k-1}$.
    For any $\epsilon>0$, let $K_{\epsilon}$ be a stopping time with respect to the filtrations $\{\mathcal{F}_k\}_{k=0}^{\infty}$.
    If, for all $k$,
    \begin{enumerate}[label=(\alph*)]
        \item There exists $\lambda>0$ and $j_{\max}\in\Z$ such that $\Delta_k \leq \Delta_{\max} := \Delta_0 e^{\lambda j_{\max}}$; \label{item_prop_supermartingale1}
        \item There exists $j_{\epsilon}\in\Z$ with $j_{\epsilon} \leq 0$ such that $K_{\epsilon} > k$ implies $\Delta_{k+1} \geq \min(\Delta_k e^{\lambda W_{k}}, \Delta_{\epsilon})$, where $\Delta_{\epsilon} := \Delta_0 e^{\lambda j_{\epsilon}}$ and $\lambda>0$ is from \ref{item_prop_supermartingale1}; and \label{item_prop_supermartingale2}
        \item There exists $\theta>0$ and non-increasing function $h:[0,\infty) \to [0,\infty)$ such that $K_{\epsilon} > k$ implies $\phi_k - \E[\phi_{k+1} | \mathcal{F}_{k-1}] \geq \theta h(\Delta_k)$; \label{item_prop_supermartingale3}
    \end{enumerate}
    then
    \begin{align}
        \E[K_{\epsilon}] \leq \frac{q \phi_0}{2(q-\frac{1}{2}) \theta h(\Delta_{\epsilon})} + 1.
    \end{align}
\end{proposition}

\propref{prop_supermartingale} can be thought of as a stochastic version of \thmref{thm_wcc_dfo}, where condition \ref{item_prop_supermartingale2} is a requirement similar to the conclusion of \lemref{lem_delta_min_dfo}.\footnote{We will use $h(\Delta)=\Delta^2$ here, but $h(\Delta)=\Delta^3$ can be used to get second-order complexity results.}

To apply \propref{prop_supermartingale} to analyze \algref{alg_storm}, we will take \revision{$\phi_k$ to be a measure of algorithm progress based on $f(\bx_k)-\flow$ (see \lemref{lem_storm_decrease} below)} and the stopping time to be 
\begin{align}
    K_{\epsilon} := \inf\{k \geq 0 : \|\grad f(\bx_k)\| < \epsilon\}   
\end{align}
and define the \revision{event $W_k$ to be when all probabilistic estimates are accurate}
\begin{align}
    W_k := 2 I_k J_k - 1 = \begin{cases} 1, & \text{if $I_k=J_k=1$}, \\ -1, & \text{otherwise}, \end{cases}
\end{align}
and so from Assumptions~\ref{ass_model_dfo_storm}\ref{ass_model_dfo_storm_g} and \ref{ass_evals_storm} we get \eqref{eq_storm_W_defn} with $q=\alpha_m \alpha_f$ (hence we will need to assume that $\alpha_m \alpha_f > \frac{1}{2}$).
Condition \ref{item_prop_supermartingale1} in \propref{prop_supermartingale} is also automatically satisfied for \algref{alg_storm} by setting $\lambda = \log(\gammainc)$.
We now need to establish that conditions \ref{item_prop_supermartingale2} and \ref{item_prop_supermartingale3} hold, for suitable choices of $j_{\epsilon}$, $\phi_k$ and $\theta$.

\begin{lemma} \label{lem_storm_delta_large}
    Suppose Assumptions~\minorrev{\ref{ass_smoothness_1}\ref{ass_smoothness_1_smooth}}, \ref{ass_cauchy_decrease},  \ref{ass_model_dfo_storm} and \ref{ass_evals_storm} hold and we run \algref{alg_storm}.
    For any $\epsilon>0$, if $j_{\epsilon}$ is the largest non-positive integer (i.e.~smallest in magnitude) such that
    \begin{align}
        \Delta_{\epsilon} := \Delta_0 e^{\lambda j_{\epsilon}} \leq \frac{\epsilon}{c_0}, \qquad \text{where} \qquad c_0 := \max\left(\frac{4\max(\kappamf, \epsilon_f)}{\kappa_s (1-\eta_S)}, \kappa_H, \mu_c\right) + \kappamg,
    \end{align}
    and where $\lambda := \log(\gammainc)$, then provided $K_{\epsilon} > k$ we have $\Delta_{k+1} \geq \min(\Delta_k e^{\lambda W_{k}}, \Delta_{\epsilon})$.
\end{lemma}
\begin{proof}
    Fix $k$ such that $K_{\epsilon} > k$, and so $\|\grad f(\bx_k)\| \geq \epsilon$.
    If $W_k=-1$ (i.e.~$I_k=0$ or $J_k=0$), then $\Delta_{k+1} = \gammainc^{-1} \Delta_k = \Delta_k e^{-\lambda} = \Delta_k e^{\lambda W_k}$ and the result holds.
    So now suppose that $W_k=1$ (i.e.~$I_k=J_k=1$).

    Since $\Delta_k = \gammainc^{i_k} \Delta_0$ for some $i_k\in\Z$ with $i_{k+1} \in\{i_k-1, i_k+1\}$ (by the trust-region updating mechanism), if $\Delta_k > \Delta_{\epsilon}$ then $i_k > j_{\epsilon}$ and so $i_{k+1}\geq j_{\epsilon}$, which means $\Delta_{k+1} \geq \Delta_{\epsilon}$ and the result holds.
    
    The only remaining case to consider is where $I_k=J_k=1$ and $\Delta_k \leq \Delta_{\epsilon}$.
    Since $\|\grad f(\bx_k)\| \geq \epsilon$ and $\Delta_k \leq \Delta_{\epsilon} \leq c_0^{-1} \epsilon$, we have $\|\grad f(\bx_k)\| \geq c_0 \Delta_k$, and so since $m_k$ is fully linear ($I_k=1$),
    \begin{align}
        \|\bg_k\| \geq \|\grad f(\bx_k)\| - \kappamg \Delta_k \geq \max\left(\frac{4\max(\kappamf, \epsilon_f)}{\kappa_s (1-\eta_S)}, \kappa_H, \mu_c\right) \Delta_k,
    \end{align} 
    and so $\|\bg_k\| \geq \mu_c \Delta_k$ and
    \begin{align}
        \Delta_k \leq \min\left(\frac{\kappa_s (1-\eta_S)}{2(\kappamf + \epsilon_f)}, \frac{1}{\kappa_H}\right) \|\bg_k\|.
    \end{align}
    From \assref{ass_cauchy_decrease} and $\Delta_k \leq \|\bg_k\| / \kappa_H$ we have $m_k(\bx_k) - m_k(\bx_k+\bs_k) \geq \kappa_s \|\bg_k\| \Delta_k$.
    Separately, since $I_k=1$ the model $m_k$ is fully linear, and since $J_k=1$ we have $\max(|f_k^0 - f(\bx_k)|, |f_k^s - f(\bx_k+\bs_k)|) \leq \epsilon_f \Delta_k^2$.
    Thus
    \begin{align}
        |\rho_k - 1| &= \frac{|(f_k^0 - m_k(\bx_k)) - (f_k^s - m_k(\bx_k+\bs_k))|}{m_k(\bx_k) - m_k(\bx_k+\bs_k)}, \\
        &\leq \frac{|f_k^0 - f(\bx_k)| + |f(\bx_k) - m_k(\bx_k)| + |f_k^s - f(\bx_k+\bs_k)| + |f(\bx_k+\bs_k) - m_k(\bx_k+\bs_k)|}{m_k(\bx_k) - m_k(\bx_k+\bs_k)}, \\
        &\leq \frac{2(\kappamf + \epsilon_f) \Delta_k^2}{\kappa_s \|\bg_k\| \Delta_k},
    \end{align}
    and so $|\rho_k-1| \leq 1-\eta_S$, or $\rho_k \geq \eta_S$, and iteration $k$ is successful.
    Hence $\Delta_{k+1} = \gammainc \Delta_k = \Delta_k e^{\lambda W_k}$ from $W_k=1$, using $\Delta_k \leq \Delta_{\epsilon} \leq \gammainc \Delta_{\max}$ since $j_{\epsilon} \leq 0 < j_{\max}$.
\end{proof}

We now establish condition \ref{item_prop_supermartingale3} in \propref{prop_supermartingale}, for a suitable choice of $\phi_k$ and $h(\Delta)$.

\begin{lemma} \label{lem_storm_decrease}
    Suppose Assumptions~\ref{ass_smoothness_1}, \ref{ass_cauchy_decrease},  \ref{ass_model_dfo_storm} and \ref{ass_evals_storm} hold, we run \algref{alg_storm} and 
    \begin{align}
         \epsilon_f < \frac{1}{2} \eta_S \kappa_s \mu_c \min\left(1, \frac{\mu_c}{\kappa_H}\right), \label{eq_storm_epsf_suff_small}
    \end{align}
    holds.
    Also assume that $\alpha_m$ and $\alpha_f$ in Assumptions~\ref{ass_model_dfo_storm} and \ref{ass_evals_storm} are sufficiently close to 1. Then there exists $\theta>0$ and $\nu\in(0,1)$ such that, for any $\epsilon>0$, provided $K_{\epsilon} > k$ we have 
    \begin{align}
        \phi_k - \E[\phi_{k+1} | \mathcal{F}_{k-1}] \geq \theta \Delta_k^2, \qquad \text{where} \qquad \phi_k := \nu [f(\bx_k)-\flow] + (1-\nu) \Delta_k^2 \geq 0.
    \end{align}
\end{lemma}
\begin{proof}
    The proof is omitted for brevity; see \cite[Theorem 3]{Blanchet2019} with further details given in \cite[Theorem 4.11]{Chen2018}.
\end{proof}

We now get our main complexity bound as an immediate consequence of \propref{prop_supermartingale}.

\begin{theorem} \label{thm_wcc_storm}
    Suppose Assumptions~\ref{ass_smoothness_1}, \ref{ass_cauchy_decrease},  \ref{ass_model_dfo_storm} and \ref{ass_evals_storm} hold, we run \algref{alg_storm} and \eqref{eq_storm_epsf_suff_small} holds.
    Also assume that $\alpha_m$ and $\alpha_f$ in Assumptions~\ref{ass_model_dfo_storm} and \ref{ass_evals_storm} are sufficiently close to 1 (the same as required by \lemref{lem_storm_decrease}). 
    Then there exists a constant $\tilde{\theta}>0$ such that
    \begin{align}
        \E[K_{\epsilon}] \leq \frac{\alpha_m \alpha_f [f(\bx_0)-\flow + \Delta_0^2]}{2(\alpha_m \alpha_f - \frac{1}{2}) \tilde{\theta} \epsilon^2} + 1.
    \end{align}
\end{theorem}
\begin{proof}
    Apply Lemmas~\ref{lem_storm_delta_large} and \ref{lem_storm_decrease} to show that \propref{prop_supermartingale} applies with $h(\Delta) = \Delta^2$.
    Then use $\phi_0 = \nu [f(\bx_0)-\flow] + (1-\nu) \Delta_0^2 \leq f(\bx_0) -\flow + \Delta_0^2$ from $\nu\in(0,1)$.
    Since $j_{\epsilon}$ in the definition of $\Delta_{\epsilon}$ in \lemref{lem_storm_delta_large} is chosen to be smallest in magnitude, we have $\Delta_{\epsilon} \geq \frac{\epsilon}{\gammainc c_0}$, and so $\theta \Delta_{\epsilon}^2 \leq \tilde{\theta} \epsilon^2$ for $\tilde{\theta} = \theta / (\gammainc^2 c_0^2)$.
\end{proof}

In general, we know that at each iteration $k$ we need to sample the stochastic objective $\bigO(\Delta_k^{-4})$ times to satisfy our probabilistic Assumptions~\ref{ass_model_dfo_storm} and \ref{ass_evals_storm}.
From \lemref{lem_storm_delta_large} we know that $\Delta_k$ can typically get as small as $\Delta_{\epsilon} = \bigO(\epsilon)$.
So, broadly speaking, we may require up to $\bigO(\epsilon^{-4})$ stochastic objective evaluations per iteration, or $\bigO(\epsilon^{-6} / (\alpha_m \alpha_f -\frac{1}{2}))$ evaluations in total.
A rigorous $\bigO(\epsilon^{-6})$ evaluation complexity bound for a variant of \algref{alg_storm} is proven in \cite[Theorem 5]{Jin2025}.

\begin{remark} \label{rem_storm_expensive}
    \revision{Using \corref{cor_chebyshev} and \thmref{thm_sample_average_fully_linear} to satisfy Assumptions~\ref{ass_model_dfo_storm} and \ref{ass_evals_storm}} relies on repeatedly sampling the objective at every iterate and interpolation point, potentially a large number of times. \revision{As such, \algref{alg_storm}} is the only method in this work that is not well-suited to the regime where objective evaluations are expensive.
    However, this is the price of seeking convergence to a stationary point (which is not generically possible in the case of deterministic noise, as the true objective value can never be known for any point).
    Convergence to a neighborhood, as in the case of deterministic noise, can be established for  \algref{alg_storm_deterministic} \cite{Cao2024}, which is more realistic for the expensive evaluation regime.
\end{remark}

\subsubsection*{Notes and References}
{\small \lemref{lem_noisy_fin_diff} is based on \cite[Lemma 9.1]{Nocedal2006}. \algref{alg_storm_deterministic} and the associated worst-case complexity analysis (\thmref{thm_deterministic_noise_wcc}) is \revision{new, but broadly follows the approach from \cite{Cao2024,Chaudhry2025}}. 
\revision{In particular, \cite{Cao2024} considers very flexible model and evaluation accuracy assumptions, including allowing for probabilistically accurate models and sufficient decrease estimates, similar to \secref{sec_stochastic_noise}.} 
It also includes second-order complexity theory.
\revision{This requires modifying the ratio test \eqref{eq_ratio_test_noisy} based on knowledge of $\epsilon_f$. If the noise were stochastic, we could estimate $\epsilon_f$ based on the standard deviation of $\tilde{f}(\bx)$ for different values of $\bx$.
For deterministic noise, $\epsilon_f$ must be estimated in a different way (e.g.~\cite{More2011}).}

A procedure to estimate the size of (especially deterministic) noise in a function was developed in \cite{More2011}. 
Accurate finite differencing schemes (i.e.~picking the perturbation size $h$ based on the estimated noise level) are developed in \cite{More2012,Shi2022a}.
Several recent works have shown how to adapt existing (derivative-based) algorithms to appropriately handle deterministic noise in objective and gradient evaluations, such as \cite{Shi2022,Sun2023,Oztoprak2023}. 
They typically also rely on a modified ratio test, similar to the use of \eqref{eq_ratio_test_noisy} in \algref{alg_storm_deterministic}.
A modified ratio test is also used in the (derivative-based) trust-region solver TRU in the GALAHAD package \cite{Gould2003} to handle roundoff errors.
If the level of deterministic noise can be controlled (i.e.~$\t{f}(\bx)$ can be computed for any $\epsilon_f>0$, but not $\epsilon_f=0$), then   \algref{alg_basic_tr_dfo} with minor adjustments can converge to any accuracy level \cite{Ehrhardt2021}.

A more widely studied setting for stochastic optimization is where stochastic gradient estimates, $\bx \mapsto \grad_{\bx} f(\bx,\omega)$ are available, which is the core of many algorithms for training machine learning models \cite{Bottou2018}.
We also note that stricter assumptions about the distribution of $f(\bx,\omega)$, e.g.~bounds on higher-order moments, can yield tighter bounds on $|\overline{f}_N(\bx,\omega) - f(\bx)|$ than Chebyshev's inequality (e.g.~Hoeffding's or Bernstein's inequality), see \cite{Boucheron2013}.

\algref{alg_storm} was introduced in \cite{Chen2018}, where almost-sure convergence of $\|\grad f(\bx_k)\|$ to zero was shown, and the worst-case complexity analysis shown here was originally given in \cite{Blanchet2019}.
The work \cite{Blanchet2019} also includes an extension of \algref{alg_storm} that converges to second-order critical points, with suitable worst-case complexity bounds.
\revision{This approach has been extended to constrained problems via stochastic SQP methods (primarily aimed at the stochastic gradient setting but also applicable to \renaming{MBDFO}) in \cite{Fang2024,Fang2026complexity,Fang2026}.}
Adaptive sample averaging can improve performance, where the number of samples $N$ is chosen dynamically at every evaluation point, based on the sample standard deviation and $\Delta_k$, to ensure that the current estimate $\overline{f}_N(\bx,\omega)$ is sufficiently accurate \cite{Shashaani2018,Ha2024}.
The sample complexity (i.e.~total number of stochastic oracle calls) of stochastic \renaming{MBDFO} algorithms can be reduced under stronger assumptions on the stochastic oracle, such as bounds on tail probabilities \cite{Rinaldi2024} or the availability of common random numbers \cite{Ha2024a}.

Instead of constructing interpolation models using averaged estimates $\overline{f}_N(\bx,\omega)$ for large $N$, an alternative is to use regression models (i.e.~$N=1$ but with very large $p$), based on the interpolation theory from \cite{Conn2008,Conn2009,Billups2013} and discussed briefly in \secref{sec_min_frob}.
In this case, the almost-sure convergence of $\|\grad f(\bx_k)\|$ to zero is shown in \cite{Larson2016}.
Regression models are useful for both deterministic and stochastic noise (unlike sample averaging), but sample averaging allows for more accurate estimates $f_k^0,f_k^s$ for deciding whether to accept a step or not.
Radial basis function models can also cope with large numbers of interpolation points $p$, and hence construct regression-type models suitable for noisy problems \cite{Augustin2017}.
Error bounds for linear regression models as the number of points $p\to\infty$, relevant for problems with stochastic noise, are derived in \cite{Hare2023}.}

{\small A growing body of recent work extends the results from \cite{Chen2018,Blanchet2019} to more general (probabilistic) conditions on function value and model accuracy \cite{Cao2024,Jin2025,Fang2024,Fang2026complexity}, \revision{and to direct search methods \cite{Audet2021,Dzahini2023}}. This is improving our understanding of how much progress can be made in \renaming{MBDFO} and similar methods for stochastic problems with given levels of accuracy. A greater understanding is useful here, since \renaming{MBDFO} algorithms for deterministic problems with minimal/no sample averaging can already perform quite well in some cases \cite[Section 5]{Cartis2019a}.}

\section{Summary and Software} \label{sec_conclusion}

The creation and analysis of \renaming{MBDFO} algorithms requires an interesting combination of optimization and approximation theory, which has been extensively developed over the last 30 years in particular.
We have introduced the most important tools used in these algorithms, in particular trust-region methods and the construction of fully linear/\minorrev{fully} quadratic interpolation models, and have demonstrated how to extend these ideas to important settings such as constrained problems and noisy objective evaluations.
As mentioned in the introduction, the intention was never to provide a comprehensive overview of developments in the field; for this, we direct the reader to \cite{Larson2019}, and also note the more detailed discussions of algorithm and interpolation theory in \cite{Conn2009}.

\subsection{Software}
We finish by providing a list of high-quality open-source software implementations of \renaming{MBDFO} methods, suitable for use in practical applications as well as a starting point for further research.\footnote{Disclosure: \minorrev{The author is} the primary developer of two packages listed here, namely Py-BOBYQA and DFO-LS.}
\begin{description}
    \item[Software of M.~J.~D.~Powell] A collection of several Fortran packages for \renaming{MBDFO}. This includes COBYLA \cite{Powell1994} (general constrained problems, linear interpolation), UOBYQA \cite{Powell2002} \& NEWUOA \cite{Powell2006} (unconstrained problems, fully and minimum Frobenius norm quadratic models respectively), BOBYQA \cite{Powell2009} (bound-constrained problems, minimum Frobenius norm quadratic models) and LINCOA \cite{Powell2015} (linearly constrained problems, minimum Frobenius norm quadratic models). They are most easily accessed through the PRIMA package \cite{Zhang2023}, which also includes C, MATLAB, Python and Julia interfaces. Py-BOBYQA \cite{Cartis2019a} is a separate, pure Python re-implementation of BOBYQA with additional heuristics to improve performance for noisy problems. \revision{An accessible overview of these algorithms can be found in \cite{Ragonneau2024}.}\footnote{Mike Powell (1936--2015) was a pioneer of both (derivative-based) trust-region and \renaming{MBDFO} methods, among many other significant contributions to optimization and approximation theory \cite{Buhmann2017}.}
    \item[COBYQA \cite{Ragonneau2022}] A general-purpose Python package \minorrev{which can handle} bound and nonlinear constraints. Nonlinear constraints are also assumed to be derivative-free. It is based on minimum Frobenius norm quadratic models. It is most readily available through SciPy's optimization module.\footnote{See \url{https://scipy.org/}.}
    \item[DFO-LS \cite{Cartis2019a}] A Python package for nonlinear least-squares problems. It can handle simple bound and convex constraints, and nonsmooth regularizers, and it includes heuristics to improve performance for noisy problems. It is based on linear interpolation models for each term in the least-squares sum (see \secref{sec_composite_models}).
    \item[IBCDFO] A collection of MATLAB and Python packages for \renaming{MBDFO} with composite models (see \secref{sec_composite_models}). Includes POUNDerS \cite{Wild2017} (nonlinear least-squares with bound constraints using minimum Frobenius norm quadratic models), manifold sampling \cite{Larson2021} \& GOOMBAH \cite{Larson2024} (both for unconstrained problems with general composite objectives \eqref{eq_composite_generic} with $h$ possibly nonsmooth, using minimum Frobenius norm quadratic models) 
    \item[ASTRO-DF \cite{Ha2021}] A Python package aimed at stochastic \renaming{MBDFO} problems. It can handle constraints and is well-suited to simulation-based optimization, where common random numbers can be used to evaluate the objective at different points in a consistent way. It is based on underdetermined quadratic interpolation models (using diagonal Hessians as per \corref{cor_wcc_dfo_min_frob}).
    \item[NOMAD \cite{Audet2022}] A widely used direct search code, written in C++ with MATLAB and Python interfaces. 
    \revision{It constructs quadratic models with minimum Frobenius norm interpolation (or an alternative surrogate) to accelerate the direct search process via a \renaming{model-based} search step, and/or by ordering the points to be checked.}
    It can handle nonsmooth problems, constraints, discrete variables, multi-objective problems and more. 
\end{description}
Table~\ref{tab_software_links} provides a list of URLs for the above packages, accurate at \minorrev{the} time of writing.
This list is not intended to be exhaustive, but hopefully can provide a starting point for the interested reader.

\begin{table}[tb]
    \centering
    {\small
    \begin{tabular}{ll}
        \hline Package & Location \\ \hline
        PRIMA (Powell's software) & \url{https://github.com/libprima/prima} \\
        Py-BOBYQA & \url{https://github.com/numericalalgorithmsgroup/pybobyqa} \\ 
        COBYQA & \url{https://github.com/cobyqa/cobyqa} \\ 
        DFO-LS & \url{https://github.com/numericalalgorithmsgroup/dfols} \\
        IBCDFO & \url{https://github.com/POptUS/IBCDFO} \\
        ASTRO-DF & \url{https://github.com/simopt-admin/simopt} \\
        NOMAD & \url{https://github.com/bbopt/nomad} \\ \hline
    \end{tabular}
    }
    \caption{URLs for selected \renaming{MBDFO} software packages.}
    \label{tab_software_links}
\end{table}

\subsection{Outlook}
\renaming{MBDFO} theory and algorithms have matured greatly over the last 30 years, but there is still much ongoing work to be done.
Compared to many other areas of optimization, the gap between theoretically studied algorithms and efficient, practical software is relatively large, and it remains to be seen how closely aligned these can be made.
As local \renaming{MBDFO} methods mature, there is more scope to incorporate them into global optimization techniques such as multistart methods \cite{Larson2018,Jaiswal2024}, since \renaming{MBDFO} methods tend to outperform global optimization methods when good estimates of minimizers are available \cite{Rios2013}.
Another consequence of the increasing maturity is that \renaming{MBDFO} can now start to be applied to optimization problems with more complex structure, such as multi-objective optimization \cite{Custodio2011,Liuzzi2025}, optimization on manifolds \cite{Najafi2026} and mixed integer problems \cite{Torres2024,Kimiaei2025}.

\paragraph{Acknowledgments}
This work was supported by the Australian Research Council Discovery Early Career Award DE240100006.
Thanks to Warren Hare for the encouragement to pursue this project. 
Thanks to \minorrev{Nicole Felice}, Warren Hare, Jeffrey Larson, Fangyu Liu, Matt Menickelly, Cl\'ement Royer and Stefan Wild \minorrev{for spotting errors and providing helpful feedback}. 
\minorrev{Thanks to the editor and anonymous referees for their suggestions on the manuscript.}

\begin{small}
\bibliographystyle{siam}
\bibliography{refs} 
\end{small}

\appendix

\section{Technical Results} \label{app_technical_results}

Here, we collect some technical results used in the main text.

\begin{lemma} \label{lem_difference_of_linear}
    Suppose we have two linear functions $m_i:\R^n\to\R$ for $i\in\{1,2\}$, defined as $m_i(\by) := c_i + \bg_i^T (\by-\bx)$, such that $|m_1(\by)-m_2(\by)| \leq \epsilon$ for all $\by\in B(\bx,\Delta)$, for some $\bx\in\R^n$ and $\Delta>0$.
    Then
    \begin{align}
        |c_1 - c_2| \leq \epsilon, \qquad \text{and} \qquad \|\bg_1 - \bg_2\| \leq \frac{2\epsilon}{\Delta}.
    \end{align}
\end{lemma}
\begin{proof}
    First, we have $|c_1-c_2| = |m_1(\bx) - m_2(\bx)| \leq \epsilon$.
    Next, if $\bg_1 \neq \bg_2$, choose $\by = \bx + \Delta\frac{\bg_1-\bg_2}{\|\bg_1-\bg_2\|}$, to ensure that $\|\by-\bx\|=\Delta$ and $\by-\bx$ is parallel to $\bg_1-\bg_2$.
    This gives
    \begin{align}
        \Delta \|\bg_1-\bg_2\| = |(\bg_1-\bg_2)^T (\by-\bx)| \leq |m_1(\by) - m_2(\by)| + |c_1 - c_2| \leq 2\epsilon.
    \end{align}
    If instead $\bg_1=\bg_2$ then the bound on $\|\bg_1-\bg_2\|$ is trivial.
\end{proof}

\begin{lemma} \label{lem_difference_of_quadratic}
    Suppose we have two quadratic functions $m_i:\R^n\to\R$ for $i\in\{1,2\}$, defined as $m_i(\by) := c_i + \bg_i^T (\by-\bx) + \frac{1}{2}(\by-\bx)^T \bH_i (\by-\bx)$, such that $|m_1(\by)-m_2(\by)| \leq \epsilon$ for all $\by\in B(\bx,\Delta)$, for some $\bx\in\R^n$ and $\Delta>0$.
    Then
    \begin{align}
        |c_1 - c_2| \leq \epsilon, \qquad \|\bg_1-\bg_2\| \leq \frac{10 \epsilon}{\Delta}, \qquad \text{and} \qquad \|\bH_1 - \bH_2\| \leq \frac{24\epsilon}{\Delta^2}.
    \end{align}
\end{lemma}
\begin{proof}
    First, we have $|c_1 - c_2| = |m_1(\bx) - m_2(\bx)| \leq \epsilon$.
    Hence for any $\by\in B(\bx,\Delta)$ we have
    \begin{align}
        \left|(\by-\bx)^T\left[\bg_1 + \frac{1}{2}\bH_1 (\by-\bx) - \bg_2 - \frac{1}{2}\bH_2 (\by-\bx)\right]\right| \leq |m_1(\by) - m_2(\by)| + |c_1 - c_2| \leq 2\epsilon. \label{eq_difference_of_quadratic_tmp1}
    \end{align}
    Now, define $\hat{\bu} := \frac{\bg_1-\bg_2}{\|\bg_1-\bg_2\|}$ if $\bg_1\neq \bg_2$, or any unit vector otherwise.
    Similarly, define $\hat{\bv}$ to be a unit eigenvector corresponding to the largest eigenvalue in magnitude of $\bH_1-\bH_2$.
    Hence we have $\hat{\bu}^T (\bg_1-\bg_2) = \|\bg_1-\bg_2\|$ and $|\hat{\bv}^T (\bH_1 - \bH_2)\hat{\bv}| = \|\bH_1 - \bH_2\|$.

    Now applying \eqref{eq_difference_of_quadratic_tmp1} to $\by=\bx+\frac{\Delta}{2}\hat{\bu}$ and $\by=\bx+\Delta\hat{\bv}$, we get 
    \begin{align}
        \left|\frac{\Delta}{2}\underbrace{\hat{\bu}^T (\bg_1-\bg_2)}_{=\|\bg_1-\bg_2\|} + \frac{\Delta^2}{8}\hat{\bu}^T (\bH_1-\bH_2)\hat{\bu} \right| &\leq 2\epsilon, \qquad \text{and} \label{eq_diff_quad_tmp1}  \\
        \left|\Delta \hat{\bv}^T (\bg_1 - \bg_2) + \frac{\Delta^2}{2}\underbrace{\hat{\bv}^T (\bH_1-\bH_2)\hat{\bv}}_{=\pm\|\bH_1-\bH_2\|}\right| &\leq 2\epsilon, \label{eq_diff_quad_tmp2}
    \end{align}
    respectively.
    The first inequality \eqref{eq_diff_quad_tmp1}, gives us
    \begin{align}
        \frac{\Delta}{2}\|\bg_1-\bg_2\| - \frac{\Delta^2}{8} |\hat{\bu}^T (\bH_1-\bH_2)\hat{\bu}| &\leq \left| \frac{\Delta}{2}\|\bg_1-\bg_2\| - \frac{\Delta^2}{8} \left|-\hat{\bu}^T (\bH_1-\bH_2)\hat{\bu}\right| \right|, \\
        &\leq \left|\frac{\Delta}{2}\|\bg_1-\bg_2\| + \frac{\Delta^2}{8} \hat{\bu}^T (\bH_1-\bH_2)\hat{\bu}\right| \leq 2\epsilon,
    \end{align}
    where the second inequality follows from the reverse triangle inequality.
    For the second condition \eqref{eq_diff_quad_tmp2}, we first suppose $\hat{\bv}^T (\bH_1-\bH_2)\hat{\bv} \geq 0$, so $\hat{\bv}^T (\bH_1-\bH_2)\hat{\bv} = \|\bH_1-\bH_2\|$. 
    In that case,
    \begin{align}
        \frac{\Delta^2}{2}\|\bH_1-\bH_2\| &= \left(\Delta \hat{\bv}^T (\bg_1-\bg_2) + \frac{\Delta^2}{2}\hat{\bv}^T (\bH_1-\bH_2)\hat{\bv}\right) - \Delta \hat{\bv}^T (\bg_1-\bg_2) \leq 2\epsilon + \Delta |\hat{\bv}^T (\bg_1-\bg_2)|,
    \end{align}
    and in the other case $\hat{\bv}^T (\bH_1-\bH_2)\hat{\bv} < 0$, for which $\hat{\bv}^T (\bH_1-\bH_2)\hat{\bv} = -\|\bH_1-\bH_2\|$, we reach the same conclusion via
    \begin{align}
        \frac{\Delta^2}{2}\|\bH_1-\bH_2\| &= -\left(\Delta \hat{\bv}^T (\bg_1-\bg_2) + \frac{\Delta^2}{2}\hat{\bv}^T (\bH_1-\bH_2)\hat{\bv}\right) + \Delta \hat{\bv}^T (\bg_1-\bg_2) \leq 2\epsilon + \Delta |\hat{\bv}^T (\bg_1-\bg_2)|.
    \end{align}
    Since $|\hat{\bv}^T (\bg_1-\bg_2)| \leq \|\bg_1-\bg_2\|$ by Cauchy-Schwarz and $|\hat{\bu}^T (\bH_1-\bH_2)\hat{\bu}| \leq \|\bH_1-\bH_2\|$ from Rayleigh quotients, we ultimately conclude
    \begin{align}
        \frac{\Delta}{2}\|\bg_1-\bg_2\| \leq 2\epsilon + \frac{\Delta^2}{8}\|\bH_1 - \bH_2\|, \qquad \text{and} \qquad \frac{\Delta^2}{2} \|\bH_1-\bH_2\| \leq 2\epsilon + \Delta \|\bg_1-\bg_2\|.
    \end{align}
    The first of these conditions implies $\frac{\Delta^2}{2}\|\bH_1-\bH_2\| \geq 2\Delta\|\bg_1-\bg_2\| - 8\epsilon$, and so $2\Delta\|\bg_1-\bg_2\| - 8\epsilon \leq \frac{\Delta^2}{2}\|\bH_1-\bH_2\| \leq 2\epsilon + \Delta\|\bg_1-\bg_2\|$.
    This gives $\Delta \|\bg_1-\bg_2\| \leq 10\epsilon$.
    Lastly, we apply $\frac{\Delta^2}{2}\|\bH_1-\bH_2\| \leq 2\epsilon + \Delta\|\bg_1-\bg_2\| \leq 12\epsilon$, \minorrev{and we get the desired result}.
\end{proof}

\section{Poisedness of Structured Fully Quadratic Models} \label{app_fully_quadratic_lebesgue}

Here, \minorrev{we explicitly estimate the poisedness constant $\Lambda$} for the structured fully quadratic interpolation set \eqref{eq_fully_quadratic_example_points}, namely
\begin{align}
    \mathcal{Y} &= \{\bx, \bx + \Delta\be_1, \ldots, \bx + \Delta\be_n, \bx - \Delta\be_1, \ldots, \bx-\Delta\be_n, \nonumber \\ 
    & \bx + \Delta(\be_1+\be_2), \ldots, \bx + \Delta (\be_1 + \be_n), \bx + \Delta(\be_2 + \be_3), \ldots, \bx + \Delta(\be_{n-1} + \be_n)\}. 
\end{align}
We may without loss of generality assume $\bx=\bm{0}$ and $\Delta=1$, since the values of the Lagrange polynomials are invariant to shifts and scalings.
That is, we take our interpolation set to be
\begin{align}
    \mathcal{Y} &= \{\bm{0}\} \cup \{\pm\be_i : i=1,\ldots,n\} \cup \{ \be_i + \be_j : i,j = 1,\ldots,n, \: j > i\}. 
\end{align}

We recall the proof of \lemref{lem_quadratic_sample_estimate}, which explicitly computes the coefficients for the associated interpolation linear system: if we wish to interpolate values $m(\by) = f(\by)$ for all $\by\in \mathcal{Y}$, then (with our normalization $\bx=\bm{0}$ and $\Delta=1$) we get $c=f(\bm{0})$, $g_i = \frac{1}{2}(f(\be_i) - f(-\be_i))$ and $H_{i,i} = f(\be_i) + f(-\be_i) - 2f(\bm{0})$ for $i=1,\ldots,n$, and $H_{i,j}=f(\be_i+\be_j) - c - g_i - g_j - \frac{1}{2} H_{i,i} - \frac{1}{2}H_{j,j}$ for $i,j=1,\ldots,n$ with $i\neq j$.
From these explicit formulae, the Lagrange polynomials for our interpolation set are:
\begin{itemize}
    \item For the interpolation point $\bm{0}$, 
    \begin{align}
        \ell_{\bm{0}}(\bx) = 1 - \sum_{i=1}^{n} x_i^2 + \frac{1}{2}\sum_{\substack{i,j=1 \\ j\neq i}}^{n} x_i x_j,
    \end{align}
    where we note that the Hessian $\grad^2 \ell_{\bm{0}}$ has entries $[\grad^2 \ell_{\bm{0}}]_{i,j} = 1$ if $i\neq j$ and $-2$ if $i=j$, and is a \emph{circulant matrix} \cite{Gray2005}; 
    \item For the interpolation points $\be_i$, $i=1,\ldots,n$,
    \begin{align}
        \ell_{\be_i}(\bx) = \frac{1}{2}x_i + \frac{1}{2} x_i^2 - \sum_{\substack{j=1 \\ j \neq i}}^{n} x_i x_j;
    \end{align}
    \item For the interpolation points $-\be_i$, $i=1,\ldots,n$,
    \begin{align}
        \ell_{-\be_i}(\bx) = -\frac{1}{2}x_i + \frac{1}{2} x_i^2; 
    \end{align}
    \item For the interpolation points $\be_i+\be_j$, $i,j=1\ldots,n$ with $j > i$,
    \begin{align}
        \ell_{\be_i+\be_j}(\bx) = x_i x_j.
    \end{align}
\end{itemize}

We \minorrev{now maximize} the magnitude of each Lagrange polynomial $\ell_{\by}(\bx)$ over $\bx\in B(\bm{0},1)$ individually for all $\by\in\mathcal{Y}$.

For $\ell_{\bm{0}}$, we observe that $\grad \ell_{\bm{0}}(\bm{0})=\bm{0}$ and so
\begin{align}
    |\ell_{\bm{0}}(\bx)| \leq 1 + \frac{1}{2}\|\bx\|^2 \|\grad^2 \ell_{\bm{0}}\| \leq 1 + \frac{1}{2} \|\grad^2 \ell_{\bm{0}}\|.
\end{align}
The Gershgorin circle theorem allows us to estimate\footnote{The exact value $\|\grad^2 \ell_{\bm{0}}\| = n-3$ (for $n$ sufficiently large) may be calculated using the explicit formula for eigenvalues of circulant matrices \cite{Gray2005}, but this is not needed for our estimate here.} $\|\grad^2 \ell_{\bm{0}}\| \leq n+1$, and so $\max_{\bx\in B(\bm{0},1)} |\ell_{\bm{0}}(\bx)| \leq 1 + \frac{n+1}{2}$.
For $\ell_{\be_i}$, we note that $\sum_{j\neq i} x_j^2 \leq 1 - x_i^2$ for all $\bx\in B(\bm{0},1)$, and for this constraint the term $\sum_{j\neq i} x_i x_j = x_i \sum_{j\neq i} x_j$ is maximized/minimized if all $x_j$ ($j\neq i$) are equal, $x_j = \pm \frac{1}{n-1}\sqrt{1-x_i^2}$.
So, for a given value of $x_i$,
\begin{align}
    |\ell_{\be_i}(\bx)| \leq \frac{1}{2}x_i + \frac{1}{2}x_i^2 \pm \frac{n}{n-1} x_i \sqrt{1-x_i^2},
\end{align}
(but where the last term is zero if $n=1$).
Since $x_i\in[-1,1]$ and $\frac{n}{n-1}\leq 2$ for all $n\geq 2$ we may estimate
\begin{align}
    \max_{\bx\in B(\bm{0},1)} |\ell_{\be_i}(\bx)| \leq \frac{1}{2} + \frac{1}{2} + 2\cdot 1\cdot 1 = 3.
\end{align}

For $\ell_{-\be_i}$, we note that we get the same Lagrange polynomial as the structured minimum Frobenius norm quadratic interpolation set in \minorrev{\eqref{eq_min_frob_example_lag_poly}, which gives} $\max_{\bx\in B(\bm{0},1)} |\ell_{-\be_i}(\bx)| \leq 1$.
Lastly, for $\ell_{\be_i+\be_j}$, we can observe that $|x_i|,|x_j| \leq 1$ for $\bx\in B(\bm{0},1)$, and so $\max_{\bx\in B(\bm{0},1)} |\ell_{-\be_i}(\bx)| \leq 1$.

All together, we have determined that this interpolation set has poisedness constant $\Lambda \leq \max(3, 1 + \frac{n+1}{2}) = \bigO(n)$.
\revision{This $\bigO(n)$ bound is tight: for example, we can compute $\ell_{\bm{0}}(\be/\sqrt{n}) = (5-n)/2$, which implies $\Lambda \geq |n-5|/2$.}

\revision{\section{Proof of \thmref{thm_deterministic_noise_wcc}} \label{sec_storm_deterministic_proofs}

\begin{lemma} \label{lem_det_noisy_small_delta_success}
    Suppose Assumptions~\ref{ass_smoothness_1}, \ref{ass_cauchy_decrease}, \ref{ass_noisy_oracle} and \ref{ass_model_dfo_storm_deterministic} hold.
    If $\Delta_k \leq \min(\frac{1}{\mu_c}, \frac{1}{\kappa_H})\|\bg_k\|$ and
    \begin{align}
        \|\bg_k\| \geq \frac{2\kappamf \Delta_k}{\kappa_s (1-\eta_S)} + \frac{2(\tkappamf+1)\epsilon_f}{\kappa_s (1-\eta_S) \Delta_k}, \label{eq_noisy_small_delta_success_condition}
    \end{align}
    then iteration $k$ is successful (i.e.~$\rho_k \geq \eta_S$ and $\|\bg_k\| \geq \mu_c\Delta_k$).
\end{lemma}
\begin{proof}
    That $\|\bg_k\| \geq \mu_c\Delta_k$ holds follows by assumption on $\Delta_k$, so it remains to show $\rho_k \geq \eta_S$.

    Since $\Delta_k \leq \|\bg_k\|/\kappa_H$ we have $m_k(\bx_k)-m_k(\bx_k+\bs_k) \geq \kappa_s \|\bg_k\| \Delta_k$ from \assref{ass_cauchy_decrease}.
    We then compute
    \begin{align}
        |\rho_k-1| &\leq \frac{|\tilde{f}(\bx_k) - f(\bx_k)| + |f(\bx_k)-m_k(\bx_k)| + |\tilde{f}(\bx_k+\bs_k) - f(\bx_k+\bs_k)| + |f(\bx_k+\bs_k) - m_k(\bx_k+\bs_k)|}{m_k(\bx_k)-m_k(\bx_k+\bs_k)}, \\
        &\leq \frac{2\epsilon_f + 2\kappamf \Delta_k^2 + 2\tkappamf \epsilon_f}{\kappa_s \|\bg_k\| \Delta_k}, \\
        &\leq 1-\eta_S,
    \end{align}
    where the last inequality follows from \eqref{eq_noisy_small_delta_success_condition}, and so $\rho_k \geq \eta_S$.
\end{proof}

\begin{lemma} \label{lem_delta_min_dfo_deterministic}
    Suppose Assumptions~\ref{ass_smoothness_1}, \ref{ass_cauchy_decrease}, \ref{ass_noisy_oracle} and \ref{ass_model_dfo_storm_deterministic} hold.
    If $\|\grad f(\bx_k)\| \geq \epsilon$ for all $k=0,\ldots,K-1$ with
    \begin{align}
        \epsilon > \epsilon_{\min} := \frac{2\sqrt{C_0 C_1}}{\sqrt{1 - \left(\frac{1-\gammadec}{1+\gammadec}\right)^2}},
    \end{align}
    where
    \begin{align}
        C_0 := \left(\max\left(\mu_c, \kappa_H, \frac{2\kappamf}{\kappa_s(1-\eta_S)}\right) + \kappamg\right), \quad \text{and} \quad C_1 := \left(\frac{2(\tkappamf+1)\epsilon_f}{\kappa_s (1-\eta_S)} + \tkappamg \epsilon_f\right),
    \end{align}
    and $\Delta_0$ is sufficiently large (specifically, $\Delta_0 \geq \Delta_{\min}(\epsilon)$), then
    \begin{align}
        \Delta_k \geq \Delta_{\min}(\epsilon) := \frac{1}{2}(1+\gammadec) \frac{\epsilon}{2C_0} - \frac{1}{2}(1-\gammadec)\frac{\sqrt{\epsilon^2-4C_0C_1}}{2C_0} > 0,
    \end{align}
    for all $k=0,\ldots,K$.
\end{lemma}
\begin{proof}
    By assumption on $\Delta_0$, the result holds for $k=0$. By induction, suppose that $\Delta_k \geq \Delta_{\min}(\epsilon)$ for some $k\in\{0,\ldots,K-1\}$ and to find a contradiction suppose that $\Delta_{k+1} < \Delta_{\min}(\epsilon)$.

    Since $\Delta_k \geq \Delta_{\min}(\epsilon) > \Delta_{k+1}$, iteration $k$ must have unsuccessful, and so $\Delta_{k+1}=\gammadec \Delta_k$.
    Thus, $\Delta_{\min}(\epsilon) \leq \Delta_k < \gammadec^{-1} \Delta_{\min}(\epsilon)$.

    From \assref{ass_model_dfo_storm_deterministic} we have
    \begin{align}
        \epsilon \leq \|\grad f(\bx_k)\| \leq \|\bg_k\| + \|\bg_k-\grad f(\bx_k)\| \leq \|\bg_k\| + \kappamg \Delta_k + \frac{\tkappamg \epsilon_f}{\Delta_k}.
    \end{align}
    Since iteration $k$ was unsuccessful, from \lemref{lem_det_noisy_small_delta_success} we have $\|\bg_k\| < \max(\mu_c, \kappa_H) \Delta_k$ or
    \begin{align}
        \|\bg_k\| < \frac{2\kappamf \Delta_k}{\kappa_s (1-\eta_S)} + \frac{2(\tkappamf+1)\epsilon_f}{\kappa_s (1-\eta_S) \Delta_k}.
    \end{align}
    In the first case, we get
    \begin{align}
        \epsilon < (\max(\mu_c,\kappa_H) + \kappamg)\Delta_k + \frac{\tkappamg \epsilon_f}{\Delta_k},
    \end{align}
    and in the second case we get
    \begin{align}
        \epsilon < \left(\frac{2\kappamf}{\kappa_s (1-\eta_S)} + \kappamg\right) \Delta_k + \left(\frac{2(\tkappamf+1)\epsilon_f}{\kappa_s (1-\eta_S)} + \tkappamg \epsilon_f\right) \frac{1}{\Delta_k}.
    \end{align}
    So, regardless of which case we are in, it must hold that
    \begin{align}
        \epsilon < C_0 \Delta_k + \frac{C_1}{\Delta_k}, \qquad \text{or} \qquad C_0 \Delta_k^2 - \epsilon \Delta_k + C_1 > 0. \label{eq_storm_det_quadratic_tmp1} 
    \end{align}
    We claim that this contradicts $\Delta_{\min}(\epsilon) \leq \Delta_k < \gammadec^{-1} \Delta_{\min}(\epsilon)$.

    Since $C_0,C_1,\epsilon>0$ and $\epsilon \geq \epsilon_{\min} > \sqrt{4C_0 C_1}$, \eqref{eq_storm_det_quadratic_tmp1} is a convex quadratic with two positive roots, say $\Delta^- < \Delta^+$.
    Our assumption $\epsilon > \epsilon_{\min}$ ensures that $\Delta^- < \gammadec \Delta^+$, and our choice of $\Delta_{\min}(\epsilon)$ gives  $\Delta_{\min}(\epsilon) = \frac{1}{2}(\Delta^- + \gammadec \Delta^+) \in (\Delta^-, \gammadec \Delta^+)$.
    So, \eqref{eq_storm_det_quadratic_tmp1} gives either $\Delta_k < \Delta^- < \Delta_{\min}(\epsilon)$ or $\Delta_k > \Delta^+ > \gammadec^{-1} \Delta_{\min}(\epsilon)$, which gives the desired contradiction.
\end{proof}

The proof of \thmref{thm_deterministic_noise_wcc} is then identical to that of \thmref{thm_wcc_dfo}, using \lemref{lem_delta_min_dfo_deterministic} in place of \lemref{lem_delta_min_dfo}.

\begin{remark} \label{rem_gammadec_noisy}
    Taking $\gammadec\to 1^-$ decreases $\epsilon_{\min}$ and increases $\Delta_{\min}(\epsilon)$ (and hence decreases the worst-case complexity bound, $K=\bigO(\Delta_{\min}(\epsilon)^{-2})$).
    That is, using $\gammadec\approx 1$ allows for higher-accuracy solutions and decreases the iteration complexity bound.
\end{remark}


}

\end{document}

%% file: macros_v2.tex
\usepackage{amsmath,amsthm,latexsym} 
\usepackage[fixamsmath]{mathtools}
\mathtoolsset{showmanualtags,mathic,centercolon} 
\usepackage{amssymb,amsfonts} 

\usepackage{booktabs} 
\usepackage{placeins} 
\usepackage{bm} 

\usepackage{graphicx}
\usepackage{algpseudocode} 
\usepackage{algorithm} 

\usepackage{bbm}  

\usepackage{enumitem}

\usepackage[labelfont=sl,textfont=sl]{subcaption} 
\usepackage[font={small,sl},labelsep=period]{caption} 

\usepackage[colorlinks=true, urlcolor=black, linkcolor=black, citecolor=black]{hyperref}
\usepackage{url}

\numberwithin{equation}{section}
\numberwithin{figure}{section}
\numberwithin{table}{section}
\numberwithin{algorithm}{section}

\theoremstyle{plain}
\newtheorem{theorem}{Theorem}[section]
\newtheorem{lemma}[theorem]{Lemma} 
\newtheorem{corollary}[theorem]{Corollary}
\newtheorem{proposition}[theorem]{Proposition}
\newtheorem{assumption}[theorem]{Assumption}
\newtheorem{definition}[theorem]{Definition}
\newtheorem{remark}[theorem]{Remark}
\newtheorem{example}[theorem]{Example}

\theoremstyle{remark}

\newcommand{\appref}[1]{Appendix~\ref{#1}}
\newcommand{\secref}[1]{Section~\ref{#1}}
\newcommand{\defref}[1]{Definition~\ref{#1}}
\newcommand{\thmref}[1]{Theorem~\ref{#1}}
\newcommand{\lemref}[1]{Lemma~\ref{#1}}
\newcommand{\propref}[1]{Proposition~\ref{#1}}
\newcommand{\corref}[1]{Corollary~\ref{#1}}
\renewcommand{\algref}[1]{Algorithm~\ref{#1}}
\newcommand{\assref}[1]{Assumption~\ref{#1}}
\newcommand{\figref}[1]{Figure~\ref{#1}}

\newcommand{\remref}[1]{Remark~\ref{#1}}

\newcommand{\exref}[1]{Example~\ref{#1}}

\newcommand{\R}{\mathbb{R}} 
\newcommand{\N}{\mathbb{N}} 
\newcommand{\Z}{\mathbb{Z}} 
\newcommand{\bigO}{\mathcal{O}} 
\DeclareMathOperator*{\argmin}{arg\,min} 
\DeclareMathOperator*{\argmax}{arg\,max} 
\newcommand{\defeq}{:=} 
\newcommand{\grad}{\nabla} 
\newcommand{\proj}{\operatorname{proj}} 
\newcommand{\E}{\mathbb{E}} 
\renewcommand{\P}{\mathbb{P}} 
\newcommand{\Var}{\operatorname{Var}} 

\renewcommand{\t}[1]{\widetilde{#1}} 

\algrenewcommand\algorithmicrequire{\textbf{Input:}}
\algrenewcommand\algorithmicensure{\textbf{Output:}}

\renewcommand{\b}[1]{\bm{#1}} 
\newcommand{\ba}{\b{a}}
\newcommand{\bb}{\b{b}}
\newcommand{\bc}{\b{c}}
\newcommand{\bd}{\b{d}}
\newcommand{\be}{\b{e}}
\newcommand{\bg}{\b{g}}
\newcommand{\bem}{\b{m}}  
\newcommand{\bn}{\b{n}}

\newcommand{\br}{\b{r}}
\newcommand{\bs}{\b{s}}
\newcommand{\bt}{\b{t}}
\newcommand{\bu}{\b{u}}
\newcommand{\bv}{\b{v}}

\newcommand{\bx}{\b{x}}
\newcommand{\by}{\b{y}}

\newcommand{\bphi}{\b{\phi}}
\newcommand{\blambda}{\b{\lambda}}
\newcommand{\bvarphi}{\b{\varphi}}

\newcommand{\bA}{\b{A}}
\newcommand{\bF}{\b{F}}
\newcommand{\bH}{\b{H}}
\newcommand{\bI}{\b{I}}
\newcommand{\bJ}{\b{J}}

\newcommand{\bM}{\b{M}}
\newcommand{\bP}{\b{P}}
\newcommand{\bQ}{\b{Q}}

\newcommand{\C}{\mathcal{C}}

\newcommand{\gammainc}{\gamma_{\textnormal{inc}}}
\newcommand{\gammadec}{\gamma_{\textnormal{dec}}}
\newcommand{\flow}{f_{\textnormal{low}}}

\newcommand{\kappam}{\kappa_{\textnormal{m}}}
\newcommand{\tkappam}{\tilde{\kappa}_{\textnormal{m}}}
\newcommand{\kappamf}{\kappa_{\textnormal{mf}}}
\newcommand{\kappamg}{\kappa_{\textnormal{mg}}}
\newcommand{\kappamh}{\kappa_{\textnormal{mh}}}
\newcommand{\tkappamf}{\tilde{\kappa}_{\textnormal{mf}}}
\newcommand{\tkappamg}{\tilde{\kappa}_{\textnormal{mg}}}

\newcommand{\upper}{\operatorname{upper}}
\newcommand{\Hprev}{\bm{H}_{\textnormal{prev}}}

\newcommand{\ind}{\mathbbm{1}}  
